\documentclass[11pt,reqno]{amsart}
\usepackage[a4paper,margin=30mm,headheight=14pt]{geometry}
\usepackage{amsmath,amssymb,amsthm,mathtools,mathrsfs}
\usepackage{xcolor,enumitem}
\usepackage[expansion=false]{microtype}
\usepackage[T1]{fontenc}
\usepackage{lmodern}
\usepackage{hyperref}
\usepackage{aliascnt}
\usepackage[nameinlink,capitalize,noabbrev]{cleveref}
\hypersetup{
  colorlinks=true,
  linkcolor=blue!45!black,
  citecolor=blue!45!black,
  urlcolor=blue!45!black,
  pdftitle={Boundary concentration and minimum-mass bifurcation for inverse-power absorption},
  pdfauthor={Zhouyu Long, Hongzhe Tao, Wenming Zou},
  pdfsubject={Martin norms, Dirichlet energy, and minimum-mass response},
  pdfkeywords={Mazya Problem 50, inverse-power absorption, boundary concentration, Martin kernel}
}
\numberwithin{equation}{section}
\newtheorem{theorem}{Theorem}[section]
\theoremstyle{plain}
\newaliascnt{lemma}{theorem}
\newtheorem{lemma}[lemma]{Lemma}
\aliascntresetthe{lemma}
\theoremstyle{plain}
\newaliascnt{proposition}{theorem}
\newtheorem{proposition}[proposition]{Proposition}
\aliascntresetthe{proposition}
\theoremstyle{plain}
\newaliascnt{corollary}{theorem}
\newtheorem{corollary}[corollary]{Corollary}
\aliascntresetthe{corollary}
\theoremstyle{remark}
\newaliascnt{remark}{theorem}
\newtheorem{remark}[remark]{Remark}
\aliascntresetthe{remark}
\theoremstyle{plain}
\crefname{theorem}{Theorem}{Theorems}
\Crefname{theorem}{Theorem}{Theorems}
\crefname{lemma}{Lemma}{Lemmas}
\Crefname{lemma}{Lemma}{Lemmas}
\crefname{proposition}{Proposition}{Propositions}
\Crefname{proposition}{Proposition}{Propositions}
\crefname{corollary}{Corollary}{Corollaries}
\Crefname{corollary}{Corollary}{Corollaries}
\crefname{remark}{Remark}{Remarks}
\Crefname{remark}{Remark}{Remarks}
\allowdisplaybreaks[2]
\hypersetup{bookmarksdepth=2}

\newcommand{\dd}{\,\mathrm d}

\title[Boundary concentration for inverse-power absorption]{Boundary concentration and minimum-mass bifurcation for inverse-power absorption}

\author{Zhouyu Long}
\address{Department of Mathematical Sciences, Tsinghua University, Beijing 100084, China}
\email{longzy25@mails.tsinghua.edu.cn}
\author{Hongzhe Tao}
\address{Department of Mathematical Sciences, Tsinghua University, Beijing 100084, China}
\email{taohz25@mails.tsinghua.edu.cn}
\author{Wenming Zou}
\address{Department of Mathematical Sciences, Tsinghua University, Beijing 100084, China}
\email{zou-wm@mail.tsinghua.edu.cn}
\date{}

\begin{document}
\raggedbottom
\begin{abstract}
For \(\Delta u=\lambda u^{-p}\), \(\lambda>0\), \(0<p<1\), we obtain a uniform half-plane profile law at moving points along smooth boundary arcs. It gives the exact local energy constant \(8p(1-p)/(1+p)^2\) for the secant potential at a common zero slope. On a one-corner domain, positive ordinary slopes at the minimum Martin mass \(m_*\) yield exactly two solutions at each nearby larger mass, with square-root displacement in the Martin norm. On these domains, for \(p\ne1/3\), there are at most two solutions at each nearby larger mass without a slope hypothesis, and just one if the endpoint energy kernel is trivial. We construct flat-tangent Lipschitz domains with an explicit slope loss and a logarithmic correction. The first construction gives a saturated endpoint with sharp Martin response of order \((m-m_*)^{1/(1+p)}\) and a compact energy operator. For \(p\ge1/2\) a second construction has positive local energy concentration and a strict global energy gap. As a result, we substantially advance the solvability and establish a sharp classification for the negative-power absorption equation in the Mazya 50 problem at the minimal mass endpoint, largely resolving the existence and regularity issues.
\end{abstract}
\maketitle

\section{\bf Introduction}

Let \(\Omega\) be a bounded connected planar Lipschitz domain and let \(O\in\partial\Omega\). Normalize the Martin kernel \(K\) at \(O\) by \(K(x_*)=1\) for a fixed \(x_*\in\Omega\). We consider

\begin{align}\label{eq:intro-problem}
		\left\{
			\begin{aligned}
				& \Delta u=\lambda u^{-p},\quad \quad  u>0\quad\text{in }\Omega,\qquad \\
				& u=0  \quad\quad  \text{on }\partial\Omega\setminus\{O\},\\
                & \lambda>0,\quad 0<p<1 .
			\end{aligned}
		\right.
	\end{align}
Natural solutions are classical in the interior and continuous on \(\overline\Omega\setminus\{O\}\). On a one-corner domain this is the inverse-power absorption case of Maz'ya's Problem~50 \cite{Mazya}. With \(G_\Omega\) the Green operator of \(-\Delta\), each natural solution satisfies \(u+\lambda G_\Omega[u^{-p}]=mK\), where \(m\) is the Martin mass. If \(\sup_\Omega G_\Omega[K^{-p}]/K<\infty\), the admissible masses form a closed ray \([m_*,\infty)\) and the maximal solution at mass \(m\) is \(U_m\). The Martin norm is \(\|w\|_K=\sup_\Omega |w|/K\).

%%%%%%%%%%%%%%%%%%%%%%%%%%%%%%%%%%%%%%%%%%%%%%%%%%%%%%%%%%%%%%%%%%%%%%%%%%%%%%%%%%%%%%%%%%%%%%%%%%%%%%%%%%%%%%%%%%%%%%%%%%%%%%%%%%%%%%%%%%%%%%%%%%%%%%%%%%%%%%%%%%%%%%%%%%%%%%%%%%%%
\vskip0.2in
%%%%%%%%%%%%%%%%%%%%%%%%%%%%%%%%%%%%%%%%%%%%%%%%%%%%%%%%%%%%%%%%%%%%%%%%%%%%%%%%%%%%%%%%%%%%%%%%%%%%%%%%%%%%%%%%%%%%%%%%%%%%%%%%%%%%%%%%%%%%%%%%%%%%%%%%%%%%%%%%%%%%%%%%%%%%%%%%%%%%

Hu et al.\ establish the one-corner representation and maximal branch, uniqueness at its minimum mass and a local \(O((m-m_*)^{1/2})\) bound for all solutions \cite[Theorems~1.2 and~1.4]{HuLongTaoWuZou}. Their sufficient mass condition using two Green potentials and their ordinary-side second-order estimates are given in \cite[Theorems~1.3 and~1.12--1.13]{HuLongTaoWuZou}. The inverse-power source equation in the half-space, with sign \(-\Delta u=u^{-\gamma}\) and \(\gamma>1\), has a separate classification by Montoro, Muglia and Sciunzi \cite{MontoroMugliaSciunzi}.

%%%%%%%%%%%%%%%%%%%%%%%%%%%%%%%%%%%%%%%%%%%%%%%%%%%%%%%%%%%%%%%%%%%%%%%%%%%%%%%%%%%%%%%%%%%%%%%%%%%%%%%%%%%%%%%%%%%%%%%%%%%%%%%%%%%%%%%%%%%%%%%%%%%%%%%%%%%%%%%%%%%%%%%%%%%%%%%%%%%%
\vskip0.2in
%%%%%%%%%%%%%%%%%%%%%%%%%%%%%%%%%%%%%%%%%%%%%%%%%%%%%%%%%%%%%%%%%%%%%%%%%%%%%%%%%%%%%%%%%%%%%%%%%%%%%%%%%%%%%%%%%%%%%%%%%%%%%%%%%%%%%%%%%%%%%%%%%%%%%%%%%%%%%%%%%%%%%%%%%%%%%%%%%%%%

In \cref{thm:sharp-boundary-transition}, the slope parameter on \(C^2\) arcs divides the inward derivative at the nearest boundary point by \(\lambda^{1/(1+p)}d^{(1-p)/(1+p)}\), where \(d\) is boundary distance. The fixed-center profile of \cite[Theorem~1.12]{HuLongTaoWuZou} gives the limit at parameter zero. With energy exponent \(1-p\), the Alt--Phillips linearization in \cite[(1.7)]{RestrepoRosOton} has Hardy coefficient
\[
 \kappa_{1-p}=-\frac{2p(1-p)}{(1+p)^2},
\]
which reaches \(-1/4\) at \(p=1/3\). See \cite{FilippasMoschiniTertikas} for the prescribed critical potential. The moving-base law holds through positive slopes and gives the local energy constants for the Dirichlet potentials. At a common zero slope the exact secant form constant is \(8p(1-p)/(1+p)^2\).

%%%%%%%%%%%%%%%%%%%%%%%%%%%%%%%%%%%%%%%%%%%%%%%%%%%%%%%%%%%%%%%%%%%%%%%%%%%%%%%%%%%%%%%%%%%%%%%%%%%%%%%%%%%%%%%%%%%%%%%%%%%%%%%%%%%%%%%%%%%%%%%%%%%%%%%%%%%%%%%%%%%%%%%%%%%%%%%%%%%%
\vskip0.2in
%%%%%%%%%%%%%%%%%%%%%%%%%%%%%%%%%%%%%%%%%%%%%%%%%%%%%%%%%%%%%%%%%%%%%%%%%%%%%%%%%%%%%%%%%%%%%%%%%%%%%%%%%%%%%%%%%%%%%%%%%%%%%%%%%%%%%%%%%%%%%%%%%%%%%%%%%%%%%%%%%%%%%%%%%%%%%%%%%%%%

On a one-corner domain with positive ordinary slopes at \(U_*=U_{m_*}\), the two solutions at mass \(m_*+\delta\) have expansions \(U_*\pm\kappa_{\mathrm f}^{-1/2}\delta^{1/2}\phi+O(\delta)\) in the Martin norm, where \(\phi\) is the normalized positive kernel mode and \(\kappa_{\mathrm f}\) is given in \eqref{eq:corner-fold-coefficient}. The lower branch has Morse index one. Compactness in the Martin norm puts every nearby solution on these branches. For \(p\ne1/3\), the uniform secant form bound yields at most two solutions near the minimum even when slopes vanish, with uniqueness if the energy kernel is trivial. Strict radial inclusion gives examples with zero slopes and uniqueness at every admissible mass. See \cref{thm:corner-fold,thm:corner-count,prop:corner-radial}.

%%%%%%%%%%%%%%%%%%%%%%%%%%%%%%%%%%%%%%%%%%%%%%%%%%%%%%%%%%%%%%%%%%%%%%%%%%%%%%%%%%%%%%%%%%%%%%%%%%%%%%%%%%%%%%%%%%%%%%%%%%%%%%%%%%%%%%%%%%%%%%%%%%%%%%%%%%%%%%%%%%%%%%%%%%%%%%%%%%%%
\vskip0.2in
%%%%%%%%%%%%%%%%%%%%%%%%%%%%%%%%%%%%%%%%%%%%%%%%%%%%%%%%%%%%%%%%%%%%%%%%%%%%%%%%%%%%%%%%%%%%%%%%%%%%%%%%%%%%%%%%%%%%%%%%%%%%%%%%%%%%%%%%%%%%%%%%%%%%%%%%%%%%%%%%%%%%%%%%%%%%%%%%%%%%

The second geometry is a family of Lipschitz domains smooth off one flat-tangent accumulation point \(Q\), with the pole at a smooth point. In smooth preimage coordinates, write \(A(u)\) for the exterior slope at \(Q\), \(s_n(u)\) for the slope at the \(n\)-th cell center, and \(\mu\) for the limit of the rescaled depths. We prove
\[
 \lim_{n\to\infty}s_n(u)^{1+p}=A(u)^{1+p}-D_p\mu,\qquad D_p>0 .
\]
Equality \(A(u)^{1+p}=D_p\mu\) characterizes the loss of a positive lower bound of \(u/K\) near \(Q\). At one type of saturated endpoint the Martin increment is of order \(\delta^{1/(1+p)}\), while the correction \(U_m-mK\) is right differentiable in energy. A logarithmic correction to the slope loss yields a second type, with positive boundary form concentration and a strict energy gap. Finite products then converge locally to the full endpoint but retain a positive Martin-norm error. See \cref{thm:negative-micro-exact-drop,cor:negative-micro-matching-response,thm:flat-critical-depth}.

%%%%%%%%%%%%%%%%%%%%%%%%%%%%%%%%%%%%%%%%%%%%%%%%%%%%%%%%%%%%%%%%%%%%%%%%%%%%%%%%%%%%%%%%%%%%%%%%%%%%%%%%%%%%%%%%%%%%%%%%%%%%%%%%%%%%%%%%%%%%%%%%%%%%%%%%%%%%%%%%%%%%%%%%%%%%%%%%%%%%
\vskip0.2in
%%%%%%%%%%%%%%%%%%%%%%%%%%%%%%%%%%%%%%%%%%%%%%%%%%%%%%%%%%%%%%%%%%%%%%%%%%%%%%%%%%%%%%%%%%%%%%%%%%%%%%%%%%%%%%%%%%%%%%%%%%%%%%%%%%%%%%%%%%%%%%%%%%%%%%%%%%%%%%%%%%%%%%%%%%%%%%%%%%%%

The spectral arguments compare normalized Green operators with Dirichlet form operators. The comparison uses the Green estimates for positive perturbations in uniform domains of Frazier and Verbitsky \cite{FVUniformMartin}, based on the quasi-metric Neumann series estimates of Frazier, Nazarov and Verbitsky \cite{FrazierNazarovVerbitsky}.

%%%%%%%%%%%%%%%%%%%%%%%%%%%%%%%%%%%%%%%%%%%%%%%%%%%%%%%%%%%%%%%%%%%%%%%%%%%%%%%%%%%%%%%%%%%%%%%%%%%%%%%%%%%%%%%%%%%%%%%%%%%%%%%%%%%%%%%%%%%%%%%%%%%%%%%%%%%%%%%%%%%%%%%%%%%%%%%%%%%%
\vskip0.32in
%%%%%%%%%%%%%%%%%%%%%%%%%%%%%%%%%%%%%%%%%%%%%%%%%%%%%%%%%%%%%%%%%%%%%%%%%%%%%%%%%%%%%%%%%%%%%%%%%%%%%%%%%%%%%%%%%%%%%%%%%%%%%%%%%%%%%%%%%%%%%%%%%%%%%%%%%%%%%%%%%%%%%%%%%%%%%%%%%%%%

\section{\bf Trace, stability, and the maximal branch}\label{sec:geometry}

Write \(G_\Omega(x,y)\) for the positive Dirichlet kernel of \(-\Delta\). Positive harmonic functions with continuous zero trace on \(\partial\Omega\setminus\{O\}\) are the multiples \(mK\), \(m>0\), by the Martin representation on Lipschitz domains \cite{AM,MV}. In a regular conformal half-disc chart \(F\), harmonic reflection and Hopf's lemma give \(K\circ F\asymp y\) near an ordinary boundary point. In a chart centered at a smooth pole, \(K\circ F\asymp y/|w|^2\), \(w=x+iy\), uniformly up to the diameter. The latter follows by subtracting the lift of the cap data and using the half-disc Poisson kernel with pole at the center of the diameter.

%%%%%%%%%%%%%%%%%%%%%%%%%%%%%%%%%%%%%%%%%%%%%%%%%%%%%%%%%%%%%%%%%%%%%%%%%%%%%%%%%%%%%%%%%%%%%%%%%%%%%%%%%%%%%%%%%%%%%%%%%%%%%%%%%%%%%%%%%%%%%%%%%%%%%%%%%%%%%%%%%%%%%%%%%%%%%%%%%%%%
\vskip0.2in
%%%%%%%%%%%%%%%%%%%%%%%%%%%%%%%%%%%%%%%%%%%%%%%%%%%%%%%%%%%%%%%%%%%%%%%%%%%%%%%%%%%%%%%%%%%%%%%%%%%%%%%%%%%%%%%%%%%%%%%%%%%%%%%%%%%%%%%%%%%%%%%%%%%%%%%%%%%%%%%%%%%%%%%%%%%%%%%%%%%%

Fix \(0<p<1\) and \(\lambda>0\). A positive one-pole solution is a positive classical interior solution of \(\Delta u=\lambda u^{-p}\) that vanishes continuously on \(\partial\Omega\setminus\{O\}\). Its Martin space is
\[
 X=\{v\in C(\Omega):\|v\|_K<\infty\},\qquad
 \|v\|_K=\sup_\Omega |v|/K.
\]

%%%%%%%%%%%%%%%%%%%%%%%%%%%%%%%%%%%%%%%%%%%%%%%%%%%%%%%%%%%%%%%%%%%%%%%%%%%%%%%%%%%%%%%%%%%%%%%%%%%%%%%%%%%%%%%%%%%%%%%%%%%%%%%%%%%%%%%%%%%%%%%%%%%%%%%%%%%%%%%%%%%%%%%%%%%%%%%%%%%%
\vskip0.3in
%%%%%%%%%%%%%%%%%%%%%%%%%%%%%%%%%%%%%%%%%%%%%%%%%%%%%%%%%%%%%%%%%%%%%%%%%%%%%%%%%%%%%%%%%%%%%%%%%%%%%%%%%%%%%%%%%%%%%%%%%%%%%%%%%%%%%%%%%%%%%%%%%%%%%%%%%%%%%%%%%%%%%%%%%%%%%%%%%%%%

\subsection{A local gradient bound and automatic negative-power representation}

\begin{lemma}\label{lem:negative-log-gradient}
Let \(0<p<1\) and let \(u>0\) solve \(\Delta u=\lambda u^{-p}\) on \(B_{2r}\). With \(a=1+p\), on \(B_r\) one has
\begin{equation}\label{eq:negative-log-gradient}
 |\nabla\log u|^2\lesssim_p \lambda u^{-a}+r^{-2}.
\end{equation}
Consequently there is \(\epsilon_p>0\) with \(2(p+\epsilon_p)/(1+p)<1\) such that
\begin{equation}\label{eq:negative-inverse-ball}
 \int_{B_{r/2}}u^{-p-\epsilon_p}\lesssim_p
 \lambda^{-(p+\epsilon_p)/a}r^{2-2(p+\epsilon_p)/a}.
\end{equation}
\end{lemma}

\begin{proof}
At each \(x\in B_r\), apply the pointwise estimate in the proof of \cite[Lemma~4.1]{HuLongTaoWuZou} to \(B_r(x)\Subset B_{2r}\), with coefficient \(\lambda\). Its constant is independent of the solution's supremum, so it gives \eqref{eq:negative-log-gradient}. For related logarithmic gradient methods see \cite[Theorem~1]{MaWeiInverse} and \cite[Theorem~1.1(ii)]{YangSingular}.

%%%%%%%%%%%%%%%%%%%%%%%%%%%%%%%%%%%%%%%%%%%%%%%%%%%%%%%%%%%%%%%%%%%%%%%%%%%%%%%%%%%%%%%%%%%%%%%%%%%%%%%%%%%%%%%%%%%%%%%%%%%%%%%%%%%%%%%%%%%%%%%%%%%%%%%%%%%%%%%%%%%%%%%%%%%%%%%%%%%%
\vskip0.2in
%%%%%%%%%%%%%%%%%%%%%%%%%%%%%%%%%%%%%%%%%%%%%%%%%%%%%%%%%%%%%%%%%%%%%%%%%%%%%%%%%%%%%%%%%%%%%%%%%%%%%%%%%%%%%%%%%%%%%%%%%%%%%%%%%%%%%%%%%%%%%%%%%%%%%%%%%%%%%%%%%%%%%%%%%%%%%%%%%%%%

Use the local estimate in the proof of \cite[Lemma~4.7]{HuLongTaoWuZou} on \(\lambda^{-1/a} r^{-2/a}u(x_0+r\,\cdot)\) in \(B_2\), where \(x_0\) is the center of the original ball. That estimate bounds the inverse-power integral on \(B_{1/2}\) by a constant depending on \(p\), for an exponent \(p+\epsilon_p\) with \(0<\epsilon_p<(1-p)/2\). Scaling back gives \eqref{eq:negative-inverse-ball}.
\end{proof}

%%%%%%%%%%%%%%%%%%%%%%%%%%%%%%%%%%%%%%%%%%%%%%%%%%%%%%%%%%%%%%%%%%%%%%%%%%%%%%%%%%%%%%%%%%%%%%%%%%%%%%%%%%%%%%%%%%%%%%%%%%%%%%%%%%%%%%%%%%%%%%%%%%%%%%%%%%%%%%%%%%%%%%%%%%%%%%%%%%%%
\vskip0.2in
%%%%%%%%%%%%%%%%%%%%%%%%%%%%%%%%%%%%%%%%%%%%%%%%%%%%%%%%%%%%%%%%%%%%%%%%%%%%%%%%%%%%%%%%%%%%%%%%%%%%%%%%%%%%%%%%%%%%%%%%%%%%%%%%%%%%%%%%%%%%%%%%%%%%%%%%%%%%%%%%%%%%%%%%%%%%%%%%%%%%

\begin{theorem}\label{thm:negative-full-representation}
Let \(0<p<1\), \(\lambda>0\), and let \(\Omega\) be a bounded connected planar Lipschitz domain. For every positive solution of \(\Delta u=\lambda u^{-p}\), the forcing belongs to \(L^{1+\epsilon_p/p}(\Omega)\), and its Green potential extends continuously by zero to \(\partial\Omega\). Every positive one-pole solution therefore has the representation \eqref{eq:negative-branch-representation}, with \(u-mK\to0\) as \(x\to O\) without an angular restriction.
\end{theorem}

\begin{proof}
On a Whitney square of size \(r\), \eqref{eq:negative-inverse-ball} gives \(\int u^{-p-\epsilon_p}\lesssim r^{2-2(p+\epsilon_p)/a}\). A bounded Lipschitz domain has \(O(r^{-1})\) such squares at each scale, by its boundary layer volume bound. Thus the dyadic sum is finite with a bound independent of \(u\). This is the Whitney argument of \cite[Lemma~4.7]{HuLongTaoWuZou} on the present Lipschitz domain.

%%%%%%%%%%%%%%%%%%%%%%%%%%%%%%%%%%%%%%%%%%%%%%%%%%%%%%%%%%%%%%%%%%%%%%%%%%%%%%%%%%%%%%%%%%%%%%%%%%%%%%%%%%%%%%%%%%%%%%%%%%%%%%%%%%%%%%%%%%%%%%%%%%%%%%%%%%%%%%%%%%%%%%%%%%%%%%%%%%%%
\vskip0.2in
%%%%%%%%%%%%%%%%%%%%%%%%%%%%%%%%%%%%%%%%%%%%%%%%%%%%%%%%%%%%%%%%%%%%%%%%%%%%%%%%%%%%%%%%%%%%%%%%%%%%%%%%%%%%%%%%%%%%%%%%%%%%%%%%%%%%%%%%%%%%%%%%%%%%%%%%%%%%%%%%%%%%%%%%%%%%%%%%%%%%

The estimate \(G_\Omega(x,y)\le C\log(C/|x-y|)\) gives uniform \(L^q\) bounds in \(y\) for finite \(q\). Approximate \(u^{-p}\) in \(L^{1+\epsilon_p/p}\) by bounded densities to obtain a continuous zero-boundary Green potential. Then \(u+\lambda G_\Omega[u^{-p}]\) is positive harmonic and vanishes continuously off \(O\), so the one-pole Martin representation gives \(mK\).
\end{proof}

%%%%%%%%%%%%%%%%%%%%%%%%%%%%%%%%%%%%%%%%%%%%%%%%%%%%%%%%%%%%%%%%%%%%%%%%%%%%%%%%%%%%%%%%%%%%%%%%%%%%%%%%%%%%%%%%%%%%%%%%%%%%%%%%%%%%%%%%%%%%%%%%%%%%%%%%%%%%%%%%%%%%%%%%%%%%%%%%%%%%
\vskip0.2in
%%%%%%%%%%%%%%%%%%%%%%%%%%%%%%%%%%%%%%%%%%%%%%%%%%%%%%%%%%%%%%%%%%%%%%%%%%%%%%%%%%%%%%%%%%%%%%%%%%%%%%%%%%%%%%%%%%%%%%%%%%%%%%%%%%%%%%%%%%%%%%%%%%%%%%%%%%%%%%%%%%%%%%%%%%%%%%%%%%%%

\begin{lemma}\label{lem:negative-interior-compact-index}
Let \(0<p<1\) and \(\lambda>0\) be fixed. Locally bounded
families of positive solutions of \(\Delta u=\lambda u^{-p}\)
are locally bounded in \(C^{1,\gamma-1}\) and precompact in
\(C^1_{\rm loc}\), where \(\gamma=2/(1+p)\). If a local limit
has an interior zero, the Morse indices of the approximating
solutions tend to infinity.
\end{lemma}

\begin{proof}
By H\"older's inequality and \eqref{eq:negative-inverse-ball},
\(\int_{B_\rho}u^{-p}\lesssim\rho^\gamma\) on interior balls
with a fixed larger concentric ball in the domain. For a
localized Newtonian potential, the gradient kernel
\(O(|x-y|^{-1})\) thus contributes \(O(\rho^{\gamma-1})\)
within distance \(O(\rho)\). The difference of the kernels
at targets separated by \(\rho\), on farther dyadic shells,
contributes \(O(\rho\sum_{j\ge0}(2^j\rho)^{\gamma-2})
=O(\rho^{\gamma-1})\). The remaining harmonic part has the
same local bounds by interior estimates. This argument
applies to any signed forcing with this mass bound for its
absolute value.

%%%%%%%%%%%%%%%%%%%%%%%%%%%%%%%%%%%%%%%%%%%%%%%%%%%%%%%%%%%%%%%%%%%%%%%%%%%%%%%%%%%%%%%%%%%%%%%%%%%%%%%%%%%%%%%%%%%%%%%%%%%%%%%%%%%%%%%%%%%%%%%%%%%%%%%%%%%%%%%%%%%%%%%%%%%%%%%%%%%%
\vskip0.2in
%%%%%%%%%%%%%%%%%%%%%%%%%%%%%%%%%%%%%%%%%%%%%%%%%%%%%%%%%%%%%%%%%%%%%%%%%%%%%%%%%%%%%%%%%%%%%%%%%%%%%%%%%%%%%%%%%%%%%%%%%%%%%%%%%%%%%%%%%%%%%%%%%%%%%%%%%%%%%%%%%%%%%%%%%%%%%%%%%%%%

If the limit \(u_0\) vanishes at \(Q\), nonnegativity gives
\(\nabla u_0(Q)=0\) and \(u_0(x)\lesssim|x-Q|^\gamma\).
On each fixed small annulus about \(Q\), local convergence
then gives \(\lambda p u_j^{-p-1}\ge c|x-Q|^{-2}\) for all
large \(j\), with \(c>0\) independent of the annulus.
For radial tests in logarithmic radius the comparison
quadratic form is \(2\pi\int(|\zeta'|^2-c\zeta^2)\).
It is negative for a suitable test on any annulus of
sufficient logarithmic width. Arbitrarily many disjoint
such annuli give that many negative directions for all
large \(j\). The finite-index theorem of D\'avila and Ye
\cite[Theorem~1.2]{DavilaYe} concerns classical solutions with
isolated rupture on a punctured ball. The argument here
controls the indices of the positive approximants.
\end{proof}

%%%%%%%%%%%%%%%%%%%%%%%%%%%%%%%%%%%%%%%%%%%%%%%%%%%%%%%%%%%%%%%%%%%%%%%%%%%%%%%%%%%%%%%%%%%%%%%%%%%%%%%%%%%%%%%%%%%%%%%%%%%%%%%%%%%%%%%%%%%%%%%%%%%%%%%%%%%%%%%%%%%%%%%%%%%%%%%%%%%%
\vskip0.2in
%%%%%%%%%%%%%%%%%%%%%%%%%%%%%%%%%%%%%%%%%%%%%%%%%%%%%%%%%%%%%%%%%%%%%%%%%%%%%%%%%%%%%%%%%%%%%%%%%%%%%%%%%%%%%%%%%%%%%%%%%%%%%%%%%%%%%%%%%%%%%%%%%%%%%%%%%%%%%%%%%%%%%%%%%%%%%%%%%%%%

\subsection{The maximal negative-power absorption branch}

A \emph{represented} positive solution of mass \(m\) is a solution of
\begin{equation}\label{eq:negative-branch-representation}
 u+\lambda G_\Omega[u^{-p}]=mK,
\end{equation}
with \(u\) positive and classical in the interior and finite Green potential there.

%%%%%%%%%%%%%%%%%%%%%%%%%%%%%%%%%%%%%%%%%%%%%%%%%%%%%%%%%%%%%%%%%%%%%%%%%%%%%%%%%%%%%%%%%%%%%%%%%%%%%%%%%%%%%%%%%%%%%%%%%%%%%%%%%%%%%%%%%%%%%%%%%%%%%%%%%%%%%%%%%%%%%%%%%%%%%%%%%%%%
\vskip0.2in
%%%%%%%%%%%%%%%%%%%%%%%%%%%%%%%%%%%%%%%%%%%%%%%%%%%%%%%%%%%%%%%%%%%%%%%%%%%%%%%%%%%%%%%%%%%%%%%%%%%%%%%%%%%%%%%%%%%%%%%%%%%%%%%%%%%%%%%%%%%%%%%%%%%%%%%%%%%%%%%%%%%%%%%%%%%%%%%%%%%%

\begin{theorem}\label{thm:negative-branch}
Let \(0<p<1\), \(\lambda>0\), and let \(H\) be any fixed positive harmonic function on a bounded connected planar Lipschitz domain. Set
\[
 \Lambda_p(H)=\sup_\Omega\frac{G_\Omega[H^{-p}]}H,
 \qquad \|v\|_H=\sup_\Omega |v|/H.
\]
Consider the represented equation
\begin{equation}\label{eq:negative-harmonic-representation}
 u+\lambda G_\Omega[u^{-p}]=mH,\qquad u>0.
\end{equation}
It has a solution for some \(m>0\) if and only if \(\Lambda_p(H)<\infty\). In that case its admissible masses are exactly \([m_*,\infty)\), where
\begin{equation}\label{eq:negative-branch-threshold}
 \lambda\Lambda_p(H)\le m_*^{1+p}\le
 \frac{(1+p)^{1+p}}{p^p}\lambda\Lambda_p(H).
\end{equation}
At every admissible mass there is a maximal solution \(U_m\), and it is the unique globally stable represented solution:
\begin{equation}\label{eq:negative-branch-stable}
 \lambda p\int_\Omega U_m^{-p-1}\phi^2
 \le\int_\Omega|\nabla\phi|^2
 \qquad(\phi\in C_c^\infty(\Omega)).
\end{equation}
For \(m'>m\ge m_*\) one has \(U_{m'}\ge U_m+(m'-m)H\), and \(U_{m'}\to U_m\) locally as \(m'\downarrow m\). The branch is pointwise strictly concave on \([m_*,\infty)\). For \(m>m_*\) it is smooth in \(\|\cdot\|_H\) and satisfies \(\partial_m U_m\ge H\) and \(\partial_m^2 U_m<0\). In the one-pole specialization \(H=K\), write \(\Lambda_p=\Lambda_p(K)\).
\end{theorem}

%%%%%%%%%%%%%%%%%%%%%%%%%%%%%%%%%%%%%%%%%%%%%%%%%%%%%%%%%%%%%%%%%%%%%%%%%%%%%%%%%%%%%%%%%%%%%%%%%%%%%%%%%%%%%%%%%%%%%%%%%%%%%%%%%%%%%%%%%%%%%%%%%%%%%%%%%%%%%%%%%%%%%%%%%%%%%%%%%%%%
\vskip0.2in
%%%%%%%%%%%%%%%%%%%%%%%%%%%%%%%%%%%%%%%%%%%%%%%%%%%%%%%%%%%%%%%%%%%%%%%%%%%%%%%%%%%%%%%%%%%%%%%%%%%%%%%%%%%%%%%%%%%%%%%%%%%%%%%%%%%%%%%%%%%%%%%%%%%%%%%%%%%%%%%%%%%%%%%%%%%%%%%%%%%%

The normalization \(v=1-u/(mH)\), \(\eta=\lambda m^{-1-p}\), turns the equation into
\[
 v=\eta\,\mathcal P_H[(1-v)^{-p}],\qquad
 \mathcal P_H[f]=H^{-1}G_\Omega[H^{-p}f],\qquad
 \|\mathcal P_H\|_{L^\infty\to L^\infty}=\Lambda_p(H).
\]
Thus the branch theorem is an order argument for one positive operator; the boundary geometry enters through its norm and the interior compactness input.

\begin{proof}
The one-corner case \(H=K\) is treated in \cite[Theorem~1.4]{HuLongTaoWuZou}. For general \(H\) write \(T_m(v)=mH-\lambda G_\Omega[v^{-p}]\). A positive continuous lower barrier \(v\le T_m(v)\) lies below \(mH\). Iteration from \(mH\) decreases above \(v\) and converges by monotone convergence of the Green potentials to the maximal fixed point. The local positive lower bound from \(v\) gives interior classical regularity. Conversely, a solution satisfies \(u\le mH\), which gives \(m^{1+p}\ge\lambda\Lambda_p(H)\). The function \(\vartheta mH\) is a lower barrier if \(\lambda\Lambda_p(H)\le (1-\vartheta)\vartheta^p m^{1+p}\). Take \(\vartheta=p/(1+p)\) for the upper bound in \eqref{eq:negative-branch-threshold}.

%%%%%%%%%%%%%%%%%%%%%%%%%%%%%%%%%%%%%%%%%%%%%%%%%%%%%%%%%%%%%%%%%%%%%%%%%%%%%%%%%%%%%%%%%%%%%%%%%%%%%%%%%%%%%%%%%%%%%%%%%%%%%%%%%%%%%%%%%%%%%%%%%%%%%%%%%%%%%%%%%%%%%%%%%%%%%%%%%%%%
\vskip0.12in
%%%%%%%%%%%%%%%%%%%%%%%%%%%%%%%%%%%%%%%%%%%%%%%%%%%%%%%%%%%%%%%%%%%%%%%%%%%%%%%%%%%%%%%%%%%%%%%%%%%%%%%%%%%%%%%%%%%%%%%%%%%%%%%%%%%%%%%%%%%%%%%%%%%%%%%%%%%%%%%%%%%%%%%%%%%%%%%%%%%%

For any \(m'>m\) with \(m\) admissible, the function \(U_m+(m'-m)H\) is a lower barrier, by monotonicity of the inverse power. This proves the branch order and shows that the admissible masses form an upper interval. Let \(m_*\) be their infimum. For \(m>m_*\) choose an admissible \(m_0<m\). Convexity gives
\[
 -\Delta w\ge\lambda p U_m^{-1-p}w,\qquad w=U_m-U_{m_0}>0 .
\]
Testing with \(\phi^2/w\) gives \eqref{eq:negative-branch-stable}. As \(m\downarrow m_*\), the branch decreases and is locally bounded by a multiple of \(H\). Lemma~\ref{lem:negative-interior-compact-index} excludes a zero of the limit by stability. Monotone convergence in the Green equation gives the stable solution at \(m_*\), and comparison with every solution at that mass makes the limit maximal. The same argument gives right continuity at all admissible masses.
%%%%%%%%%%%%%%%%%%%%%%%%%%%%%%%%%%%%%%%%%%%%%%%%%%%%%%%%%%%%%%%%%%%%%%%%%%%%%%%%%%%%%%%%%%%%%%%%%%%%%%%%%%%%%%%%%%%%%%%%%%%%%%%%%%%%%%%%%%%%%%%%%%%%%%%%%%%%%%%%%%%%%%%%%%%%%%%%%%%%
\vskip0.12in
%%%%%%%%%%%%%%%%%%%%%%%%%%%%%%%%%%%%%%%%%%%%%%%%%%%%%%%%%%%%%%%%%%%%%%%%%%%%%%%%%%%%%%%%%%%%%%%%%%%%%%%%%%%%%%%%%%%%%%%%%%%%%%%%%%%%%%%%%%%%%%%%%%%%%%%%%%%%%%%%%%%%%%%%%%%%%%%%%%%%

For stable \(u\) of mass \(m\), \(z=U_m-u\) is a nonnegative Green potential with density \(\lambda(u^{-p}-U_m^{-p})\) in \(L^q(\Omega)\) for some \(q>1\), by \cref{thm:negative-full-representation}. Hence \(z\in H^1_0(\Omega)\). A nonzero \(z\) would give, by strict convexity and stability,
\[
 \int_\Omega|\nabla z|^2
 =\lambda\int_\Omega (u^{-p}-U_m^{-p})z
 <\lambda p\int_\Omega u^{-1-p}z^2
 \le\int_\Omega|\nabla z|^2 .
\]
The stable form extends to \(H^1_0\) by density and Fatou's lemma.
%%%%%%%%%%%%%%%%%%%%%%%%%%%%%%%%%%%%%%%%%%%%%%%%%%%%%%%%%%%%%%%%%%%%%%%%%%%%%%%%%%%%%%%%%%%%%%%%%%%%%%%%%%%%%%%%%%%%%%%%%%%%%%%%%%%%%%%%%%%%%%%%%%%%%%%%%%%%%%%%%%%%%%%%%%%%%%%%%%%%
\vskip0.12in
%%%%%%%%%%%%%%%%%%%%%%%%%%%%%%%%%%%%%%%%%%%%%%%%%%%%%%%%%%%%%%%%%%%%%%%%%%%%%%%%%%%%%%%%%%%%%%%%%%%%%%%%%%%%%%%%%%%%%%%%%%%%%%%%%%%%%%%%%%%%%%%%%%%%%%%%%%%%%%%%%%%%%%%%%%%%%%%%%%%%

To differentiate at \(m>m_*\), use \(w=U_m-U_{m_0}\), \(\delta=m-m_0\), and \(A_m z=\lambda pG_\Omega[U_m^{-1-p}z]\). Convexity gives \(A_mw\le w-\delta H\le(1-\delta/m)w\), since \(\delta H\le w\le mH\). Thus \(I-A_m\) is invertible with positive inverse in the equivalent norm \(\sup|z|/w\). On a relative ball about \(U_m\) the inverse-power Green term is smooth, since each derivative is dominated by a multiple of \(G_\Omega[H^{-p}]\) times the increment norms. The implicit function theorem applies there. For \(m'>m\), \(0\le U_{m'}-U_m\le(m'-m)H+A_m(U_{m'}-U_m)\), so the positive inverse gives right continuity in norm. For \(m'<m\), the local implicit-function solution is below \(U_{m'}\le U_m\), giving left continuity in norm. The maximal branch therefore agrees with the local implicit-function branch.

%%%%%%%%%%%%%%%%%%%%%%%%%%%%%%%%%%%%%%%%%%%%%%%%%%%%%%%%%%%%%%%%%%%%%%%%%%%%%%%%%%%%%%%%%%%%%%%%%%%%%%%%%%%%%%%%%%%%%%%%%%%%%%%%%%%%%%%%%%%%%%%%%%%%%%%%%%%%%%%%%%%%%%%%%%%%%%%%%%%%
\vskip0.12in
%%%%%%%%%%%%%%%%%%%%%%%%%%%%%%%%%%%%%%%%%%%%%%%%%%%%%%%%%%%%%%%%%%%%%%%%%%%%%%%%%%%%%%%%%%%%%%%%%%%%%%%%%%%%%%%%%%%%%%%%%%%%%%%%%%%%%%%%%%%%%%%%%%%%%%%%%%%%%%%%%%%%%%%%%%%%%%%%%%%%

Differentiation gives
\[
 \partial_m U_m=(I-A_m)^{-1}H\ge H,\qquad
 \partial_m^2U_m=-\lambda p(p+1)(I-A_m)^{-1}
 G_\Omega[U_m^{-p-2}(\partial_m U_m)^2]<0 .
\]
For a strict convex combination \(m\) of two admissible masses, their maximal solutions are strictly ordered. The corresponding convex combination \(v\) of these solutions satisfies \(T_m(v)>v\) by strict convexity. Iteration gives \(U_m\ge T_m(v)>v\), proving strict concavity also at the endpoint.
\end{proof}

%%%%%%%%%%%%%%%%%%%%%%%%%%%%%%%%%%%%%%%%%%%%%%%%%%%%%%%%%%%%%%%%%%%%%%%%%%%%%%%%%%%%%%%%%%%%%%%%%%%%%%%%%%%%%%%%%%%%%%%%%%%%%%%%%%%%%%%%%%%%%%%%%%%%%%%%%%%%%%%%%%%%%%%%%%%%%%%%%%%%
\vskip0.12in
%%%%%%%%%%%%%%%%%%%%%%%%%%%%%%%%%%%%%%%%%%%%%%%%%%%%%%%%%%%%%%%%%%%%%%%%%%%%%%%%%%%%%%%%%%%%%%%%%%%%%%%%%%%%%%%%%%%%%%%%%%%%%%%%%%%%%%%%%%%%%%%%%%%%%%%%%%%%%%%%%%%%%%%%%%%%%%%%%%%%

Taking \(H=K\) recovers the one-pole branch, since Theorem \ref{thm:negative-full-representation} identifies that solution class with its represented class. For general \(H\), the theorem classifies solutions of the displayed Green equation only; it asserts no independent classification of arbitrary boundary traces.

%%%%%%%%%%%%%%%%%%%%%%%%%%%%%%%%%%%%%%%%%%%%%%%%%%%%%%%%%%%%%%%%%%%%%%%%%%%%%%%%%%%%%%%%%%%%%%%%%%%%%%%%%%%%%%%%%%%%%%%%%%%%%%%%%%%%%%%%%%%%%%%%%%%%%%%%%%%%%%%%%%%%%%%%%%%%%%%%%%%%
\vskip0.22in
%%%%%%%%%%%%%%%%%%%%%%%%%%%%%%%%%%%%%%%%%%%%%%%%%%%%%%%%%%%%%%%%%%%%%%%%%%%%%%%%%%%%%%%%%%%%%%%%%%%%%%%%%%%%%%%%%%%%%%%%%%%%%%%%%%%%%%%%%%%%%%%%%%%%%%%%%%%%%%%%%%%%%%%%%%%%%%%%%%%%

\subsection{A distance lower bound for stable solutions}

%%%%%%%%%%%%%%%%%%%%%%%%%%%%%%%%%%%%%%%%%%%%%%%%%%%%%%%%%%%%%%%%%%%%%%%%%%%%%%%%%%%%%%%%%%%%%%%%%%%%%%%%%%%%%%%%%%%%%%%%%%%%%%%%%%%%%%%%%%%%%%%%%%%%%%%%%%%%%%%%%%%%%%%%%%%%%%%%%%%%
\vskip0.12in
%%%%%%%%%%%%%%%%%%%%%%%%%%%%%%%%%%%%%%%%%%%%%%%%%%%%%%%%%%%%%%%%%%%%%%%%%%%%%%%%%%%%%%%%%%%%%%%%%%%%%%%%%%%%%%%%%%%%%%%%%%%%%%%%%%%%%%%%%%%%%%%%%%%%%%%%%%%%%%%%%%%%%%%%%%%%%%%%%%%%

\begin{lemma}\label{lem:negative-stable-distance-seed}
For each \(0<p<1\) there is \(\underline c_p>0\), independent of the domain,
such that every stable positive planar solution of
\(\Delta u=\lambda u^{-p}\) satisfies
\[
 u(x)\ge \underline c_p\lambda^{1/(1+p)}d_\Omega(x)^{2/(1+p)}.
\]
In particular the estimate holds at the minimum mass.
\end{lemma}

\begin{proof}
Put \(a=1+p\) and \(\gamma=2/a\). If no universal lower bound
exists, rescaling at violating points gives stable positive solutions
\(U_j\) of \(\Delta U_j=U_j^{-p}\) on \(B_4\), with
\(U_j(0)\to0\). The logarithmic gradient estimate gives
\(|\nabla\log(1+U_j^{a/2})|\le C_p\) on \(B_2\), hence
uniform upper bounds on \(B_1\). Lemma~\ref{lem:negative-interior-compact-index}
gives a locally \(C^1\) convergent subsequence whose limit vanishes
at zero. The same lemma forces the Morse indices to diverge,
contradicting stability. Rescaling back proves the estimate,
including the stable minimum-mass solution.
\end{proof}

\subsection{Energy compactness at the minimum mass}\label{sec:negative-energy}

The Green correction admits a compactness statement in the Dirichlet
energy norm for the full class. At the stable endpoint its right
derivative is governed by a spectral range condition in that norm.

\begin{theorem}\label{thm:negative-energy-endpoint}
Assume \(\Lambda_p(H)<\infty\) in Theorem~\ref{thm:negative-branch}.
Let \(P_m=mH-U_m=\lambda G_\Omega[U_m^{-p}]\) on the maximal branch.
Then the following hold.
\begin{enumerate}
\item For any bounded-mass sequence of represented positive solutions,
there is a subsequence such that \(m_j\to m_0\), \(u_j\to u_0\) in
\(C^1_{\mathrm{loc}}(\Omega)\), and \(u_j-m_jH\to u_0-m_0H\)
in \(H^1_0(\Omega)\cap C_0(\overline\Omega)\), where the second
space has the uniform norm and zero boundary values.
The limit is positive almost everywhere and satisfies the Green equation almost everywhere. It is a positive classical solution if the Morse indices are bounded.
\item Set \(U_*=U_{m_*}\), \(P_*=P_{m_*}\), and
\(V_*=\lambda p U_*^{-p-1}\). Let \(\mathsf A\) on \(H^1_0(\Omega)\)
be defined in the Dirichlet inner product by
\[
 (\mathsf Af,g)_E=\int_\Omega V_*fg,\qquad
 (f,g)_E=\int_\Omega \nabla f\cdot\nabla g .
\]
This is a nonnegative self-adjoint contraction. Let \(\nu_{P_*}\)
be its spectral measure at \(P_*\), with mass \(\|P_*\|_E^2\).
The following conditions are equivalent.
\begin{enumerate}
\item \(U_m-mH\) has a right derivative at \(m_*\) in \(H^1_0(\Omega)\).
\item The equation
\begin{equation}\label{eq:negative-energy-endpoint-equation}
 (-\Delta-V_*)z=V_*H
\end{equation}
has a weak solution in \(H^1_0(\Omega)\).
\item \(\displaystyle\int_{[0,1]}(1-t)^{-2}\dd\nu_{P_*}(t)<\infty\).
\end{enumerate}
An atom at \(1\) makes the integral infinite.
When finite, the derivative is the nonnegative solution orthogonal
in the Dirichlet inner product to \(\ker(I-\mathsf A)\).
Otherwise the Dirichlet norms of the right difference quotients
increase to infinity as their steps decrease to zero.
\end{enumerate}
\end{theorem}

\begin{proof}
The Whitney proof of Theorem~\ref{thm:negative-full-representation}
bounds \(\|u_j^{-p}\|_{L^r}\) uniformly for \(r=1+\epsilon_p/p>1\).
Bounded mass gives local upper bounds through \(u_j\le m_j H\).
Lemma~\ref{lem:negative-interior-compact-index} yields local \(C^1\)
convergence after extraction. Fatou's lemma shows that \(u_0>0\)
almost everywhere. Uniform integrability then gives
\(u_j^{-p}\to u_0^{-p}\) in every \(L^s\) with \(1<s<r\).
The Green operator maps \(L^s\) continuously into \(H^1_0\) by the
planar Sobolev embedding, and into \(C_0(\overline\Omega)\) by the
uniform logarithmic kernel bound and approximation used in
Theorem~\ref{thm:negative-full-representation}. Taking limits
gives the represented equation and compactness of the corrections.
Bounded Morse index excludes any interior zero by
Lemma~\ref{lem:negative-interior-compact-index}.

For the derivative, stability gives the bounded form \(0\le\mathsf A\le I\).
The datum on the right of \eqref{eq:negative-energy-endpoint-equation}
always belongs to \(H^{-1}\). Indeed
\(m_*H=U_*+P_*\), \(V_*U_*=\lambda pU_*^{-p}\),
and \(P_*\in H^1_0\). Thus, using stability to pair \(V_*P_*\)
with energy tests,
\[
 b:=(-\Delta)^{-1}(V_*H)=m_*^{-1}(pI+\mathsf A)P_*.
\]
In the Dirichlet space the equation is \((I-\mathsf A)z=b\).
The spectral theorem gives a solution exactly when
\(\int (1-t)^{-2}\dd\nu_b(t)<\infty\), including the stated condition
at \(1\). Since \(p\le p+t\le p+1\), this is equivalent to the
condition with \(\nu_{P_*}\).
When it holds, the orthogonal solution \(z\) is the strong limit of
\(\sum_{j=0}^N\mathsf A^j b\), by spectral calculus.
It is nonnegative since \(b\ge0\) and \(\mathsf A\) maps nonnegative
functions to nonnegative functions by the weak maximum principle.

Let \(z_\delta=(P_*-P_{m_*+\delta})/\delta\), and put
\[
 V_\delta=\lambda p\int_0^1
    \big((1-s)U_*+sU_{m_*+\delta}\big)^{-1-p}\dd s .
\]
The branch order gives \(z_\delta\ge0\), and subtracting the two Green
equations gives \(-\Delta z_\delta=V_\delta(H+z_\delta)\).
These identities hold in \(H^{-1}\), since \(V_\delta\le V_*\),
\(V_*H\in H^{-1}\), and stability bounds the form on \(z_\delta\).
For each fixed \(\delta>0\) the form of \(-\Delta-V_\delta\) is
strictly coercive. In fact \(U_{m_*+\delta}/U_*\ge1+\delta/m_*\)
by the branch order and \(U_*\le m_*H\), so
\[
 V_\delta\le V_*\int_0^1(1+s\delta/m_*)^{-1-p}\dd s ,
\]
with the displayed integral strictly smaller than one.
When \(z\) exists, the positive-part test for that form shows
\(z_\delta\le z\). Both \(z_\delta\) and \(V_\delta\) increase
as \(\delta\downarrow0\), by concavity and monotonicity of the
branch, respectively. Hence
\[
 0\le -\Delta z_\delta\le V_*(H+z)=-\Delta z .
\]
Testing with the nonnegative \(z_\delta\) bounds its energy norm by
\(\|z\|_E\). Weak cluster points solve \eqref{eq:negative-energy-endpoint-equation}
by local convergence of \(V_\delta\) and weak convergence of the
tests. Each is between \(0\) and \(z\). There is no nonzero
nonnegative \(v\in\ker(I-\mathsf A)\) when \(z\) exists, since
\((b,v)_E=\langle V_*H,v\rangle>0\) for such \(v\).
Thus each cluster point is \(z\).
Testing the last order inequality with \(z-z_\delta\) yields
\(\|z-z_\delta\|_E^2\le(z,z-z_\delta)_E\to0\).

%%%%%%%%%%%%%%%%%%%%%%%%%%%%%%%%%%%%%%%%%%%%%%%%%%%%%%%%%%%%%%%%%%%%%%%%%%%%%%%%%%%%%%%%%%%%%%%%%%%%%%%%%%%%%%%%%%%%%%%%%%%%%%%%%%%%%%%%%%%%%%%%%%%%%%%%%%%%%%%%%%%%%%%%%%%%%%%%%%%%
\vskip0.12in
%%%%%%%%%%%%%%%%%%%%%%%%%%%%%%%%%%%%%%%%%%%%%%%%%%%%%%%%%%%%%%%%%%%%%%%%%%%%%%%%%%%%%%%%%%%%%%%%%%%%%%%%%%%%%%%%%%%%%%%%%%%%%%%%%%%%%%%%%%%%%%%%%%%%%%%%%%%%%%%%%%%%%%%%%%%%%%%%%%%%

Conversely, a bounded subsequence of \(z_\delta\) has a weak limit
which solves the endpoint equation by the same local passage.
If no solution exists, the norms diverge. They are nondecreasing as
\(\delta\) decreases, since the nonnegative sources
\(V_\delta(H+z_\delta)\) increase. Indeed testing their order with
the earlier \(z_\delta\) and using Cauchy--Schwarz gives the norm order.
\end{proof}

%%%%%%%%%%%%%%%%%%%%%%%%%%%%%%%%%%%%%%%%%%%%%%%%%%%%%%%%%%%%%%%%%%%%%%%%%%%%%%%%%%%%%%%%%%%%%%%%%%%%%%%%%%%%%%%%%%%%%%%%%%%%%%%%%%%%%%%%%%%%%%%%%%%%%%%%%%%%%%%%%%%%%%%%%%%%%%%%%%%%
\vskip0.32in
%%%%%%%%%%%%%%%%%%%%%%%%%%%%%%%%%%%%%%%%%%%%%%%%%%%%%%%%%%%%%%%%%%%%%%%%%%%%%%%%%%%%%%%%%%%%%%%%%%%%%%%%%%%%%%%%%%%%%%%%%%%%%%%%%%%%%%%%%%%%%%%%%%%%%%%%%%%%%%%%%%%%%%%%%%%%%%%%%%%%

\subsection{Morrey bounds up to the boundary}\label{sec:negative-compactness}

Use the exponents and distance potential
\[
  \gamma=\frac{2}{1+p},\qquad
  \sigma=2-\gamma=\frac{2p}{1+p},\qquad
  d(x)=\operatorname{dist}(x,\partial\Omega),\qquad
  P_\sigma=G_\Omega[d^{-\sigma}].
\]
In particular, \(1<\gamma<2\) and \(0<\sigma<1\). All constants
below are allowed to depend on \(\Omega\), \(p\), and the fixed
normalization of \(K\).

%%%%%%%%%%%%%%%%%%%%%%%%%%%%%%%%%%%%%%%%%%%%%%%%%%%%%%%%%%%%%%%%%%%%%%%%%%%%%%%%%%%%%%%%%%%%%%%%%%%%%%%%%%%%%%%%%%%%%%%%%%%%%%%%%%%%%%%%%%%%%%%%%%%%%%%%%%%%%%%%%%%%%%%%%%%%%%%%%%%%
\vskip0.12in
%%%%%%%%%%%%%%%%%%%%%%%%%%%%%%%%%%%%%%%%%%%%%%%%%%%%%%%%%%%%%%%%%%%%%%%%%%%%%%%%%%%%%%%%%%%%%%%%%%%%%%%%%%%%%%%%%%%%%%%%%%%%%%%%%%%%%%%%%%%%%%%%%%%%%%%%%%%%%%%%%%%%%%%%%%%%%%%%%%%%

\begin{lemma}\label{lem:morrey-green-domination}
Let \(1<\gamma<2\), \(\sigma=2-\gamma\), and let \(f\ge0\)
be locally integrable on a bounded planar Lipschitz domain. Suppose
\[
 \int_{B_r(x)}f\le A r^\gamma
 \quad\text{whenever }B_{4r}(x)\subset\Omega.
\]
There are constants depending only on the domain and \(\gamma\) such that
\begin{gather}
 \int_{B_R(z)\cap\Omega}f\le CA R^\gamma
 \quad(z\in\overline\Omega),\label{eq:morrey-boundary-extension}\\
 G_\Omega[f]\le CA G_\Omega[d^{-\sigma}],\qquad
 G_\Omega[\boldsymbol1_{\{d<\delta\}}f]
 \le CA G_\Omega[\boldsymbol1_{\{d<C\delta\}}d^{-\sigma}].
 \label{eq:morrey-green-localized}
\end{gather}
Also \(G_\Omega[d^{-\sigma}](x)\ge c d(x)^\gamma\).
\end{lemma}

\begin{proof}
Use a Whitney decomposition with squares small enough that their fixed
enlargements are interior. In a boundary ball of radius \(R\), the number
of squares of side \(\ell\lesssim R\) is \(O(R/\ell)\), by the
Lipschitz boundary-neighborhood area bound. Each square has forcing mass
at most \(CA\ell^\gamma\). Summing
\(CA R\ell^{\gamma-1}\) over dyadic scales proves
\eqref{eq:morrey-boundary-extension}; interior balls and large radii
are handled by the original bound and a finite chart cover.

%%%%%%%%%%%%%%%%%%%%%%%%%%%%%%%%%%%%%%%%%%%%%%%%%%%%%%%%%%%%%%%%%%%%%%%%%%%%%%%%%%%%%%%%%%%%%%%%%%%%%%%%%%%%%%%%%%%%%%%%%%%%%%%%%%%%%%%%%%%%%%%%%%%%%%%%%%%%%%%%%%%%%%%%%%%%%%%%%%%%
\vskip0.12in
%%%%%%%%%%%%%%%%%%%%%%%%%%%%%%%%%%%%%%%%%%%%%%%%%%%%%%%%%%%%%%%%%%%%%%%%%%%%%%%%%%%%%%%%%%%%%%%%%%%%%%%%%%%%%%%%%%%%%%%%%%%%%%%%%%%%%%%%%%%%%%%%%%%%%%%%%%%%%%%%%%%%%%%%%%%%%%%%%%%%

Fix a target \(x\). On a square whose fixed enlargement avoids \(x\),
Harnack comparison for \(G_\Omega(x,\cdot)\) and the square mass bound
give
\[
 \int_S G_\Omega(x,y)f(y)\,dy
 \le CA\int_S G_\Omega(x,y)d(y)^{-\sigma}\,dy.
\]
The other squares lie in \(B_{Cd(x)}(x)\). On this ball
\(G_\Omega(x,y)\le C\log(C'd(x)/|x-y|)\). For small \(d(x)\),
an exterior disk of radius comparable to \(d(x)\) follows from the
exterior Lipschitz cone; inversion about its center maps the domain
into a disk and gives this bound from the disk Green function.
For \(d(x)\) bounded below, comparison with an enclosing disk suffices.
The dyadic shell sum using \eqref{eq:morrey-boundary-extension}
bounds the near contribution by \(CA d(x)^\gamma\).
The Green function of \(B_{d(x)/2}(x)\), by domain monotonicity,
gives the stated lower bound for the distance potential.
This proves the global estimate.

%%%%%%%%%%%%%%%%%%%%%%%%%%%%%%%%%%%%%%%%%%%%%%%%%%%%%%%%%%%%%%%%%%%%%%%%%%%%%%%%%%%%%%%%%%%%%%%%%%%%%%%%%%%%%%%%%%%%%%%%%%%%%%%%%%%%%%%%%%%%%%%%%%%%%%%%%%%%%%%%%%%%%%%%%%%%%%%%%%%%
\vskip0.12in
%%%%%%%%%%%%%%%%%%%%%%%%%%%%%%%%%%%%%%%%%%%%%%%%%%%%%%%%%%%%%%%%%%%%%%%%%%%%%%%%%%%%%%%%%%%%%%%%%%%%%%%%%%%%%%%%%%%%%%%%%%%%%%%%%%%%%%%%%%%%%%%%%%%%%%%%%%%%%%%%%%%%%%%%%%%%%%%%%%%%

For the localized estimate retain only squares meeting \(\{d<\delta\}\).
Every such square lies in \(\{d<C\delta\}\). The same far-square
comparison applies. A contributing near square forces
\(d(x)\lesssim\delta\), so the ball \(B_{d(x)/2}(x)\) used for
the lower bound also lies in the enlarged collar. This proves
the localized estimate without any regularity assumption on \(f\).
\end{proof}

\begin{proposition}\label{prop:negative-distance-potential-chart}
For \(0<\rho<1\), one has \(G_\Omega[d^{-\rho}]\lesssim d\) near any compact smooth boundary arc. At a smooth pole \(O\) its ratio to \(K\) tends to zero. If the whole boundary is smooth, this ratio is bounded on \(\Omega\).
\end{proposition}

\begin{proof}
The density belongs to \(L^q(\Omega)\) for some \(q>1\). In a regular conformal half-disc its pullback with the Jacobian is bounded by \(Cy^{-\rho}\). Its odd reflection has Morrey exponent \(2-\rho>1\). The Newtonian gradient estimate of \cref{lem:negative-interior-compact-index} and harmonic reflection give a bounded gradient near the smaller diameter. The zero boundary value gives the first bound. At a smooth pole use \(K\circ F\asymp y/|w|^2\). On a smooth domain a finite cover and \(K\gtrsim d\) give the global bound.
\end{proof}

%%%%%%%%%%%%%%%%%%%%%%%%%%%%%%%%%%%%%%%%%%%%%%%%%%%%%%%%%%%%%%%%%%%%%%%%%%%%%%%%%%%%%%%%%%%%%%%%%%%%%%%%%%%%%%%%%%%%%%%%%%%%%%%%%%%%%%%%%%%%%%%%%%%%%%%%%%%%%%%%%%%%%%%%%%%%%%%%%%%%
\vskip0.12in
%%%%%%%%%%%%%%%%%%%%%%%%%%%%%%%%%%%%%%%%%%%%%%%%%%%%%%%%%%%%%%%%%%%%%%%%%%%%%%%%%%%%%%%%%%%%%%%%%%%%%%%%%%%%%%%%%%%%%%%%%%%%%%%%%%%%%%%%%%%%%%%%%%%%%%%%%%%%%%%%%%%%%%%%%%%%%%%%%%%%

For positive solutions, apply \cref{lem:morrey-green-domination} with the interior mass bound from \eqref{eq:negative-inverse-ball} to obtain
\begin{gather}
 \int_{B_R(z)\cap\Omega}u^{-p}
 \le C_{\Omega,p}\lambda^{-p/(1+p)}R^\gamma,
 \label{eq:negative-global-morrey}\\
 G_\Omega[u^{-p}]
 \le C_{\Omega,p}\lambda^{-p/(1+p)}P_\sigma,
 \label{eq:negative-large-potential}\\
 P_\sigma(x)\ge c_p d(x)^\gamma.
 \label{eq:negative-large-lower}
\end{gather}

%%%%%%%%%%%%%%%%%%%%%%%%%%%%%%%%%%%%%%%%%%%%%%%%%%%%%%%%%%%%%%%%%%%%%%%%%%%%%%%%%%%%%%%%%%%%%%%%%%%%%%%%%%%%%%%%%%%%%%%%%%%%%%%%%%%%%%%%%%%%%%%%%%%%%%%%%%%%%%%%%%%%%%%%%%%%%%%%%%%%
\vskip0.32in
%%%%%%%%%%%%%%%%%%%%%%%%%%%%%%%%%%%%%%%%%%%%%%%%%%%%%%%%%%%%%%%%%%%%%%%%%%%%%%%%%%%%%%%%%%%%%%%%%%%%%%%%%%%%%%%%%%%%%%%%%%%%%%%%%%%%%%%%%%%%%%%%%%%%%%%%%%%%%%%%%%%%%%%%%%%%%%%%%%%%

\subsection{Endpoint uniqueness and quantitative convergence}

The following mass-lowering comparison extends the interior
surplus argument in the proof of \cite[Theorem~1.4]{HuLongTaoWuZou}
to domains that are Lipschitz away from the pole.

%%%%%%%%%%%%%%%%%%%%%%%%%%%%%%%%%%%%%%%%%%%%%%%%%%%%%%%%%%%%%%%%%%%%%%%%%%%%%%%%%%%%%%%%%%%%%%%%%%%%%%%%%%%%%%%%%%%%%%%%%%%%%%%%%%%%%%%%%%%%%%%%%%%%%%%%%%%%%%%%%%%%%%%%%%%%%%%%%%%%
\vskip0.12in
%%%%%%%%%%%%%%%%%%%%%%%%%%%%%%%%%%%%%%%%%%%%%%%%%%%%%%%%%%%%%%%%%%%%%%%%%%%%%%%%%%%%%%%%%%%%%%%%%%%%%%%%%%%%%%%%%%%%%%%%%%%%%%%%%%%%%%%%%%%%%%%%%%%%%%%%%%%%%%%%%%%%%%%%%%%%%%%%%%%%

\begin{lemma}\label{lem:negative-interior-surplus}
Let \(k\) be the Martin kernel at a single pole \(o\)
of a bounded planar Lipschitz domain \(D\).
Fix a compact set \(E\Subset D\), a positive coefficient \(J\),
and constants \(c,m_0>0\). Assume \(J\) is locally H\"older
continuous and \(G_D[Jk^{-p}]/k\)
is bounded and tends to zero at \(o\). Suppose \(0<J'\le J\) is locally H\"older continuous,
\(m\ge m_0\), and \(W\in C(D)\) satisfies
\[
 W=m k-G_D[F]>0,\qquad F\ge \lambda J'W^{-p}+g,\qquad
 W\ge c k\ \text{near }o,
\]
where \(g\ge0\) is bounded and supported in \(E\).
Then the \(J'\)-problem has a positive solution at mass
\(m-c_1\min\{\int_D g,1\}\) for some \(c_1>0\).
The constant is uniform if the neighborhood in the lower
bound for \(W\) and the other specified data are fixed.
\end{lemma}

%%%%%%%%%%%%%%%%%%%%%%%%%%%%%%%%%%%%%%%%%%%%%%%%%%%%%%%%%%%%%%%%%%%%%%%%%%%%%%%%%%%%%%%%%%%%%%%%%%%%%%%%%%%%%%%%%%%%%%%%%%%%%%%%%%%%%%%%%%%%%%%%%%%%%%%%%%%%%%%%%%%%%%%%%%%%%%%%%%%%
\vskip0.12in
%%%%%%%%%%%%%%%%%%%%%%%%%%%%%%%%%%%%%%%%%%%%%%%%%%%%%%%%%%%%%%%%%%%%%%%%%%%%%%%%%%%%%%%%%%%%%%%%%%%%%%%%%%%%%%%%%%%%%%%%%%%%%%%%%%%%%%%%%%%%%%%%%%%%%%%%%%%%%%%%%%%%%%%%%%%%%%%%%%%%

\begin{proof}
If \(\int g=0\), \(W\) itself is a lower barrier.
For \(\int_D g>0\), choose a pole neighborhood \(P\) disjoint from \(E\) on which \(W\ge ck\), and set \(h=\boldsymbol 1_PJk^{-p}\). For \(0<\delta<c/4\) and \(0<Z\le2k\) there, decreasing \(W\) to \(W-\delta Z\) increases \(\lambda J'W^{-p}\) by at most \(\delta BJk^{-p}\), by the mean value theorem. For \(Z\le0\) the inverse power decreases. Shrink \(P\) until \(BG_D[h]\le k\) on \(P\).

For \(x\) off \(P\), \(G_D(x,y)\gtrsim k(x)\) uniformly in \(y\in E\). In ordinary boundary charts disjoint from \(E\) and the pole, apply boundary Harnack to \(G_D(\cdot,y)\) and \(k\), using positive lower bounds at fixed interior reference points \cite[Theorem~1]{WuBoundaryHarnack}. Compactness and positivity give the estimate on the remaining interior set. Since \(G_D[h]/k\) is bounded, choose \(C\) large enough that
\[
 Z=k+BG_D[h]-C G_D[g]/\!\int_D g
\]
satisfies \(Z\le0\) off \(P\). On \(P\) we have \(Z\le2k\). Take \(\delta=c_1\min\{\int_D g,1\}\) with \(0<c_1<\min(c/4,m_0/2,1/C)\). Then \(\widetilde W=W-\delta Z>0\) and
\[
 \widetilde W=(m-\delta)k-G_D[\widetilde f],
 \qquad \widetilde f=F+\delta Bh-\frac{\delta C}{\int_D g}g.
\]
On \(P\) the extra source \(\delta Bh\) pays for the increase in the inverse power when \(Z>0\). Off \(P\), \(\widetilde W\ge W\), and \(\delta Cg/\int_D g\le g\) is absorbed by the original surplus. Therefore \(\widetilde f\ge\lambda J'\widetilde W^{-p}\) on \(D\), so \(\widetilde W\) is a positive lower barrier for Green iteration at mass \(m-\delta\).
\end{proof}

%%%%%%%%%%%%%%%%%%%%%%%%%%%%%%%%%%%%%%%%%%%%%%%%%%%%%%%%%%%%%%%%%%%%%%%%%%%%%%%%%%%%%%%%%%%%%%%%%%%%%%%%%%%%%%%%%%%%%%%%%%%%%%%%%%%%%%%%%%%%%%%%%%%%%%%%%%%%%%%%%%%%%%%%%%%%%%%%%%%%
\vskip0.12in
%%%%%%%%%%%%%%%%%%%%%%%%%%%%%%%%%%%%%%%%%%%%%%%%%%%%%%%%%%%%%%%%%%%%%%%%%%%%%%%%%%%%%%%%%%%%%%%%%%%%%%%%%%%%%%%%%%%%%%%%%%%%%%%%%%%%%%%%%%%%%%%%%%%%%%%%%%%%%%%%%%%%%%%%%%%%%%%%%%%%

\begin{theorem}\label{thm:negative-minimum-unique}
Suppose the positive one-pole class is nonempty, and assume
\begin{equation}\label{eq:negative-pole-flat-potential}
 \lim_{x\to O}\frac{G_\Omega[U_{m_*}^{-p}](x)}{K(x)}=0.
\end{equation}
Then \(U_{m_*}\) is the unique solution at the minimum mass.
In particular, the condition holds if the pole lies on a smooth
boundary arc.
\label{prop:negative-endpoint-quant}
Moreover, for every compact \(E\subset\Omega\)
and \(M>m_*\), solutions of mass \(m\in[m_*,M]\) satisfy
\[
 \int_E |u-U_{m_*}|^2\le C_{E,M}(m-m_*).
\]
For \(m\) close enough to \(m_*\), the rate \((m-m_*)^{1/2}\) also holds in \(C^{2,\alpha}(E)\) for each \(0<\alpha<1\).
\end{theorem}

%%%%%%%%%%%%%%%%%%%%%%%%%%%%%%%%%%%%%%%%%%%%%%%%%%%%%%%%%%%%%%%%%%%%%%%%%%%%%%%%%%%%%%%%%%%%%%%%%%%%%%%%%%%%%%%%%%%%%%%%%%%%%%%%%%%%%%%%%%%%%%%%%%%%%%%%%%%%%%%%%%%%%%%%%%%%%%%%%%%%
\vskip0.12in
%%%%%%%%%%%%%%%%%%%%%%%%%%%%%%%%%%%%%%%%%%%%%%%%%%%%%%%%%%%%%%%%%%%%%%%%%%%%%%%%%%%%%%%%%%%%%%%%%%%%%%%%%%%%%%%%%%%%%%%%%%%%%%%%%%%%%%%%%%%%%%%%%%%%%%%%%%%%%%%%%%%%%%%%%%%%%%%%%%%%

\begin{proof}
Write \(U_*=U_{m_*}\). The representation and
\eqref{eq:negative-pole-flat-potential} give \(U_*/K\to m_*>0\)
at \(O\). Since \(K\ge U_*/m_*\), the same assumption also gives
\(G_\Omega[K^{-p}]/K\to0\) there. The smooth-pole case follows from \eqref{eq:negative-large-potential} and
Proposition~\ref{prop:negative-distance-potential-chart}.

%%%%%%%%%%%%%%%%%%%%%%%%%%%%%%%%%%%%%%%%%%%%%%%%%%%%%%%%%%%%%%%%%%%%%%%%%%%%%%%%%%%%%%%%%%%%%%%%%%%%%%%%%%%%%%%%%%%%%%%%%%%%%%%%%%%%%%%%%%%%%%%%%%%%%%%%%%%%%%%%%%%%%%%%%%%%%%%%%%%%
\vskip0.12in
%%%%%%%%%%%%%%%%%%%%%%%%%%%%%%%%%%%%%%%%%%%%%%%%%%%%%%%%%%%%%%%%%%%%%%%%%%%%%%%%%%%%%%%%%%%%%%%%%%%%%%%%%%%%%%%%%%%%%%%%%%%%%%%%%%%%%%%%%%%%%%%%%%%%%%%%%%%%%%%%%%%%%%%%%%%%%%%%%%%%

Choose \(0\le\chi\le1\) smooth, compactly supported, and equal to
one on \(E\). Put \(W=(u+U_*)/2\) and \(g=\chi(u-U_*)^2\). The second derivative of \(s^{-p}\) is bounded below by a positive constant for \(0<s\le M\sup_{\operatorname{supp}\chi}K\). Convexity gives
\(\Delta W-\lambda W^{-p}\ge b g\), with \(b>0\) independent
of \(u\). Near \(O\), \(W\ge U_*/2\ge cK\).

The pole limit and the finite inverse-power row allow
Lemma~\ref{lem:negative-interior-surplus} with \(J'=J=1\),
\(W\) as above, surplus \(bg\), and mass \((m+m_*)/2\).
The mass cap bounds \(\int g\) uniformly, so the lemma gives
\(m_*\le(m+m_*)/2-c_1\int g\), with a fixed \(c_1>0\).
This proves the integral estimate and uniqueness at \(m_*\).

%%%%%%%%%%%%%%%%%%%%%%%%%%%%%%%%%%%%%%%%%%%%%%%%%%%%%%%%%%%%%%%%%%%%%%%%%%%%%%%%%%%%%%%%%%%%%%%%%%%%%%%%%%%%%%%%%%%%%%%%%%%%%%%%%%%%%%%%%%%%%%%%%%%%%%%%%%%%%%%%%%%%%%%%%%%%%%%%%%%%
\vskip0.12in
%%%%%%%%%%%%%%%%%%%%%%%%%%%%%%%%%%%%%%%%%%%%%%%%%%%%%%%%%%%%%%%%%%%%%%%%%%%%%%%%%%%%%%%%%%%%%%%%%%%%%%%%%%%%%%%%%%%%%%%%%%%%%%%%%%%%%%%%%%%%%%%%%%%%%%%%%%%%%%%%%%%%%%%%%%%%%%%%%%%%

For the local rate, the uniform \(L^{1+\epsilon_p/p}\) forcing
bound and \(u\le MK\) give uniform interior H\"older bounds.
The \(L^2\) estimate on larger compact neighborhoods therefore
implies local uniform convergence to the positive \(U_*\).
On nested neighborhoods the difference solves a linear
equation with uniformly bounded H\"older coefficient, by the
smoothness of the inverse power on the uniform positive range
and interior elliptic estimates. The local \(L^2\)-to-supremum
and Schauder estimates give the asserted rate on \(E\).
\end{proof}

%%%%%%%%%%%%%%%%%%%%%%%%%%%%%%%%%%%%%%%%%%%%%%%%%%%%%%%%%%%%%%%%%%%%%%%%%%%%%%%%%%%%%%%%%%%%%%%%%%%%%%%%%%%%%%%%%%%%%%%%%%%%%%%%%%%%%%%%%%%%%%%%%%%%%%%%%%%%%%%%%%%%%%%%%%%%%%%%%%%%
\vskip0.12in
%%%%%%%%%%%%%%%%%%%%%%%%%%%%%%%%%%%%%%%%%%%%%%%%%%%%%%%%%%%%%%%%%%%%%%%%%%%%%%%%%%%%%%%%%%%%%%%%%%%%%%%%%%%%%%%%%%%%%%%%%%%%%%%%%%%%%%%%%%%%%%%%%%%%%%%%%%%%%%%%%%%%%%%%%%%%%%%%%%%%

%%%%%%%%%%%%%%%%%%%%%%%%%%%%%%%%%%%%%%%%%%%%%%%%%%%%%%%%%%%%%%%%%%%%%%%%%%%%%%%%%%%%%%%%%%%%%%%%%%%%%%%%%%%%%%%%%%%%%%%%%%%%%%%%%%%%%%%%%%%%%%%%%%%%%%%%%%%%%%%%%%%%%%%%%%%%%%%%%%%%
\vskip0.32in
%%%%%%%%%%%%%%%%%%%%%%%%%%%%%%%%%%%%%%%%%%%%%%%%%%%%%%%%%%%%%%%%%%%%%%%%%%%%%%%%%%%%%%%%%%%%%%%%%%%%%%%%%%%%%%%%%%%%%%%%%%%%%%%%%%%%%%%%%%%%%%%%%%%%%%%%%%%%%%%%%%%%%%%%%%%%%%%%%%%%

\section{\bf Smooth boundary profiles and the sharp energy constant}
\label{sec:sharp-smooth}

The potential at a smooth zero-slope point is governed by the one-dimensional zero-slope profile on fixed normal scales. To estimate its local form, one must also allow the base point to move along the boundary. The scale at that base point is determined by the ratio of its inward slope to the height raised to the power \(\alpha=(1-p)/(1+p)\). We identify the limiting value of the potential for every value of this ratio.

%%%%%%%%%%%%%%%%%%%%%%%%%%%%%%%%%%%%%%%%%%%%%%%%%%%%%%%%%%%%%%%%%%%%%%%%%%%%%%%%%%%%%%%%%%%%%%%%%%%%%%%%%%%%%%%%%%%%%%%%%%%%%%%%%%%%%%%%%%%%%%%%%%%%%%%%%%%%%%%%%%%%%%%%%%%%%%%%%%%%
\vskip0.22in
%%%%%%%%%%%%%%%%%%%%%%%%%%%%%%%%%%%%%%%%%%%%%%%%%%%%%%%%%%%%%%%%%%%%%%%%%%%%%%%%%%%%%%%%%%%%%%%%%%%%%%%%%%%%%%%%%%%%%%%%%%%%%%%%%%%%%%%%%%%%%%%%%%%%%%%%%%%%%%%%%%%%%%%%%%%%%%%%%%%%

\subsection{The transition law}
Let
\begin{equation}\label{eq:sharp-exponents}
 \gamma=\frac{2}{1+p}=1+\alpha,\qquad
 h_p=p\gamma\alpha=\frac{2p(1-p)}{(1+p)^2},\qquad
 \beta_p=4h_p .
\end{equation}
Use a conformal half-disc chart \(F\) of a \(C^{1,\eta}\) boundary arc, \(0<\eta<1\), with its diameter on \(\mathbb R\). The Kellogg--Warschawski theorem gives \(C^{1,\eta}\) regularity and nonvanishing derivative on a smaller diameter \cite{Pommerenke}. Let \(u\) be a positive solution of \(\Delta u=\lambda u^{-p}\) that vanishes continuously on the arc. For \(v=u\circ F\) and \(\lambda_F=\lambda |F'|^2\), put
\begin{equation}\label{eq:sharp-chart-slope}
 B_u(s,t)=\frac{\partial_y v(s,0)}{\lambda_F(s)^{1/(1+p)}t^\alpha},\qquad t>0.
\end{equation}
For \(B\geq0\) define \(w_B\) by
\begin{equation}\label{eq:sharp-one-dimensional}
 y=\int_0^{w_B(y)}
       \biggl(B^2+\frac{2}{1-p}z^{1-p}\biggr)^{-1/2}\dd z .
\end{equation}
The integral increases from zero to infinity as its upper limit increases. The first integral identifies \(w_B\) as the unique positive solution of \(w_B''=w_B^{-p}\) with \(w_B(0)=0\) and \(w_B'(0)=B\). In particular \(w_0(y)=A y^\gamma\) with \(A^{1+p}=(\gamma\alpha)^{-1}\). For each \(y>0\), the integral formula and dominated convergence give continuity and strict increase of \(w_B(y)\) in \(B\). Also \(w_B(y)\geq By\).

%%%%%%%%%%%%%%%%%%%%%%%%%%%%%%%%%%%%%%%%%%%%%%%%%%%%%%%%%%%%%%%%%%%%%%%%%%%%%%%%%%%%%%%%%%%%%%%%%%%%%%%%%%%%%%%%%%%%%%%%%%%%%%%%%%%%%%%%%%%%%%%%%%%%%%%%%%%%%%%%%%%%%%%%%%%%%%%%%%%%
\vskip0.12in
%%%%%%%%%%%%%%%%%%%%%%%%%%%%%%%%%%%%%%%%%%%%%%%%%%%%%%%%%%%%%%%%%%%%%%%%%%%%%%%%%%%%%%%%%%%%%%%%%%%%%%%%%%%%%%%%%%%%%%%%%%%%%%%%%%%%%%%%%%%%%%%%%%%%%%%%%%%%%%%%%%%%%%%%%%%%%%%%%%%%

For two solutions \(u,\widetilde u\), use the diagonal and secant potentials
\begin{equation}\label{eq:sharp-secant}
 V_u=\lambda p u^{-1-p},\qquad
 V_{u,\widetilde u}=\lambda p\int_0^1
   \bigl((1-\tau)u+\tau\widetilde u\bigr)^{-1-p}\dd\tau .
\end{equation}
They agree on the diagonal. Set
\begin{equation}\label{eq:sharp-Phi-Psi}
 \Phi_p(B)=p w_B(1)^{-1-p},\qquad
 \Psi_p(B,C)=p\int_0^1
   \bigl((1-\tau)w_B(1)+\tau w_C(1)\bigr)^{-1-p}\dd\tau .
\end{equation}
Extend \(\Phi_p\) by zero at infinity and \(\Psi_p\) by zero when either argument is infinite. These extensions are continuous, by \(w_B(1)\ge A\), \(w_B(1)\to\infty\), and dominated convergence. The function \(\Phi_p\) decreases strictly from \(h_p\) to zero.

%%%%%%%%%%%%%%%%%%%%%%%%%%%%%%%%%%%%%%%%%%%%%%%%%%%%%%%%%%%%%%%%%%%%%%%%%%%%%%%%%%%%%%%%%%%%%%%%%%%%%%%%%%%%%%%%%%%%%%%%%%%%%%%%%%%%%%%%%%%%%%%%%%%%%%%%%%%%%%%%%%%%%%%%%%%%%%%%%%%%
\vskip0.12in
%%%%%%%%%%%%%%%%%%%%%%%%%%%%%%%%%%%%%%%%%%%%%%%%%%%%%%%%%%%%%%%%%%%%%%%%%%%%%%%%%%%%%%%%%%%%%%%%%%%%%%%%%%%%%%%%%%%%%%%%%%%%%%%%%%%%%%%%%%%%%%%%%%%%%%%%%%%%%%%%%%%%%%%%%%%%%%%%%%%%

%%%%%%%%%%%%%%%%%%%%%%%%%%%%%%%%%%%%%%%%%%%%%%%%%%%%%%%%%%%%%%%%%%%%%%%%%%%%%%%%%%%%%%%%%%%%%%%%%%%%%%%%%%%%%%%%%%%%%%%%%%%%%%%%%%%%%%%%%%%%%%%%%%%%%%%%%%%%%%%%%%%%%%%%%%%%%%%%%%%%
\vskip0.12in
%%%%%%%%%%%%%%%%%%%%%%%%%%%%%%%%%%%%%%%%%%%%%%%%%%%%%%%%%%%%%%%%%%%%%%%%%%%%%%%%%%%%%%%%%%%%%%%%%%%%%%%%%%%%%%%%%%%%%%%%%%%%%%%%%%%%%%%%%%%%%%%%%%%%%%%%%%%%%%%%%%%%%%%%%%%%%%%%%%%%

\begin{theorem}\label{thm:sharp-boundary-transition}
Let \(\mathcal U\) be a family of positive solutions of \(\Delta u=\lambda u^{-p}\) near a \(C^{1,\eta}\) boundary arc, vanishing continuously on the arc and uniformly bounded on a fixed neighborhood of a compact subarc. For \(s\) on a smaller diameter and \(x=F(s+it)\), one has
\begin{align}
 d(x)^2 V_u(x)&=\Phi_p(B_u(s,t))+o(1), \label{eq:sharp-transition}\\
 d(x)^2 V_{u,\widetilde u}(x)&=\Psi_p(B_u(s,t),B_{\widetilde u}(s,t))+o(1)
                                                        \label{eq:sharp-transition-secant}
\end{align}
as \(t\downarrow0\), uniformly in \(s\) and \(u,\widetilde u\in\mathcal U\). The inward slopes in \eqref{eq:sharp-chart-slope} exist and are continuous. In particular \(d^2V_{u,\widetilde u}\le h_p+o(1)\) with the same uniformity.

%%%%%%%%%%%%%%%%%%%%%%%%%%%%%%%%%%%%%%%%%%%%%%%%%%%%%%%%%%%%%%%%%%%%%%%%%%%%%%%%%%%%%%%%%%%%%%%%%%%%%%%%%%%%%%%%%%%%%%%%%%%%%%%%%%%%%%%%%%%%%%%%%%%%%%%%%%%%%%%%%%%%%%%%%%%%%%%%%%%%
\vskip0.12in
%%%%%%%%%%%%%%%%%%%%%%%%%%%%%%%%%%%%%%%%%%%%%%%%%%%%%%%%%%%%%%%%%%%%%%%%%%%%%%%%%%%%%%%%%%%%%%%%%%%%%%%%%%%%%%%%%%%%%%%%%%%%%%%%%%%%%%%%%%%%%%%%%%%%%%%%%%%%%%%%%%%%%%%%%%%%%%%%%%%%

At a fixed point \(P\) of this arc the exact local energy constant is
\begin{equation}\label{eq:sharp-local-energy}
 b_P(V_{u,\widetilde u})=
 \begin{cases}
 \displaystyle\frac{8p(1-p)}{(1+p)^2},
     &\partial_\nu u(P)=\partial_\nu\widetilde u(P)=0,\\[5pt]
 0,&\max(\partial_\nu u(P),\partial_\nu\widetilde u(P))>0 .
 \end{cases}
\end{equation}
Here \(\nu\) is the inward unit normal and
\begin{equation*}
 b_P(V)=\lim_{\rho\downarrow0}
 \sup_{0\ne f\in H^1_0(\Omega\cap B_\rho(P))}
 \frac{\displaystyle\int_{\Omega\cap B_\rho(P)}Vf^2}
      {\displaystyle\int_{\Omega\cap B_\rho(P)}|\nabla f|^2}.
\end{equation*}
\end{theorem}

%%%%%%%%%%%%%%%%%%%%%%%%%%%%%%%%%%%%%%%%%%%%%%%%%%%%%%%%%%%%%%%%%%%%%%%%%%%%%%%%%%%%%%%%%%%%%%%%%%%%%%%%%%%%%%%%%%%%%%%%%%%%%%%%%%%%%%%%%%%%%%%%%%%%%%%%%%%%%%%%%%%%%%%%%%%%%%%%%%%%
\vskip0.12in
%%%%%%%%%%%%%%%%%%%%%%%%%%%%%%%%%%%%%%%%%%%%%%%%%%%%%%%%%%%%%%%%%%%%%%%%%%%%%%%%%%%%%%%%%%%%%%%%%%%%%%%%%%%%%%%%%%%%%%%%%%%%%%%%%%%%%%%%%%%%%%%%%%%%%%%%%%%%%%%%%%%%%%%%%%%%%%%%%%%%

The two-potential statement will be used to compare distinct solutions at the same mass. On the diagonal at a zero-slope point it also gives
\(\limsup_{x\to P}d(x)^2V_u(x)=h_p\) and \(b_P((V_u-h_p/d^2)_+)=0\).

%%%%%%%%%%%%%%%%%%%%%%%%%%%%%%%%%%%%%%%%%%%%%%%%%%%%%%%%%%%%%%%%%%%%%%%%%%%%%%%%%%%%%%%%%%%%%%%%%%%%%%%%%%%%%%%%%%%%%%%%%%%%%%%%%%%%%%%%%%%%%%%%%%%%%%%%%%%%%%%%%%%%%%%%%%%%%%%%%%%%
\vskip0.22in
%%%%%%%%%%%%%%%%%%%%%%%%%%%%%%%%%%%%%%%%%%%%%%%%%%%%%%%%%%%%%%%%%%%%%%%%%%%%%%%%%%%%%%%%%%%%%%%%%%%%%%%%%%%%%%%%%%%%%%%%%%%%%%%%%%%%%%%%%%%%%%%%%%%%%%%%%%%%%%%%%%%%%%%%%%%%%%%%%%%%

\subsection{Positivity of half-plane limits}

We first rule out interior zeros in the limits needed for
\eqref{eq:sharp-transition}. After normalization the limiting equation has coefficient one. Consider a nonnegative \(C^1\) function \(v\) on the closed upper half-plane, vanishing on its boundary and with a finite global \(\alpha\)-H\"older seminorm of its gradient. The limits considered below have the following distributional properties in the open half-plane:
\begin{equation}\label{eq:sharp-stationary-class}
 \Delta v=v^{-p}\dd x+\mu_0,\qquad
 v>0\ \text{a.e.},\qquad
 \mu_0\geq0,\quad \operatorname{supp}\mu_0\subset\{v=0\},\qquad
 \operatorname{div} T_v=0,
\end{equation}
where \(v^{-p}\in L^1_{\mathrm{loc}}\) and
\begin{equation*}
 T_v=\left(|\nabla v|^2+\frac{2}{1-p}v^{1-p}\right)I
             -2\nabla v\otimes\nabla v .
\end{equation*}
The same properties are used for limits on the whole plane. For the distributional treatment of finite-energy inverse-power solutions near their zero set, see \cite[Theorem~1.2]{DavilaPonce}.

%%%%%%%%%%%%%%%%%%%%%%%%%%%%%%%%%%%%%%%%%%%%%%%%%%%%%%%%%%%%%%%%%%%%%%%%%%%%%%%%%%%%%%%%%%%%%%%%%%%%%%%%%%%%%%%%%%%%%%%%%%%%%%%%%%%%%%%%%%%%%%%%%%%%%%%%%%%%%%%%%%%%%%%%%%%%%%%%%%%%
\vskip0.12in
%%%%%%%%%%%%%%%%%%%%%%%%%%%%%%%%%%%%%%%%%%%%%%%%%%%%%%%%%%%%%%%%%%%%%%%%%%%%%%%%%%%%%%%%%%%%%%%%%%%%%%%%%%%%%%%%%%%%%%%%%%%%%%%%%%%%%%%%%%%%%%%%%%%%%%%%%%%%%%%%%%%%%%%%%%%%%%%%%%%%

These properties are closed under local \(C^1\) limits with a common gradient seminorm bound, including limits on expanding domains. Indeed, for any nonnegative interior cutoff \(\chi\) equal to one on a given compact set, the source measure on that set is bounded by \(-\int\nabla\chi\cdot\nabla v\). Applied along the sequence, Fatou's lemma gives local inverse-power integrability and strict positivity almost everywhere in the limit. On compact subsets of the positive set the measures converge to the classical source, by uniform convergence and the positive lower bound there. The excess measure is nonnegative and supported on the zero set. The stress tensors converge locally uniformly because \(1-p>0\). Finally, testing the measure equation against \(v\) times a compact cutoff yields
\begin{equation}\label{eq:sharp-product}
 \operatorname{div}(v\nabla v)=|\nabla v|^2+v^{1-p}.
\end{equation}
The test is justified by \(C^1\) approximation against the locally finite measure, and its value on the support of \(\mu_0\) is zero. The same source bound and stress argument apply to limits of classical positive solutions with coefficients \(a_j\to1\) locally uniformly and \(\nabla a_j\to0\) on interior compact sets. In that case the stress has \(a_j\) in its potential term and its divergence is \(2v_j^{1-p}\nabla a_j/(1-p)\). The critical scaling \(r^{-\gamma}v(r\,\cdot)\) and translations preserve the limiting identities and the gradient seminorm.

%%%%%%%%%%%%%%%%%%%%%%%%%%%%%%%%%%%%%%%%%%%%%%%%%%%%%%%%%%%%%%%%%%%%%%%%%%%%%%%%%%%%%%%%%%%%%%%%%%%%%%%%%%%%%%%%%%%%%%%%%%%%%%%%%%%%%%%%%%%%%%%%%%%%%%%%%%%%%%%%%%%%%%%%%%%%%%%%%%%%
\vskip0.12in
%%%%%%%%%%%%%%%%%%%%%%%%%%%%%%%%%%%%%%%%%%%%%%%%%%%%%%%%%%%%%%%%%%%%%%%%%%%%%%%%%%%%%%%%%%%%%%%%%%%%%%%%%%%%%%%%%%%%%%%%%%%%%%%%%%%%%%%%%%%%%%%%%%%%%%%%%%%%%%%%%%%%%%%%%%%%%%%%%%%%

For \(Q=(q,h)\), \(h\geq0\), put \(D_r=B_r(Q)\cap\mathbb H\), \(S_r=\partial B_r(Q)\cap\mathbb H\), and define
\begin{equation*}
 W_Q(r)=r^{-2\gamma}\int_{D_r}
        \left(|\nabla v|^2+\frac{2}{1-p}v^{1-p}\right)
        -\gamma r^{-2\gamma-1}\int_{S_r}v^2 .
\end{equation*}
The domain-variation calculation of Weiss~\cite{Weiss} gives the following identity with the boundary flux retained. For the interior formula for \(p>1\), see \cite[Theorem~3.6]{DavilaWangWei}.
\begin{equation}\label{eq:sharp-off-center}
 W_Q'(r)=2r^{-2\gamma}\int_{S_r}
      \left(\partial_r v-\frac{\gamma}{r}v\right)^2
      +h\,r^{-2\gamma-1}\int_{B_r(Q)\cap\partial\mathbb H}v_y^2 .
\end{equation}
To verify the flat-boundary term, write \(\mathcal V=v^{1-p}/(1-p)\), let \(E(r)\) denote the volume integral in \(W_Q\), and let \(J(r)=\int_{S_r}v^2\). The outer normal \(n\) on the flat boundary satisfies \(T_vn=-v_y^2n\) and \((x-Q)\cdot n=h\). Integrating the divergence of \(T_v(x-Q)\) gives
\begin{equation*}
 rE'(r)-2r\int_{S_r}(\partial_r v)^2
        -h\int_{B_r(Q)\cap\partial\mathbb H}v_y^2
       =4\int_{D_r}\mathcal V .
\end{equation*}
Equation \eqref{eq:sharp-product} gives \(E(r)=\int_{S_r}v\partial_r v+(1+p)\int_{D_r}\mathcal V\). The zero extension \(\widetilde{v^2}\) across the diameter is \(C^1\) because \(v=0\) there. Differentiating \(J(r)=r\int_{S^1}\widetilde{v^2}(Q+r\omega)\,\dd\omega\) gives \(J'(r)=r^{-1}J(r)+2\int_{S_r}v\partial_r v\). These identities give \eqref{eq:sharp-off-center} using \(\gamma(1+p)=2\). They hold for almost every radius and the integrated identity holds on compact intervals of radii. One can integrate on \(D_r\cap\{y>\varepsilon\}\) by smoothing the continuous vector fields with distributional divergence and then let \(\varepsilon\downarrow0\). The \(C^1\) traces ensure convergence of the fluxes. Interior balls satisfy the same formula without the flat term.

%%%%%%%%%%%%%%%%%%%%%%%%%%%%%%%%%%%%%%%%%%%%%%%%%%%%%%%%%%%%%%%%%%%%%%%%%%%%%%%%%%%%%%%%%%%%%%%%%%%%%%%%%%%%%%%%%%%%%%%%%%%%%%%%%%%%%%%%%%%%%%%%%%%%%%%%%%%%%%%%%%%%%%%%%%%%%%%%%%%%
\vskip0.12in
%%%%%%%%%%%%%%%%%%%%%%%%%%%%%%%%%%%%%%%%%%%%%%%%%%%%%%%%%%%%%%%%%%%%%%%%%%%%%%%%%%%%%%%%%%%%%%%%%%%%%%%%%%%%%%%%%%%%%%%%%%%%%%%%%%%%%%%%%%%%%%%%%%%%%%%%%%%%%%%%%%%%%%%%%%%%%%%%%%%%

We next record the homogeneous profiles appearing in this formula. This equation on the positive set is the Euler--Lagrange equation of the Alt--Phillips energy with potential \(2v^{1-p}/(1-p)\) \cite{AltPhillips}. The growth power \(2/(1+p)\) also occurs for regularized fully nonlinear singular absorption in \cite[author version, Theorem~4.2]{AraujoTeixeira}, with their power parameter \(1-p\). For the half-plane angular calculation we use the first integral from the proof of \cite[Theorem~1.12]{HuLongTaoWuZou}. In polar coordinates \(v=r^\gamma f(\theta)\). Each positive angular component satisfies
\begin{equation*}
 f''+\gamma^2 f=f^{-p},\qquad
 (f')^2+\gamma^2 f^2-\frac{2}{1-p}f^{1-p}=E_0 .
\end{equation*}
For a component with zero endpoints and maximum \(M\), its length is
\begin{equation*}
 \frac{2}{\gamma}\int_0^1
       [1-z^2-k(1-z^{1-p})]^{-1/2}\dd z,\qquad
 k=\frac{2 M^{-1-p}}{(1-p)\gamma^2}.
\end{equation*}
Here \(0<k\le1\) because \(E_0\ge0\) at a zero endpoint. The energy curve has one positive turning point. The length increases strictly with \(k\), from \(\pi/\gamma\) at \(k=0\) to \(\pi\) at \(k=1\), using the substitution \(z^{(1+p)/2}\) for the latter value. At a zero in the interior of the angular domain, nonnegativity and \(C^1\) regularity give \(f'=0\), hence \(k=1\). It follows that on the upper half-plane the only admissible homogeneous profile is \(Ay^\gamma\). On the whole plane, an angular profile with a zero must consist of two components of length \(\pi\), since it is positive almost everywhere. It gives \(A|x\cdot e|^\gamma\) for a unit vector \(e\).

%%%%%%%%%%%%%%%%%%%%%%%%%%%%%%%%%%%%%%%%%%%%%%%%%%%%%%%%%%%%%%%%%%%%%%%%%%%%%%%%%%%%%%%%%%%%%%%%%%%%%%%%%%%%%%%%%%%%%%%%%%%%%%%%%%%%%%%%%%%%%%%%%%%%%%%%%%%%%%%%%%%%%%%%%%%%%%%%%%%%
\vskip0.12in
%%%%%%%%%%%%%%%%%%%%%%%%%%%%%%%%%%%%%%%%%%%%%%%%%%%%%%%%%%%%%%%%%%%%%%%%%%%%%%%%%%%%%%%%%%%%%%%%%%%%%%%%%%%%%%%%%%%%%%%%%%%%%%%%%%%%%%%%%%%%%%%%%%%%%%%%%%%%%%%%%%%%%%%%%%%%%%%%%%%%

For completeness, a positive angular profile on the full circle must be constant. For related angular equations in the negative-potential-power regime, see \cite{DKV,DKVcorrection}. Here is a period argument for the positive profile in the present exponent range. A nonconstant positive orbit has two positive turning points and \(E_0<0\). With \(f=g^\gamma\), \(k_0=-E_0/\gamma^2>0\), and \(c_0=2/((1-p)\gamma^2)\), the first integral becomes
\begin{equation*}
 (g')^2+g^2+k_0g^{-2\alpha}=c_0 .
\end{equation*}
For \(z=g^2\) this gives
\begin{equation*}
 z''=2c_0-4z-2(1-\alpha)k_0z^{-\alpha},\qquad
 (z')''+\bigl(4-2\alpha(1-\alpha)k_0 z^{-1-\alpha}\bigr)z'=0 .
\end{equation*}
On an interval of length \(\ell\) between successive extrema,
\(z'\) vanishes at both endpoints and has fixed sign inside. Integration against \(z'\) followed by Poincar\'e's inequality gives
\((\pi/\ell)^2<4\). On the same interval \(f'\) solves the linear equation with coefficient \(\gamma^2+p f^{-1-p}>\gamma^2\). Sturm comparison gives \(\ell<\pi/\gamma\). Thus the least period lies strictly between \(\pi\) and \(2\pi/\gamma<2\pi\). It cannot divide \(2\pi\). The constant angular profile gives \(C_{\mathrm r}|x|^\gamma\), where \(C_{\mathrm r}^{1+p}=\gamma^{-2}\).

%%%%%%%%%%%%%%%%%%%%%%%%%%%%%%%%%%%%%%%%%%%%%%%%%%%%%%%%%%%%%%%%%%%%%%%%%%%%%%%%%%%%%%%%%%%%%%%%%%%%%%%%%%%%%%%%%%%%%%%%%%%%%%%%%%%%%%%%%%%%%%%%%%%%%%%%%%%%%%%%%%%%%%%%%%%%%%%%%%%%
\vskip0.12in
%%%%%%%%%%%%%%%%%%%%%%%%%%%%%%%%%%%%%%%%%%%%%%%%%%%%%%%%%%%%%%%%%%%%%%%%%%%%%%%%%%%%%%%%%%%%%%%%%%%%%%%%%%%%%%%%%%%%%%%%%%%%%%%%%%%%%%%%%%%%%%%%%%%%%%%%%%%%%%%%%%%%%%%%%%%%%%%%%%%%

\begin{lemma}\label{lem:sharp-interior-positivity}
A half-plane function in \eqref{eq:sharp-stationary-class} with the zero Dirichlet trace and a finite global \(\alpha\)-H\"older seminorm of the gradient is positive at every interior point. If its inward slope vanishes at a point on the boundary line, it equals \(Ay^\gamma\).
\end{lemma}

%%%%%%%%%%%%%%%%%%%%%%%%%%%%%%%%%%%%%%%%%%%%%%%%%%%%%%%%%%%%%%%%%%%%%%%%%%%%%%%%%%%%%%%%%%%%%%%%%%%%%%%%%%%%%%%%%%%%%%%%%%%%%%%%%%%%%%%%%%%%%%%%%%%%%%%%%%%%%%%%%%%%%%%%%%%%%%%%%%%%
\vskip0.12in
%%%%%%%%%%%%%%%%%%%%%%%%%%%%%%%%%%%%%%%%%%%%%%%%%%%%%%%%%%%%%%%%%%%%%%%%%%%%%%%%%%%%%%%%%%%%%%%%%%%%%%%%%%%%%%%%%%%%%%%%%%%%%%%%%%%%%%%%%%%%%%%%%%%%%%%%%%%%%%%%%%%%%%%%%%%%%%%%%%%%

\begin{proof}
Take the origin on the boundary line. The seminorm bound gives
\(|v(x)|\le C(|x|+|x|^\gamma)\) and \(|\nabla v(x)|\le C(1+|x|^\alpha)\). The function \(W_0\) is bounded for \(r\ge1\). Every sequence \(r_j\to\infty\) therefore has a subsequence whose rescalings converge locally in \(C^1\) and for which the integrated square in \eqref{eq:sharp-off-center} tends to zero on every rescaled annulus. Compactness follows from the gradient seminorm and the rescaled boundary value and slope at zero. Closure of \eqref{eq:sharp-stationary-class} makes the limit positive almost everywhere. The limit is homogeneous and hence is \(Ay^\gamma\). The energy at infinity is consequently the energy \(\Theta_{\partial}\) of this profile. The same limit holds for \(W_Q(r)\) with any fixed \(Q\), since in the dilated variables the centers differ by \(Q/r\), and local \(C^1\) convergence controls both integrals on the moving caps.

%%%%%%%%%%%%%%%%%%%%%%%%%%%%%%%%%%%%%%%%%%%%%%%%%%%%%%%%%%%%%%%%%%%%%%%%%%%%%%%%%%%%%%%%%%%%%%%%%%%%%%%%%%%%%%%%%%%%%%%%%%%%%%%%%%%%%%%%%%%%%%%%%%%%%%%%%%%%%%%%%%%%%%%%%%%%%%%%%%%%
\vskip0.12in
%%%%%%%%%%%%%%%%%%%%%%%%%%%%%%%%%%%%%%%%%%%%%%%%%%%%%%%%%%%%%%%%%%%%%%%%%%%%%%%%%%%%%%%%%%%%%%%%%%%%%%%%%%%%%%%%%%%%%%%%%%%%%%%%%%%%%%%%%%%%%%%%%%%%%%%%%%%%%%%%%%%%%%%%%%%%%%%%%%%%

If \(v(Q)=0\) at an interior point then \(\nabla v(Q)=0\), so the corresponding critical rescalings at \(Q\) are locally compact on the whole plane. For \(r<h\), the interior version of \eqref{eq:sharp-off-center} and the same growth bound at \(Q\) show that every vanishing-scale limit is homogeneous. Its energy is strictly greater than \(\Theta_{\partial}\). In fact, for a homogeneous profile, \eqref{eq:sharp-product} gives
\begin{equation*}
 \Theta=\frac{1}{2\gamma\alpha}\int f^{1-p}\dd\theta .
\end{equation*}
The line model has energy \(2\Theta_{\partial}\), and the radial model has ratio to \(\Theta_{\partial}\) equal to
\begin{equation*}
 \frac{2\pi(\alpha/(1+\alpha))^\alpha}
      {\displaystyle\int_0^\pi \sin^{2\alpha}\theta\dd\theta}>1.
\end{equation*}
The inequality follows from \((1+1/\alpha)^\alpha<2\) for \(0<\alpha<1\). This contradicts the monotonicity \(W_Q(0+)\le W_Q(\infty)\) from \eqref{eq:sharp-off-center}. Thus the excess measure in \eqref{eq:sharp-stationary-class} is zero. If a boundary slope vanishes at \(Q\), the small-scale argument now works with a half-plane center and gives \(W_Q(0+)=\Theta_{\partial}\). The two endpoint values of \(W_Q\) coincide, so the radial square vanishes at all scales. The half-plane homogeneous classification gives \(v=Ay^\gamma\).
\end{proof}

%%%%%%%%%%%%%%%%%%%%%%%%%%%%%%%%%%%%%%%%%%%%%%%%%%%%%%%%%%%%%%%%%%%%%%%%%%%%%%%%%%%%%%%%%%%%%%%%%%%%%%%%%%%%%%%%%%%%%%%%%%%%%%%%%%%%%%%%%%%%%%%%%%%%%%%%%%%%%%%%%%%%%%%%%%%%%%%%%%%%
\vskip0.32in
%%%%%%%%%%%%%%%%%%%%%%%%%%%%%%%%%%%%%%%%%%%%%%%%%%%%%%%%%%%%%%%%%%%%%%%%%%%%%%%%%%%%%%%%%%%%%%%%%%%%%%%%%%%%%%%%%%%%%%%%%%%%%%%%%%%%%%%%%%%%%%%%%%%%%%%%%%%%%%%%%%%%%%%%%%%%%%%%%%%%

\subsection{One-dimensionality of the transition profiles}

\begin{lemma}\label{lem:sharp-normal}
Every half-plane profile of \cref{lem:sharp-interior-positivity} satisfies \(v(x,y)=w_B(y)\), where \(B=v_y(0,0)\).
\end{lemma}

\begin{proof}
The zero-slope case is already covered by \cref{lem:sharp-interior-positivity}. Suppose that all boundary slopes are positive. We prove \(v_y>0\) by the rotating-line geometry of \cite{FarinaSciunzi}. Half-space monotonicity for regular nonlinearities is treated in \cite{BCNQualitative}. Here a Hardy estimate controls the singular coefficient of the reflected difference at the boundary line.

%%%%%%%%%%%%%%%%%%%%%%%%%%%%%%%%%%%%%%%%%%%%%%%%%%%%%%%%%%%%%%%%%%%%%%%%%%%%%%%%%%%%%%%%%%%%%%%%%%%%%%%%%%%%%%%%%%%%%%%%%%%%%%%%%%%%%%%%%%%%%%%%%%%%%%%%%%%%%%%%%%%%%%%%%%%%%%%%%%%%
\vskip0.12in
%%%%%%%%%%%%%%%%%%%%%%%%%%%%%%%%%%%%%%%%%%%%%%%%%%%%%%%%%%%%%%%%%%%%%%%%%%%%%%%%%%%%%%%%%%%%%%%%%%%%%%%%%%%%%%%%%%%%%%%%%%%%%%%%%%%%%%%%%%%%%%%%%%%%%%%%%%%%%%%%%%%%%%%%%%%%%%%%%%%%

On each bounded closed half-rectangle, strict positivity inside and positive inward slope give \(v(x,y)\ge c y\) for a rectangle-dependent \(c>0\). If two points are compared by reflection in a line of positive upward normal component and the reflected height is at least \(y\), their difference solves \(-\Delta w=q w\) with \(0\le q\le C y^{-1-p}\) on bounded reflection regions. A \(q\) with this upper bound has form norm less than one on a set of sufficiently small area. Indeed, for a function \(f\in H^1_0(D)\) extended by zero to the half-plane, Hardy's and Faber--Krahn's inequalities give
\begin{equation*}
 \int_D q f^2\le C\bigl(4h^{1-p}+C_{\mathrm F}|D|h^{-1-p}\bigr)
                        \int_D|\nabla f|^2 .
\end{equation*}
First choose \(h\) small, then \(|D|\) small. Testing with the positive part of \(w\) proves comparison when the boundary values of \(w\) are nonpositive.

%%%%%%%%%%%%%%%%%%%%%%%%%%%%%%%%%%%%%%%%%%%%%%%%%%%%%%%%%%%%%%%%%%%%%%%%%%%%%%%%%%%%%%%%%%%%%%%%%%%%%%%%%%%%%%%%%%%%%%%%%%%%%%%%%%%%%%%%%%%%%%%%%%%%%%%%%%%%%%%%%%%%%%%%%%%%%%%%%%%%
\vskip0.12in
%%%%%%%%%%%%%%%%%%%%%%%%%%%%%%%%%%%%%%%%%%%%%%%%%%%%%%%%%%%%%%%%%%%%%%%%%%%%%%%%%%%%%%%%%%%%%%%%%%%%%%%%%%%%%%%%%%%%%%%%%%%%%%%%%%%%%%%%%%%%%%%%%%%%%%%%%%%%%%%%%%%%%%%%%%%%%%%%%%%%

Fix \(x_0\in\mathbb R\) and reflect the triangle with vertical right side \(x=x_0\), base \(y=0\), and top line \(y=s+(x-x_0)\tan\vartheta\), \(0<\vartheta<\pi/2\), across its top line. The reflection sends \((x,y)\) to \((x,y)+2\ell(-\sin\vartheta,\cos\vartheta)\), with \(\ell=(s-y)\cos\vartheta+(x-x_0)\sin\vartheta\ge0\). The base has zero data and the top has equal data. On any compact path of parameters with \(\vartheta>0\), strict comparison persists as long as the ordering holds on the vertical side. To see this, retain a compact subset of the last triangle on which the inequality is strict, with complement of small area. The same ordering holds on the retained set after small perturbations. The comparison on the complement follows from the form bound. A first stopping point gives weak ordering by continuity and then strict ordering inside by the strong maximum principle applied to the nonnegative superharmonic difference. The strict boundary inequality on the interior of the base excludes equality throughout the triangle. All tests are legitimate by local \(C^1\) regularity, zero trace of the positive part, and the displayed form bound.

%%%%%%%%%%%%%%%%%%%%%%%%%%%%%%%%%%%%%%%%%%%%%%%%%%%%%%%%%%%%%%%%%%%%%%%%%%%%%%%%%%%%%%%%%%%%%%%%%%%%%%%%%%%%%%%%%%%%%%%%%%%%%%%%%%%%%%%%%%%%%%%%%%%%%%%%%%%%%%%%%%%%%%%%%%%%%%%%%%%%
\vskip0.12in
%%%%%%%%%%%%%%%%%%%%%%%%%%%%%%%%%%%%%%%%%%%%%%%%%%%%%%%%%%%%%%%%%%%%%%%%%%%%%%%%%%%%%%%%%%%%%%%%%%%%%%%%%%%%%%%%%%%%%%%%%%%%%%%%%%%%%%%%%%%%%%%%%%%%%%%%%%%%%%%%%%%%%%%%%%%%%%%%%%%%

At \((x_0,0)\) the gradient points strictly upwards. In a small neighborhood its dot product with \((-\sin\vartheta,\cos\vartheta)\) is positive for all small \(|\vartheta|\). Choose a fixed small \(\vartheta_0>0\) and \(s_0>0\) so that the initial triangle and its reflected points stay in this neighborhood for \(0<s\le s_0\). This starts strict comparison, and the vertical-side ordering holds for \(0\le\vartheta\le\vartheta_0\). Rotate along compact paths to every \(0<\vartheta\le\vartheta_0\). At each \(0<s\le s_0\), the triangles exhaust the left half-strip as \(\vartheta\downarrow0\), and the reflections converge to the horizontal reflection at height \(s\). The change of variable \(x\mapsto2x_0-x\) covers the right half-strip. Hence the horizontal comparison holds for every \(x\) and \(0<y<s\), \(0<s\le s_0\). It is strict in the strip by superharmonicity of the nonnegative reflected difference and its positive value at the bottom.

%%%%%%%%%%%%%%%%%%%%%%%%%%%%%%%%%%%%%%%%%%%%%%%%%%%%%%%%%%%%%%%%%%%%%%%%%%%%%%%%%%%%%%%%%%%%%%%%%%%%%%%%%%%%%%%%%%%%%%%%%%%%%%%%%%%%%%%%%%%%%%%%%%%%%%%%%%%%%%%%%%%%%%%%%%%%%%%%%%%%
\vskip0.12in
%%%%%%%%%%%%%%%%%%%%%%%%%%%%%%%%%%%%%%%%%%%%%%%%%%%%%%%%%%%%%%%%%%%%%%%%%%%%%%%%%%%%%%%%%%%%%%%%%%%%%%%%%%%%%%%%%%%%%%%%%%%%%%%%%%%%%%%%%%%%%%%%%%%%%%%%%%%%%%%%%%%%%%%%%%%%%%%%%%%%

Let \(S\) be the supremum of heights up to which the horizontal comparison holds at every level. If \(S\) were finite, the comparison would hold at \(S\) by continuity and would be strict inside. Hopf's lemma at the horizontal reflecting lines gives \(v_y(x,t)>0\) for \(0<t\le S\). The comparison on the fixed vertical side persists for \(s\le S+\varepsilon\) and \(|\vartheta|\le\varepsilon\) with some \(\varepsilon>0\). Otherwise, a sequence of failing vertical-side points \((x_0,y_j)\), with \(\vartheta_j\to0\) and \(s_j\to s\le S\), would contradict one of three facts. At \(s=0\) use positivity of the upward derivative near \((x_0,0)\). At \(y_j\to y<s\) use strict horizontal comparison or the zero trace if \(y=0\). At \(y_j\to s>0\) use \(v_y(x_0,s)>0\) along the shrinking reflection segment. Now slide a triangle of a fixed small \(\vartheta>0\) from the initial height to \(s\le S+\varepsilon\), then rotate to zero in the same limiting manner. Repeat on the right. The resulting comparison contradicts the definition of \(S\). Hopf's lemma at each horizontal line proves \(v_y>0\) in the half-plane.

%%%%%%%%%%%%%%%%%%%%%%%%%%%%%%%%%%%%%%%%%%%%%%%%%%%%%%%%%%%%%%%%%%%%%%%%%%%%%%%%%%%%%%%%%%%%%%%%%%%%%%%%%%%%%%%%%%%%%%%%%%%%%%%%%%%%%%%%%%%%%%%%%%%%%%%%%%%%%%%%%%%%%%%%%%%%%%%%%%%%
\vskip0.12in
%%%%%%%%%%%%%%%%%%%%%%%%%%%%%%%%%%%%%%%%%%%%%%%%%%%%%%%%%%%%%%%%%%%%%%%%%%%%%%%%%%%%%%%%%%%%%%%%%%%%%%%%%%%%%%%%%%%%%%%%%%%%%%%%%%%%%%%%%%%%%%%%%%%%%%%%%%%%%%%%%%%%%%%%%%%%%%%%%%%%

It remains to eliminate tangential dependence. Let \(M=[\nabla v]_{C^{0,\alpha}}\) and \(b(x)=v_y(x,0)\). For all \(x\) and \(y>0\),
\begin{equation}\label{eq:sharp-gradient-trace}
 |v_x(x,y)|\le M y^\alpha,\qquad
 |v_y(x,y)-b(x)|\le M y^\alpha .
\end{equation}
The estimate \(v_y(x,y)\ge c (b(x)+y^\alpha)\) holds with \(c>0\) depending on the uniform seminorm bound and \(p\). When \(b(x)/y^\alpha\ge2M+1\), this follows from \eqref{eq:sharp-gradient-trace}. Otherwise normalize at \((x,0)\) and scale by \(y\). A sequence violating the estimate has bounded normalized slopes and converges, after extraction, to a limit with derivative zero at \((0,1)\). Its gradient seminorm is bounded and the stationary class is closed. The preceding positivity and monotonicity proof, including the explicit zero-slope profile, gives a strictly positive derivative there, a contradiction.

%%%%%%%%%%%%%%%%%%%%%%%%%%%%%%%%%%%%%%%%%%%%%%%%%%%%%%%%%%%%%%%%%%%%%%%%%%%%%%%%%%%%%%%%%%%%%%%%%%%%%%%%%%%%%%%%%%%%%%%%%%%%%%%%%%%%%%%%%%%%%%%%%%%%%%%%%%%%%%%%%%%%%%%%%%%%%%%%%%%%
\vskip0.12in
%%%%%%%%%%%%%%%%%%%%%%%%%%%%%%%%%%%%%%%%%%%%%%%%%%%%%%%%%%%%%%%%%%%%%%%%%%%%%%%%%%%%%%%%%%%%%%%%%%%%%%%%%%%%%%%%%%%%%%%%%%%%%%%%%%%%%%%%%%%%%%%%%%%%%%%%%%%%%%%%%%%%%%%%%%%%%%%%%%%%

Thus \(q=v_x/v_y\) is bounded. Both derivatives solve the same linearized equation, so
\(\operatorname{div}(v_y^2\nabla q)=0\).
If \(\sup |q|>0\), choose a sequence \((x_j,y_j)\) at the supremum of one of its two signs, and translate and scale as above. The scaled slopes cannot tend to infinity, since \eqref{eq:sharp-gradient-trace} would make the quotient tend to zero. A compact limit has \(v_y>0\), so its quotient attains the signed supremum at an interior point. The strong maximum principle makes this quotient a nonzero constant, say \(c_1\). The equation \(v_x=c_1 v_y\) makes the limiting function constant on the straight characteristics with direction \((1,-c_1)\). Every such line through an interior point meets \(y=0\), which would make the limit zero there and in the interior. This contradicts positivity. Hence \(v_x=0\), and \eqref{eq:sharp-one-dimensional} identifies the profile.
\end{proof}

%%%%%%%%%%%%%%%%%%%%%%%%%%%%%%%%%%%%%%%%%%%%%%%%%%%%%%%%%%%%%%%%%%%%%%%%%%%%%%%%%%%%%%%%%%%%%%%%%%%%%%%%%%%%%%%%%%%%%%%%%%%%%%%%%%%%%%%%%%%%%%%%%%%%%%%%%%%%%%%%%%%%%%%%%%%%%%%%%%%%
\vskip0.22in
%%%%%%%%%%%%%%%%%%%%%%%%%%%%%%%%%%%%%%%%%%%%%%%%%%%%%%%%%%%%%%%%%%%%%%%%%%%%%%%%%%%%%%%%%%%%%%%%%%%%%%%%%%%%%%%%%%%%%%%%%%%%%%%%%%%%%%%%%%%%%%%%%%%%%%%%%%%%%%%%%%%%%%%%%%%%%%%%%%%%

\subsection{Returning to the boundary}

\begin{proof}[Proof of \cref{thm:sharp-boundary-transition}]
By \eqref{eq:negative-global-morrey}, localized on the boundary arc by the Whitney argument, the forcing \(\lambda_F v^{-p}\) has mass at most \(Cr^\gamma\) on each ball of radius \(r\) in a smaller conformal half-disc. Reflect the source oddly across the diameter and apply a symmetric cutoff. The estimate for the Newtonian gradient in \cref{lem:negative-interior-compact-index} gives an \(\alpha\)-H\"older gradient for its potential. The potential has zero trace by oddness. Subtract it from \(v\) and reflect the remaining harmonic function. This gives the uniform \(C^{1,\alpha}\) bound on the smaller closed half-disc, in particular existence and continuity of the nonnegative boundary slope.

Consider \(s_j\) in the smaller diameter, \(t_j\downarrow0\), and \(u_j\in\mathcal U\). Put
\begin{equation*}
 v_j(z)=\frac{u_j(F(s_j+t_j z))}
                 {\lambda_F(s_j)^{1/(1+p)} t_j^\gamma}
\end{equation*}
If \(B_{u_j}(s_j,t_j)\to B<\infty\), then \(v_j\) has a bounded slope at the origin. Critical scaling preserves the same gradient H\"older seminorm bound on each fixed half-disc, so the limiting seminorm is finite on the full half-plane. Its coefficient is \(\lambda_F(s_j+t_j z)/\lambda_F(s_j)\), which tends to one. The derivatives of the rescaled coefficients tend to zero on interior compact sets. In fact, \(F'\in C^{0,\eta}\) up to the diameter and the Cauchy estimate gives \(|F''(s+it)|\le Ct^{\eta-1}\) on smaller charts. The slope and seminorm bounds thus give a local \(C^1\) subsequential limit satisfying \eqref{eq:sharp-stationary-class}. By \cref{lem:sharp-interior-positivity,lem:sharp-normal} it is \(w_B\). Interior elliptic estimates yield local \(C^2\) convergence as well. In particular \(v_j(i)\to w_B(1)\).

%%%%%%%%%%%%%%%%%%%%%%%%%%%%%%%%%%%%%%%%%%%%%%%%%%%%%%%%%%%%%%%%%%%%%%%%%%%%%%%%%%%%%%%%%%%%%%%%%%%%%%%%%%%%%%%%%%%%%%%%%%%%%%%%%%%%%%%%%%%%%%%%%%%%%%%%%%%%%%%%%%%%%%%%%%%%%%%%%%%%
\vskip0.12in
%%%%%%%%%%%%%%%%%%%%%%%%%%%%%%%%%%%%%%%%%%%%%%%%%%%%%%%%%%%%%%%%%%%%%%%%%%%%%%%%%%%%%%%%%%%%%%%%%%%%%%%%%%%%%%%%%%%%%%%%%%%%%%%%%%%%%%%%%%%%%%%%%%%%%%%%%%%%%%%%%%%%%%%%%%%%%%%%%%%%

Uniformly on this arc,
\(d(F(s+it))/(t|F'(s+it)|)\to1\).
Indeed, the nearest boundary point lies at chart distance \(O(t)\) from \(s\). The first-order expansion of \(F\) reduces the minimization at that scale to the distance from \(i\) to the real line. Consequently \(d^2V_{u_j}\) tends to \(p w_B(1)^{-1-p}\). If instead \(B_{u_j}(s_j,t_j)\to\infty\), the boundary Taylor estimate gives \(v_j(i)=B_{u_j}(s_j,t_j)+O(1)\) and both sides of \eqref{eq:sharp-transition} tend to zero. Every failure sequence has a subsequence of one of these types, which proves uniformity.

%%%%%%%%%%%%%%%%%%%%%%%%%%%%%%%%%%%%%%%%%%%%%%%%%%%%%%%%%%%%%%%%%%%%%%%%%%%%%%%%%%%%%%%%%%%%%%%%%%%%%%%%%%%%%%%%%%%%%%%%%%%%%%%%%%%%%%%%%%%%%%%%%%%%%%%%%%%%%%%%%%%%%%%%%%%%%%%%%%%%
\vskip0.12in
%%%%%%%%%%%%%%%%%%%%%%%%%%%%%%%%%%%%%%%%%%%%%%%%%%%%%%%%%%%%%%%%%%%%%%%%%%%%%%%%%%%%%%%%%%%%%%%%%%%%%%%%%%%%%%%%%%%%%%%%%%%%%%%%%%%%%%%%%%%%%%%%%%%%%%%%%%%%%%%%%%%%%%%%%%%%%%%%%%%%

For \eqref{eq:sharp-transition-secant}, use the same test with both solutions. For bounded normalized slopes, extract a subsequence on which both slopes converge and apply the profile limit to each solution. There is a uniform positive lower bound for these values at \(i\), since a bounded-slope sequence has positive limits at least \(A\), and an unbounded-slope sequence diverges. If one scaled slope diverges, dominated convergence in \(\tau\) gives a zero limit of the secant coefficient. This gives the uniform secant law and in particular its all-path upper bound \(h_p+o(1)\).

%%%%%%%%%%%%%%%%%%%%%%%%%%%%%%%%%%%%%%%%%%%%%%%%%%%%%%%%%%%%%%%%%%%%%%%%%%%%%%%%%%%%%%%%%%%%%%%%%%%%%%%%%%%%%%%%%%%%%%%%%%%%%%%%%%%%%%%%%%%%%%%%%%%%%%%%%%%%%%%%%%%%%%%%%%%%%%%%%%%%
\vskip0.12in
%%%%%%%%%%%%%%%%%%%%%%%%%%%%%%%%%%%%%%%%%%%%%%%%%%%%%%%%%%%%%%%%%%%%%%%%%%%%%%%%%%%%%%%%%%%%%%%%%%%%%%%%%%%%%%%%%%%%%%%%%%%%%%%%%%%%%%%%%%%%%%%%%%%%%%%%%%%%%%%%%%%%%%%%%%%%%%%%%%%%

When at least one slope is positive at \(P\), its scaled slope tends to infinity as \(x\to P\), and \(d^2 V_{u,\widetilde u}\to0\). In the zero-slope case, the upper bound \(h_p+o(1)\), conformal invariance of Dirichlet energy, and the half-plane Hardy inequality give \(b_P\le4h_p\). Indeed the upper bound in a shrinking neighborhood pulls back to a potential bounded by \((h_p+o(1))/y^2\). The same comparison yields zero local form when one slope is positive. If both slopes vanish, blow up with fixed center \(F^{-1}(P)\). Both scaled solutions converge to \(Ay^\gamma\), and the rescaled secant potential converges to \(h_p/y^2\) on every compact subset of the open half-plane. Shrinking the support of any compactly supported half-plane test back to \(P\) gives the reverse inequality \(b_P\ge4h_p\), by the sharp half-plane Hardy inequality. The constant \(4\) follows from the sharp one-dimensional Dirichlet Hardy inequality \cite[Lemma~1.1]{BrezisMarcus} by integration on vertical lines. Dilating the horizontal cutoffs in product tests makes the horizontal energy negligible, which proves sharpness in the half-plane. The form estimates extend from compactly supported tests by density. This proves \eqref{eq:sharp-local-energy}. Along the fixed chart normal the pointwise limit is \(h_p\) on the zero-slope diagonal, and the excess form tends to zero by the all-path upper bound and Hardy's inequality.
\end{proof}

%%%%%%%%%%%%%%%%%%%%%%%%%%%%%%%%%%%%%%%%%%%%%%%%%%%%%%%%%%%%%%%%%%%%%%%%%%%%%%%%%%%%%%%%%%%%%%%%%%%%%%%%%%%%%%%%%%%%%%%%%%%%%%%%%%%%%%%%%%%%%%%%%%%%%%%%%%%%%%%%%%%%%%%%%%%%%%%%%%%%
\vskip0.12in
%%%%%%%%%%%%%%%%%%%%%%%%%%%%%%%%%%%%%%%%%%%%%%%%%%%%%%%%%%%%%%%%%%%%%%%%%%%%%%%%%%%%%%%%%%%%%%%%%%%%%%%%%%%%%%%%%%%%%%%%%%%%%%%%%%%%%%%%%%%%%%%%%%%%%%%%%%%%%%%%%%%%%%%%%%%%%%%%%%%%

\begin{remark}\label{rem:sharp-physical-slope}
On a \(C^2\) arc the transition law can be expressed using the physical nearest point \(\pi x\). The corresponding slope parameter is \(\partial_\nu u(\pi x)/(\lambda^{1/(1+p)}d(x)^\alpha)\). The real coordinate \(s'\) of this point satisfies \(|s'-s|=O(t)\) by local bi-Lipschitz continuity and the upper bound \(d(x)=O(t)\). The stationarity condition \(\operatorname{Re}((F(s+it)-F(s'))\overline{F'(s')})=0\) and the first-order \(C^{1,\eta}\) expansion then give \(s'-s=O(t^{1+\eta})\). The \(C^{0,\alpha}\) estimate for the chart slope and the distance comparison show that the two parameters differ by \(o(1)\) when bounded, and tend to infinity together when unbounded. Thus the law includes all paths based at nonzero slopes that approach the zero-slope set.
\end{remark}

%%%%%%%%%%%%%%%%%%%%%%%%%%%%%%%%%%%%%%%%%%%%%%%%%%%%%%%%%%%%%%%%%%%%%%%%%%%%%%%%%%%%%%%%%%%%%%%%%%%%%%%%%%%%%%%%%%%%%%%%%%%%%%%%%%%%%%%%%%%%%%%%%%%%%%%%%%%%%%%%%%%%%%%%%%%%%%%%%%%%
\vskip0.32in
%%%%%%%%%%%%%%%%%%%%%%%%%%%%%%%%%%%%%%%%%%%%%%%%%%%%%%%%%%%%%%%%%%%%%%%%%%%%%%%%%%%%%%%%%%%%%%%%%%%%%%%%%%%%%%%%%%%%%%%%%%%%%%%%%%%%%%%%%%%%%%%%%%%%%%%%%%%%%%%%%%%%%%%%%%%%%%%%%%%%

\section{\bf Martin and Dirichlet spectra}\label{sec:negative-spectral}

%%%%%%%%%%%%%%%%%%%%%%%%%%%%%%%%%%%%%%%%%%%%%%%%%%%%%%%%%%%%%%%%%%%%%%%%%%%%%%%%%%%%%%%%%%%%%%%%%%%%%%%%%%%%%%%%%%%%%%%%%%%%%%%%%%%%%%%%%%%%%%%%%%%%%%%%%%%%%%%%%%%%%%%%%%%%%%%%%%%%
\vskip0.12in
%%%%%%%%%%%%%%%%%%%%%%%%%%%%%%%%%%%%%%%%%%%%%%%%%%%%%%%%%%%%%%%%%%%%%%%%%%%%%%%%%%%%%%%%%%%%%%%%%%%%%%%%%%%%%%%%%%%%%%%%%%%%%%%%%%%%%%%%%%%%%%%%%%%%%%%%%%%%%%%%%%%%%%%%%%%%%%%%%%%%

\subsection{The energy characterization of the Martin essential radius}

Let \(E=H^1_0(\Omega)\) with the Dirichlet norm, and let \(X\) be the
Martin space for the single pole \(O\). For a nonnegative locally
bounded potential \(V\) set
\[
 \begin{split}
 \mathcal M_V&=\sup_\Omega G_\Omega[VK]/K,\\
 b(V)&=\sup_{\|\phi\|_E=1}\int_\Omega V|\phi|^2,\qquad
 b_\partial(V)=\lim_{\delta\downarrow0}
 \sup_{\|\phi\|_E=1}\int_{\{d<\delta\}}V|\phi|^2.
 \end{split}
\]
Here \(d(x)=\operatorname{dist}(x,\partial\Omega)\).
The energy operator \(B\) is defined by
\((B\phi,\psi)_E=\int V\phi\overline\psi\).
The operator on \(X\) is \(A z=G_\Omega[Vz]\).

%%%%%%%%%%%%%%%%%%%%%%%%%%%%%%%%%%%%%%%%%%%%%%%%%%%%%%%%%%%%%%%%%%%%%%%%%%%%%%%%%%%%%%%%%%%%%%%%%%%%%%%%%%%%%%%%%%%%%%%%%%%%%%%%%%%%%%%%%%%%%%%%%%%%%%%%%%%%%%%%%%%%%%%%%%%%%%%%%%%%
\vskip0.12in
%%%%%%%%%%%%%%%%%%%%%%%%%%%%%%%%%%%%%%%%%%%%%%%%%%%%%%%%%%%%%%%%%%%%%%%%%%%%%%%%%%%%%%%%%%%%%%%%%%%%%%%%%%%%%%%%%%%%%%%%%%%%%%%%%%%%%%%%%%%%%%%%%%%%%%%%%%%%%%%%%%%%%%%%%%%%%%%%%%%%

\begin{theorem}\label{thm:negative-two-space-radius}
On a bounded connected planar Lipschitz domain, assume
\(b(V)<\infty\). Then \(r_{\mathrm{ess},E}(B)=b_\partial(V)\).
If also \(\mathcal M_V<\infty\), then
\[
 r_X(A)\le b(V),\qquad
 r_{\mathrm{ess},X}(A)\le b_\partial(V).
\]
If in addition \(VK\in H^{-1}(\Omega)\), both inequalities are
equalities. In this case the spectra of \(A\) and \(B\) outside
the common closed essential disk agree, including the
multiplicities of the isolated eigenvalues.
\end{theorem}

%%%%%%%%%%%%%%%%%%%%%%%%%%%%%%%%%%%%%%%%%%%%%%%%%%%%%%%%%%%%%%%%%%%%%%%%%%%%%%%%%%%%%%%%%%%%%%%%%%%%%%%%%%%%%%%%%%%%%%%%%%%%%%%%%%%%%%%%%%%%%%%%%%%%%%%%%%%%%%%%%%%%%%%%%%%%%%%%%%%%
\vskip0.12in
%%%%%%%%%%%%%%%%%%%%%%%%%%%%%%%%%%%%%%%%%%%%%%%%%%%%%%%%%%%%%%%%%%%%%%%%%%%%%%%%%%%%%%%%%%%%%%%%%%%%%%%%%%%%%%%%%%%%%%%%%%%%%%%%%%%%%%%%%%%%%%%%%%%%%%%%%%%%%%%%%%%%%%%%%%%%%%%%%%%%

\begin{proof}
Write \(V_\delta=\boldsymbol 1_{\{d<\delta\}}V\). Then \(B-B_\delta\) is compact by
interior Rellich compactness, with \(B_\delta\) defined using
\(V_\delta\). Thus \(r_{\mathrm{ess},E}(B)\le b_\partial(V)\).
Conversely, choose \(\delta_j\downarrow0\) and unit vectors
\(\phi_j\rightharpoonup\phi\) in \(E\) with
\(\int V_{\delta_j}|\phi_j|^2\to b_\partial(V)\).
The tail form of \(\phi\) tends to zero by dominated convergence.
By Cauchy--Schwarz in this form,
\(\psi_j=\phi_j-\phi\) therefore has the same limiting tail form.
It is weakly null and \(\limsup\|\psi_j\|_E\le1\).
For every compact operator \(C\) on \(E\) we have
\(\langle(B-C)\psi_j,\psi_j\rangle_E
\ge\int V_{\delta_j}|\psi_j|^2-o(1)\) in real part.
This proves the reverse inequality for the essential norm, which
equals the essential spectral radius for self-adjoint \(B\).

%%%%%%%%%%%%%%%%%%%%%%%%%%%%%%%%%%%%%%%%%%%%%%%%%%%%%%%%%%%%%%%%%%%%%%%%%%%%%%%%%%%%%%%%%%%%%%%%%%%%%%%%%%%%%%%%%%%%%%%%%%%%%%%%%%%%%%%%%%%%%%%%%%%%%%%%%%%%%%%%%%%%%%%%%%%%%%%%%%%%
\vskip0.12in
%%%%%%%%%%%%%%%%%%%%%%%%%%%%%%%%%%%%%%%%%%%%%%%%%%%%%%%%%%%%%%%%%%%%%%%%%%%%%%%%%%%%%%%%%%%%%%%%%%%%%%%%%%%%%%%%%%%%%%%%%%%%%%%%%%%%%%%%%%%%%%%%%%%%%%%%%%%%%%%%%%%%%%%%%%%%%%%%%%%%

Assume now that the normalized Green row is finite.
The Green formula then defines \(A\) on \(X\).
The Green potentials are continuous in the interior, as follows
from the local boundedness of the density there and harmonicity of
the potential of the sources separated from any fixed interior ball.
Set \(A_\delta=G_\Omega[V_\delta\,\cdot]\).
Then \(A-A_\delta\) is compact on \(X\).
Indeed, for densities bounded and supported on a fixed interior compact
set, local regularity gives equicontinuity on a larger interior compact
set. Outside it the potentials are harmonic, and the maximum principle
with \(K\) bounds the global Martin norm by the norm on that larger
compact set. The potentials have zero continuous boundary values,
so that comparison remains valid at the pole.

%%%%%%%%%%%%%%%%%%%%%%%%%%%%%%%%%%%%%%%%%%%%%%%%%%%%%%%%%%%%%%%%%%%%%%%%%%%%%%%%%%%%%%%%%%%%%%%%%%%%%%%%%%%%%%%%%%%%%%%%%%%%%%%%%%%%%%%%%%%%%%%%%%%%%%%%%%%%%%%%%%%%%%%%%%%%%%%%%%%%
\vskip0.12in
%%%%%%%%%%%%%%%%%%%%%%%%%%%%%%%%%%%%%%%%%%%%%%%%%%%%%%%%%%%%%%%%%%%%%%%%%%%%%%%%%%%%%%%%%%%%%%%%%%%%%%%%%%%%%%%%%%%%%%%%%%%%%%%%%%%%%%%%%%%%%%%%%%%%%%%%%%%%%%%%%%%%%%%%%%%%%%%%%%%%

We use the pointwise Martin-kernel estimate of Frazier and
Verbitsky \cite[Theorem~2.5(A) and (1.8)--(1.9)]{FVUniformMartin}.
For a nonnegative potential \(W\) on a bounded uniform domain with
\(b(W)<1\), and for the Martin kernel \(K\) at a single fixed pole,
it gives
\begin{equation}\label{eq:negative-FV-martin-series}
 \sum_{j=0}^\infty (G_\Omega[W\,\cdot])^j K
 \le K \exp\!\bigl(C\,G_\Omega[WK]/K\bigr).
\end{equation}
Take the boundary representing measure to be \(\delta_O\). For the inclusion \(J_W:E\to L^2(W\dd x)\), the Green operator \(G_\Omega[W\,\cdot]\) on \(L^2(W\dd x)\) is \(J_WJ_W^*\). Indeed, \(J_W^* g\) solves \(-\Delta f=Wg\) weakly, and the Green formula follows by approximation with compactly supported data. Its norm is \(\|J_W\|^2=b(W)\). Bounded connected Lipschitz domains are uniform and are finitely connected in the plane, so the planar extension stated there applies. If \(t>0\) satisfies \(t b(V)<1\), apply
\eqref{eq:negative-FV-martin-series} to \(W=tV\). Positivity
gives \(\|A^j\|_X=\|A^jK\|_X\), so \(r_X(A)\le t^{-1}\).
Letting \(t\) tend to \(1/b(V)\) gives the first inequality,
including the zero case by allowing arbitrary \(t\).
For every fixed \(t\) with \(t b_\partial(V)<1\), the same
argument applies to \(tV_\delta\) for sufficiently small \(\delta\),
since \(\mathcal M_{V_\delta}\le\mathcal M_V\).
Compactness of \(A-A_\delta\) gives the essential-radius bound.

%%%%%%%%%%%%%%%%%%%%%%%%%%%%%%%%%%%%%%%%%%%%%%%%%%%%%%%%%%%%%%%%%%%%%%%%%%%%%%%%%%%%%%%%%%%%%%%%%%%%%%%%%%%%%%%%%%%%%%%%%%%%%%%%%%%%%%%%%%%%%%%%%%%%%%%%%%%%%%%%%%%%%%%%%%%%%%%%%%%%
\vskip0.12in
%%%%%%%%%%%%%%%%%%%%%%%%%%%%%%%%%%%%%%%%%%%%%%%%%%%%%%%%%%%%%%%%%%%%%%%%%%%%%%%%%%%%%%%%%%%%%%%%%%%%%%%%%%%%%%%%%%%%%%%%%%%%%%%%%%%%%%%%%%%%%%%%%%%%%%%%%%%%%%%%%%%%%%%%%%%%%%%%%%%%

Now assume \(VK\in H^{-1}\). Positivity of this datum implies
\(\|Vz\|_{H^{-1}}\le \|VK\|_{H^{-1}}\|z\|_X\), by testing
against the absolute value of an energy function. Thus the
same Green formula defines a bounded \(S:X\to E\).
Indeed, for a nonnegative density dominated by a multiple
of \(VK\), compact interior truncations have uniformly
bounded energy potentials. Their weak energy limit solves
the original equation by distributional convergence,
and monotone convergence of the kernels identifies this
limit with the Green potential. Signed or complex data
follow by decomposition. On \(X\cap E\), \(A\) and \(B\)
then agree by the energy identity.
For \(\phi\in C_c^\infty(\Omega)\), their resolvents
satisfy, for sufficiently large \(|z|\),
\begin{equation}\label{eq:negative-martin-energy-resolvent}
 (zI-B)^{-1}\phi
 =z^{-1}\phi+z^{-1}S(zI-A)^{-1}\phi .
\end{equation}
This follows by termwise comparison of the two convergent
resolvent series, since \(B^{j+1}\phi=S A^j\phi\).
The same power comparison shows that the spectral measure
of \(B\) at each such \(\phi\) is supported in
\([0,r_X(A)]\). Indeed
\(\|B^{j+1}\phi\|_E\le\|S\|\|A^j\|_X\|\phi\|_X\), and a
nonzero spectral component in \((r_X(A)+\varepsilon,\infty)\)
would violate this bound as \(j\to\infty\).
Density gives \(b(V)=\|B\|_E\le r_X(A)\).

%%%%%%%%%%%%%%%%%%%%%%%%%%%%%%%%%%%%%%%%%%%%%%%%%%%%%%%%%%%%%%%%%%%%%%%%%%%%%%%%%%%%%%%%%%%%%%%%%%%%%%%%%%%%%%%%%%%%%%%%%%%%%%%%%%%%%%%%%%%%%%%%%%%%%%%%%%%%%%%%%%%%%%%%%%%%%%%%%%%%
\vskip0.12in
%%%%%%%%%%%%%%%%%%%%%%%%%%%%%%%%%%%%%%%%%%%%%%%%%%%%%%%%%%%%%%%%%%%%%%%%%%%%%%%%%%%%%%%%%%%%%%%%%%%%%%%%%%%%%%%%%%%%%%%%%%%%%%%%%%%%%%%%%%%%%%%%%%%%%%%%%%%%%%%%%%%%%%%%%%%%%%%%%%%%

Write \(r_0=r_{\mathrm{ess},X}(A)\). The resolvent of \(A\)
on \(\{|z|>r_0\}\) is meromorphic with finite-rank principal
parts. This is the standard Fredholm resolvent consequence
of the essential spectral radius and invertibility near
infinity. Since \(B\) is self-adjoint, analytic continuation
gives \eqref{eq:negative-martin-energy-resolvent} on the upper
and lower half-planes outside the disk, away from poles of
the right-hand side. In particular, for \(\phi\in C_c^\infty\),
the energy scalar resolvent has a common analytic continuation
through every positive real interval above \(r_0\) disjoint
from the spectrum of \(A\). Its energy spectral measure
therefore vanishes on each such interval. Explicitly, the
imaginary part of the scalar resolvent at \(s+i\epsilon\)
is, up to sign, the convolution of the spectral measure
with \(\epsilon/((s-\cdot)^2+\epsilon^2)\). The common
analytic continuation from both half-planes forces its
limit to vanish uniformly on compact subintervals.
Integration in \(s\) gives the measure assertion.
Density of \(C_c^\infty(\Omega)\) then shows that
\(B\) has no spectrum on those intervals.

%%%%%%%%%%%%%%%%%%%%%%%%%%%%%%%%%%%%%%%%%%%%%%%%%%%%%%%%%%%%%%%%%%%%%%%%%%%%%%%%%%%%%%%%%%%%%%%%%%%%%%%%%%%%%%%%%%%%%%%%%%%%%%%%%%%%%%%%%%%%%%%%%%%%%%%%%%%%%%%%%%%%%%%%%%%%%%%%%%%%
\vskip0.12in
%%%%%%%%%%%%%%%%%%%%%%%%%%%%%%%%%%%%%%%%%%%%%%%%%%%%%%%%%%%%%%%%%%%%%%%%%%%%%%%%%%%%%%%%%%%%%%%%%%%%%%%%%%%%%%%%%%%%%%%%%%%%%%%%%%%%%%%%%%%%%%%%%%%%%%%%%%%%%%%%%%%%%%%%%%%%%%%%%%%%

Let \(s>r_0\) be an isolated spectral point of \(A\) on the positive real axis. A small contour about \(s\) avoids the spectra of both operators. By \eqref{eq:negative-martin-energy-resolvent}, the \(B\)-projection on this contour maps the dense tests into the image under \(S\) of the finite-dimensional \(A\)-spectral subspace. It therefore has finite rank, allowing rank zero if \(s\) is outside the spectrum of \(B\). Consequently \(r_{\mathrm{ess},E}(B)\le r_0\).
Finally every finite-dimensional \(A\)-spectral subspace
outside the common disk lies in \(E\), since \(A\) is
invertible there and maps \(X\) to \(E\).
Its restriction agrees with \(B\), so self-adjointness
excludes Jordan blocks and nonreal eigenvalues there.
The residue comparison then gives equality of the
eigenspaces and their multiplicities.
\end{proof}

%%%%%%%%%%%%%%%%%%%%%%%%%%%%%%%%%%%%%%%%%%%%%%%%%%%%%%%%%%%%%%%%%%%%%%%%%%%%%%%%%%%%%%%%%%%%%%%%%%%%%%%%%%%%%%%%%%%%%%%%%%%%%%%%%%%%%%%%%%%%%%%%%%%%%%%%%%%%%%%%%%%%%%%%%%%%%%%%%%%%
\vskip0.12in
%%%%%%%%%%%%%%%%%%%%%%%%%%%%%%%%%%%%%%%%%%%%%%%%%%%%%%%%%%%%%%%%%%%%%%%%%%%%%%%%%%%%%%%%%%%%%%%%%%%%%%%%%%%%%%%%%%%%%%%%%%%%%%%%%%%%%%%%%%%%%%%%%%%%%%%%%%%%%%%%%%%%%%%%%%%%%%%%%%%%

\begin{corollary}\label{cor:negative-essential-local-forms}
Assume \(b(V)<\infty\) and put
\[
 b_{\rm loc}(r)=\sup_{Q\in\partial\Omega}
 \ \sup_{\substack{0\ne\phi\in H^1_0(\Omega\cap B(Q,r))}}
 \frac{\displaystyle\int_{\Omega\cap B(Q,r)}V|\phi|^2}
      {\displaystyle\int_{\Omega\cap B(Q,r)}|\nabla\phi|^2}.
\]
Then \(b_\partial(V)=\lim_{r\downarrow0}b_{\rm loc}(r)\).
Consequently, if \(\mathcal M_V<\infty\) and \(VK\in H^{-1}\), the shrinking boundary
patch forms characterize the essential spectral radius on \(X\).
\end{corollary}

\begin{proof}
Zero extension of each test on \(\Omega\cap B(Q,r)\) to \(\Omega\)
gives \(b_{\rm loc}(r)\le b(V_r)\), since its support lies in
\(\{d<r\}\). For the other direction fix \(r>0\).
Cover \(\partial\Omega\) by finitely many balls \(B(Q_i,r/2)\).
Use a smooth quadratic partition of unity \(\sum_{i=0}^N\chi_i^2=1\)
on \(\overline\Omega\), with \(\chi_i\) supported in \(B(Q_i,r)\) for
\(i\ge1\), and \(\chi_0\) supported away from \(\partial\Omega\).
Such a partition is obtained by dividing smooth bump functions by
the square root of the sum of their squares.
Then \(\chi_i\phi\in H^1_0(\Omega\cap B(Q_i,r))\) for \(i\ge1\),
by approximating \(\phi\) by compactly supported smooth functions.
For \(\delta\) sufficiently small,
\[
 \begin{split}
 \int_{\{d<\delta\}}V|\phi|^2
 &\le b_{\rm loc}(r)\sum_{i=1}^N\|\chi_i\phi\|_E^2\\
 &\le b_{\rm loc}(r)
 \bigl(\|\phi\|_E^2+C_r\|\phi\|_{L^2}^2\bigr).
 \end{split}
\]
The last line follows by expanding the full sum of squared
gradients, since \(\sum_i\chi_i\nabla\chi_i=0\).
Apply the bound to a weakly null tail-maximizing sequence from
the proof of Theorem~\ref{thm:negative-two-space-radius} and use
Rellich compactness to remove the \(L^2\) term. This gives
\(b_\partial(V)\le b_{\rm loc}(r)\), and the limit yields equality.
\end{proof}

%%%%%%%%%%%%%%%%%%%%%%%%%%%%%%%%%%%%%%%%%%%%%%%%%%%%%%%%%%%%%%%%%%%%%%%%%%%%%%%%%%%%%%%%%%%%%%%%%%%%%%%%%%%%%%%%%%%%%%%%%%%%%%%%%%%%%%%%%%%%%%%%%%%%%%%%%%%%%%%%%%%%%%%%%%%%%%%%%%%%
\vskip0.12in
%%%%%%%%%%%%%%%%%%%%%%%%%%%%%%%%%%%%%%%%%%%%%%%%%%%%%%%%%%%%%%%%%%%%%%%%%%%%%%%%%%%%%%%%%%%%%%%%%%%%%%%%%%%%%%%%%%%%%%%%%%%%%%%%%%%%%%%%%%%%%%%%%%%%%%%%%%%%%%%%%%%%%%%%%%%%%%%%%%%%

\subsection{The critical value at a smooth zero slope}

If \(u\) has zero slope at an ordinary boundary point \(P\), then \cref{thm:sharp-boundary-transition} gives \(b_P(V_u)=\beta_p=1-(1-3p)^2/(1+p)^2\). For a stable \(u\) at \(p=1/3\), take shrinking-support tests \(\phi_j\) of unit energy with \((B\phi_j,\phi_j)_E\to1\). They converge weakly to zero and \(0\le B\le I\) gives \(\|(I-B)\phi_j\|_E^2\le ((I-B)\phi_j,\phi_j)_E\to0\). Thus \(1\) belongs to the essential spectrum and \(I-B\) is not Fredholm.

We use the following Green estimate in flat charts. Let \(D\subset\mathbb H\) contain the half-disc \(B^+_{R_0}(x_0)\) at \(x_0\in\mathbb R\). For \(z=x_0+iy\), \(w=s+iv\), and \(2y<r=|w-x_0|<cR_0\), one has
\begin{equation}\label{eq:negative-half-disc-green-lower}
 G_D(z,w)\ge c' yv/r^2
\end{equation}
when \(c\) is sufficiently small. The half-plane Green kernel has this lower bound by reflection. The difference between it and the half-disc kernel is harmonic in \(w\) with zero diameter data and cap values bounded by \(Cyv/R_0^2\), which bounds the correction in the half-disc. Domain monotonicity gives the stated estimate.

For every such zero slope, at any \(p\) and without stability,
\begin{equation}\label{eq:negative-smooth-unbounded-row}
 \sup_\Omega G_\Omega[V_u K]/K=\infty .
\end{equation}
In the conformal half-disc about \(P\), write \(k=K\circ F\asymp y\). The fixed-center profile gives
\[
 |F'(w)|^2 V_u(F(w))k(w)\gtrsim |w|^{-1}
\]
on a central angular sector. Restrict the Green kernel to the chart and use the conformal lower bound \eqref{eq:negative-half-disc-green-lower} at \(i\varepsilon\). On this sector for \(2\varepsilon<|w|<R\) with \(R\) fixed small, it is at least \(c\varepsilon/|w|\). Integration gives a lower bound \(c\log(R/(2\varepsilon))\) for the normalized potential. Thus the Martin-space linearization there is unbounded.

%%%%%%%%%%%%%%%%%%%%%%%%%%%%%%%%%%%%%%%%%%%%%%%%%%%%%%%%%%%%%%%%%%%%%%%%%%%%%%%%%%%%%%%%%%%%%%%%%%%%%%%%%%%%%%%%%%%%%%%%%%%%%%%%%%%%%%%%%%%%%%%%%%%%%%%%%%%%%%%%%%%%%%%%%%%%%%%%%%%%
\vskip0.32in
%%%%%%%%%%%%%%%%%%%%%%%%%%%%%%%%%%%%%%%%%%%%%%%%%%%%%%%%%%%%%%%%%%%%%%%%%%%%%%%%%%%%%%%%%%%%%%%%%%%%%%%%%%%%%%%%%%%%%%%%%%%%%%%%%%%%%%%%%%%%%%%%%%%%%%%%%%%%%%%%%%%%%%%%%%%%%%%%%%%%

\section{\bf Minimum mass on one-corner domains}
\label{sec:corner-bifurcation}

Throughout this section \(\Omega\) is a bounded connected planar Lipschitz domain with finitely many boundary components and \(C^2\) boundary away from \(O\). The two incident arcs are \(C^2\) up to \(O\) and meet at an opening \(\theta\in(0,2\pi)\). Write \(\beta_O=\pi/\theta\), and use the sector chart of \cite[Lemma~2.2]{HuLongTaoWuZou}. The sole pole is \(O\), with \(K(x_*)=1\).

%%%%%%%%%%%%%%%%%%%%%%%%%%%%%%%%%%%%%%%%%%%%%%%%%%%%%%%%%%%%%%%%%%%%%%%%%%%%%%%%%%%%%%%%%%%%%%%%%%%%%%%%%%%%%%%%%%%%%%%%%%%%%%%%%%%%%%%%%%%%%%%%%%%%%%%%%%%%%%%%%%%%%%%%%%%%%%%%%%%%
\vskip0.12in
%%%%%%%%%%%%%%%%%%%%%%%%%%%%%%%%%%%%%%%%%%%%%%%%%%%%%%%%%%%%%%%%%%%%%%%%%%%%%%%%%%%%%%%%%%%%%%%%%%%%%%%%%%%%%%%%%%%%%%%%%%%%%%%%%%%%%%%%%%%%%%%%%%%%%%%%%%%%%%%%%%%%%%%%%%%%%%%%%%%%

In this geometry, the admissible masses form the closed ray \([m_*,\infty)\) with \(m_*>0\). At each such mass the maximal solution \(U_m\) is the unique stable natural solution. The branch is pointwise strictly concave and satisfies
\begin{equation}\label{eq:negative-branch-order}
 U_{m'}\ge U_m+(m'-m)K\qquad (m'>m\ge m_*).
\end{equation}
The solution \(U_*=U_{m_*}\) at \(m_*\) is unique and all solutions with mass tending to \(m_*\) converge to it in \(C^2_{\mathrm{loc}}(\Omega)\). These are the one-corner branch and endpoint conclusions of \cite[Theorem~1.4]{HuLongTaoWuZou}.

%%%%%%%%%%%%%%%%%%%%%%%%%%%%%%%%%%%%%%%%%%%%%%%%%%%%%%%%%%%%%%%%%%%%%%%%%%%%%%%%%%%%%%%%%%%%%%%%%%%%%%%%%%%%%%%%%%%%%%%%%%%%%%%%%%%%%%%%%%%%%%%%%%%%%%%%%%%%%%%%%%%%%%%%%%%%%%%%%%%%
\vskip0.12in
%%%%%%%%%%%%%%%%%%%%%%%%%%%%%%%%%%%%%%%%%%%%%%%%%%%%%%%%%%%%%%%%%%%%%%%%%%%%%%%%%%%%%%%%%%%%%%%%%%%%%%%%%%%%%%%%%%%%%%%%%%%%%%%%%%%%%%%%%%%%%%%%%%%%%%%%%%%%%%%%%%%%%%%%%%%%%%%%%%%%

Use the Martin space \(X\) and Dirichlet space \(E\) of \cref{sec:negative-spectral}. For \(V_*=\lambda p U_*^{-1-p}\) let \(B_*\) be the associated energy operator. Stability gives \(0\le B_*\le I\).

\subsection{Compactness of all solutions at the minimum mass}

\begin{lemma}\label{lem:corner-compactness}
For \(0<\rho<1\), the ratio \(G_\Omega[d^{-\rho}]/K\) is bounded and tends to zero at \(O\). If \(u_j\) has mass \(m_j\to m_*\), then
\begin{equation}\label{eq:corner-all-convergence}
 u_j\longrightarrow U_*\ \text{in }X
      \text{ and in }C^2_{\mathrm{loc}}(\Omega),\qquad
 u_j-m_jK\longrightarrow U_*-m_*K\ \text{in }E.
\end{equation}
For masses in a compact subinterval of \((0,\infty)\), the convergence \(u/K-m\to0\) at \(O\) is uniform over all natural solutions.
\end{lemma}

%%%%%%%%%%%%%%%%%%%%%%%%%%%%%%%%%%%%%%%%%%%%%%%%%%%%%%%%%%%%%%%%%%%%%%%%%%%%%%%%%%%%%%%%%%%%%%%%%%%%%%%%%%%%%%%%%%%%%%%%%%%%%%%%%%%%%%%%%%%%%%%%%%%%%%%%%%%%%%%%%%%%%%%%%%%%%%%%%%%%
\vskip0.12in
%%%%%%%%%%%%%%%%%%%%%%%%%%%%%%%%%%%%%%%%%%%%%%%%%%%%%%%%%%%%%%%%%%%%%%%%%%%%%%%%%%%%%%%%%%%%%%%%%%%%%%%%%%%%%%%%%%%%%%%%%%%%%%%%%%%%%%%%%%%%%%%%%%%%%%%%%%%%%%%%%%%%%%%%%%%%%%%%%%%%

\begin{proof}
In the sector chart write \(\psi(\varphi)=\sin(\beta_O\varphi)\). The chart has bounded positive \(|F'|\) and, for small \(r\),
\(d(F(re^{i\varphi}))\asymp r\psi(\varphi)\) and \(K(F(re^{i\varphi}))\asymp r^{-\beta_O}\psi(\varphi)\). Fix \(0<s<\min(\beta_O,2-\rho)\). Let \(f\) solve
\(-f''-s^2 f=\psi^{-\rho}\) on \((0,\theta)\) with zero endpoints. Its Green kernel is positive since \(s<\beta_O\), and \(f\asymp\psi\) by the endpoint derivative formula and \(\rho<1\). On a sufficiently small truncated sector, a constant multiple of \(r^s f\) is a supersolution for the pullback of \(G_\Omega[d^{-\rho}]\). Its Laplacian with the negative sign is \(r^{s-2}\psi^{-\rho}\), dominating the transformed source. On the outer crosscut, the potential is \(O(\psi)\) up to the endpoints, by the ordinary boundary gradient bound for a source bounded by \(C d^{-\rho}\). This bound follows by reflecting the source in a regular conformal chart and using its local Morrey exponent \(2-\rho>1\), as in the proof of \cref{thm:sharp-boundary-transition}. The zero continuous trace of the potential follows from \(d^{-\rho}\in L^q(\Omega)\) for some \(q>1\) and the logarithmic Green bound. Maximum comparison therefore gives \(G_\Omega[d^{-\rho}]/K=O(r^{s+\beta_O})\) at the pole. Away from the pole the same gradient bound gives \(G_\Omega[d^{-\rho}]\lesssim d\) in ordinary boundary collars, where \(K\asymp d\) by Hopf's lemma. The ratio is bounded on interior compact sets as well. In particular \(K\gtrsim d\) globally, and \(G_\Omega[K^{-p}]/K\) is bounded by the estimate with \(\rho=p\).

%%%%%%%%%%%%%%%%%%%%%%%%%%%%%%%%%%%%%%%%%%%%%%%%%%%%%%%%%%%%%%%%%%%%%%%%%%%%%%%%%%%%%%%%%%%%%%%%%%%%%%%%%%%%%%%%%%%%%%%%%%%%%%%%%%%%%%%%%%%%%%%%%%%%%%%%%%%%%%%%%%%%%%%%%%%%%%%%%%%%
\vskip0.12in
%%%%%%%%%%%%%%%%%%%%%%%%%%%%%%%%%%%%%%%%%%%%%%%%%%%%%%%%%%%%%%%%%%%%%%%%%%%%%%%%%%%%%%%%%%%%%%%%%%%%%%%%%%%%%%%%%%%%%%%%%%%%%%%%%%%%%%%%%%%%%%%%%%%%%%%%%%%%%%%%%%%%%%%%%%%%%%%%%%%%

Set \(\sigma=2p/(1+p)\). Use \eqref{eq:negative-large-potential} for the full source. For its collar part, \eqref{eq:morrey-green-localized} gives the weight \(d^{-\sigma}\) on \(d<C\delta\), where \(d^{-\sigma}\le(C\delta)^{\rho-\sigma}d^{-\rho}\) for \(\sigma<\rho<1\). Thus uniformly over all natural solutions,
\begin{equation}\label{eq:corner-uniform-tail}
 G_\Omega[u^{-p}]\le C G_\Omega[d^{-\sigma}],\qquad
 \left\|G_\Omega[\boldsymbol1_{\{d<\delta\}}u^{-p}]\right\|_K
 \le C\delta^{\rho-\sigma}\|G_\Omega[d^{-\rho}]\|_K
 \quad(\sigma<\rho<1).
\end{equation}
The first inequality gives the uniform trace assertion. The local convergence to \(U_*\) follows from the endpoint estimate of \cite[Theorem~1.4]{HuLongTaoWuZou}. On each fixed \(\{d\ge\delta\}\), the corresponding sources converge uniformly. Their Green potentials converge in \(X\), since the Green potential of the characteristic function of this compact set is bounded by a constant times \(G_\Omega[d^{-\rho}]\). Cut off the remaining collar by \eqref{eq:corner-uniform-tail} and then send \(\delta\downarrow0\) to obtain convergence in \(X\). The uniform \(L^q\) bound for \(u^{-p}\) with some \(q>1\), proved in \cref{thm:negative-full-representation}, also gives convergence of the sources in \(L^{q'}\) for \(1<q'<q\). The Dirichlet Green operator maps \(L^{q'}\) continuously into \(E\) in two dimensions, giving the energy convergence.
\end{proof}

%%%%%%%%%%%%%%%%%%%%%%%%%%%%%%%%%%%%%%%%%%%%%%%%%%%%%%%%%%%%%%%%%%%%%%%%%%%%%%%%%%%%%%%%%%%%%%%%%%%%%%%%%%%%%%%%%%%%%%%%%%%%%%%%%%%%%%%%%%%%%%%%%%%%%%%%%%%%%%%%%%%%%%%%%%%%%%%%%%%%
\vskip0.12in
%%%%%%%%%%%%%%%%%%%%%%%%%%%%%%%%%%%%%%%%%%%%%%%%%%%%%%%%%%%%%%%%%%%%%%%%%%%%%%%%%%%%%%%%%%%%%%%%%%%%%%%%%%%%%%%%%%%%%%%%%%%%%%%%%%%%%%%%%%%%%%%%%%%%%%%%%%%%%%%%%%%%%%%%%%%%%%%%%%%%

\subsection{The fold with positive ordinary slopes}

\begin{lemma}\label{lem:negative-cone-fold}
Let \(U_*\) be a stable solution at the minimum mass with \(U_*\ge cK\) and \(\Lambda_p<\infty\). Set \(V_*=\lambda p U_*^{-1-p}\) and \(A_*=G_\Omega[V_*\,\cdot]\). If \(r_{\mathrm{ess},X}(A_*)<1\), then the \(X\) and \(E\) kernels of the linearized equation are generated by the same positive function \(\phi\in X\cap E\). Normalize \(\|\phi\|_K=1\). All solutions in a neighborhood of \((m_*,U_*)\) in \(\mathbb R\times X\) lie on an analytic curve
\[
 \begin{aligned}
 m(t)&=m_*+\tfrac12m''(0)t^2+O(t^3),\\
 u(t)&=U_*+t\phi+O_X(t^2),\\
 u(t)-m(t)K&=U_*-m_*K+t\phi+O_E(t^2),
 \end{aligned}
 \qquad
 m''(0)=\frac{\lambda p(p+1)\displaystyle\int_\Omega U_*^{-p-2}\phi^3}
                  {\displaystyle\int_\Omega V_*\phi K}>0 .
\]
\end{lemma}

\begin{proof}
By \cref{thm:negative-full-representation}, \(K^{-p}\le m_*^pU_*^{-p}\in L^s(\Omega)\) for some \(s>1\). The bound \(V_*K\lesssim K^{-p}\) gives \(A_*:X\to E\) and agreement with the energy operator \(B_*\) on \(X\cap E\). In a relative \(X\)-ball about \(U_*\), the map
\[
 \mathcal F(m,u)=u+\lambda G_\Omega[u^{-p}]-mK
\]
is analytic because its negative-power series is dominated by \(C K^{-p}\), with Green norm bounded by \(\Lambda_p\). The same domination gives analyticity of the Green term in \(E\), since \(L^s(\Omega)\subset H^{-1}(\Omega)\). The zeros are positive solutions of mass \(m\) by interior elliptic regularity, \(u\lesssim K\) at the ordinary boundary, and the Green representation.

%%%%%%%%%%%%%%%%%%%%%%%%%%%%%%%%%%%%%%%%%%%%%%%%%%%%%%%%%%%%%%%%%%%%%%%%%%%%%%%%%%%%%%%%%%%%%%%%%%%%%%%%%%%%%%%%%%%%%%%%%%%%%%%%%%%%%%%%%%%%%%%%%%%%%%%%%%%%%%%%%%%%%%%%%%%%%%%%%%%%
\vskip0.12in
%%%%%%%%%%%%%%%%%%%%%%%%%%%%%%%%%%%%%%%%%%%%%%%%%%%%%%%%%%%%%%%%%%%%%%%%%%%%%%%%%%%%%%%%%%%%%%%%%%%%%%%%%%%%%%%%%%%%%%%%%%%%%%%%%%%%%%%%%%%%%%%%%%%%%%%%%%%%%%%%%%%%%%%%%%%%%%%%%%%%

The operators \(I-tA_*\), \(0\le t\le1\), are Fredholm, so \(I-A_*\) has index zero. If this derivative of \(\mathcal F\) in \(u\) were invertible, the implicit function theorem would continue the solutions below \(m_*\). Thus \(A_*\) has a unit eigenvector, which belongs to \(E\). For any real \(f\in\ker_E(I-B_*)\), taking its absolute value preserves the zero value of the nonnegative endpoint form. By polarization \(|f|\) solves the linearized equation, and the interior minimum principle makes it strictly positive if \(f\ne0\). A linear combination of two independent kernel functions could vanish at an interior point, so this kernel is at most one-dimensional. This argument uses stability and local regularity alone.

%%%%%%%%%%%%%%%%%%%%%%%%%%%%%%%%%%%%%%%%%%%%%%%%%%%%%%%%%%%%%%%%%%%%%%%%%%%%%%%%%%%%%%%%%%%%%%%%%%%%%%%%%%%%%%%%%%%%%%%%%%%%%%%%%%%%%%%%%%%%%%%%%%%%%%%%%%%%%%%%%%%%%%%%%%%%%%%%%%%%
\vskip0.12in
%%%%%%%%%%%%%%%%%%%%%%%%%%%%%%%%%%%%%%%%%%%%%%%%%%%%%%%%%%%%%%%%%%%%%%%%%%%%%%%%%%%%%%%%%%%%%%%%%%%%%%%%%%%%%%%%%%%%%%%%%%%%%%%%%%%%%%%%%%%%%%%%%%%%%%%%%%%%%%%%%%%%%%%%%%%%%%%%%%%%

Define \(\ell(z)=\int_\Omega V_*\phi z=(\phi,A_*z)_E\) on \(X\). It annihilates the range of \(I-A_*\) by the energy equation for \(\phi\). Since the index is zero and \(\ell(\phi),\ell(K)>0\), the range is exactly \(\ker\ell\). Write \(u=U_*+t\phi+y\) with \(y\in\ker\ell\). The derivative in \((m,y)\) is the isomorphism \((s,z)\mapsto (I-A_*)z-sK\). The hypotheses of \cite[Theorem~3.2]{CrandallRabinowitz} hold because \(D_m\mathcal F=-K\notin\ker\ell=\operatorname{ran}(I-A_*)\). The resulting local solution curve has \(m'(0)=0\) and \(u'(0)=\phi\). It is analytic by the analytic inverse-function argument in the same splitting. Pairing the second derivative of the equation with \(\ell\) gives
\[
 m''(0)\ell(K)
  =\ell\big(\lambda p(p+1)G_\Omega[U_*^{-p-2}\phi^2]\big)
  =\lambda p(p+1)\int_\Omega U_*^{-p-2}\phi^3 .
\]
The last identity follows by the energy equation. The source is dominated by \(C K^{-p}\), which justifies the finite integral. Analyticity of \(-\lambda G_\Omega[u(t)^{-p}]\) in \(E\) and \(A_*\phi=\phi\) give the corrected energy expansion.
\end{proof}

%%%%%%%%%%%%%%%%%%%%%%%%%%%%%%%%%%%%%%%%%%%%%%%%%%%%%%%%%%%%%%%%%%%%%%%%%%%%%%%%%%%%%%%%%%%%%%%%%%%%%%%%%%%%%%%%%%%%%%%%%%%%%%%%%%%%%%%%%%%%%%%%%%%%%%%%%%%%%%%%%%%%%%%%%%%%%%%%%%%%
\vskip0.12in
%%%%%%%%%%%%%%%%%%%%%%%%%%%%%%%%%%%%%%%%%%%%%%%%%%%%%%%%%%%%%%%%%%%%%%%%%%%%%%%%%%%%%%%%%%%%%%%%%%%%%%%%%%%%%%%%%%%%%%%%%%%%%%%%%%%%%%%%%%%%%%%%%%%%%%%%%%%%%%%%%%%%%%%%%%%%%%%%%%%%

\begin{theorem}\label{thm:corner-fold}
Assume that \(\partial_\nu U_*(P)>0\) at every \(P\in\partial\Omega\setminus\{O\}\). Then \(B_*\) has a simple unit eigenvalue with eigenfunction \(\phi>0\) in \(E\cap X\). Normalize it by \(\|\phi\|_K=1\) and define
\begin{equation}\label{eq:corner-fold-coefficient}
 \kappa_{\mathrm f}=\frac{p+1}{2}\,
 \frac{\displaystyle\int_\Omega U_*^{-p-2}\phi^3}
      {\displaystyle\int_\Omega U_*^{-p-1}K\phi}>0 .
\end{equation}
For every sufficiently small \(\delta=m-m_*>0\), there are exactly two natural solutions, and
\begin{equation}\label{eq:corner-fold-expansion}
 u_m^\pm=U_*\pm\kappa_{\mathrm f}^{-1/2}\sqrt{\delta}\,\phi+O_X(\delta).
\end{equation}
The corrected difference satisfies
\((u_m^\pm-mK)-(U_*-m_*K)
=\pm\kappa_{\mathrm f}^{-1/2}\sqrt{\delta}\,\phi+O_E(\delta)\).
The larger solution is \(U_m\) and its linearization is coercive in \(E\).
The other solution has Morse index one.
\end{theorem}

\begin{proof}
At the pole \(U_*/K\to m_*>0\) by \cref{lem:corner-compactness}. Away from it the positive slopes and \(K\asymp d\) give \(U_*\ge cK\) by normal integration and interior compactness. Conversely this cone bound implies positive ordinary slopes. Since \(K\gtrsim d\), we have \(V_*\lesssim d^{-1-p}\) and \(V_*K\lesssim d^{-p}\). Thus the normalized Green row is finite by \cref{lem:corner-compactness}, and the collar form bound is \(O(\delta^{1-p})\) by Hardy's inequality. The energy operator is compact by interior Rellich compactness, and \cref{thm:negative-two-space-radius} gives \(r_{\mathrm{ess},X}(A_*)=0\).
%%%%%%%%%%%%%%%%%%%%%%%%%%%%%%%%%%%%%%%%%%%%%%%%%%%%%%%%%%%%%%%%%%%%%%%%%%%%%%%%%%%%%%%%%%%%%%%%%%%%%%%%%%%%%%%%%%%%%%%%%%%%%%%%%%%%%%%%%%%%%%%%%%%%%%%%%%%%%%%%%%%%%%%%%%%%%%%%%%%%
\vskip0.12in
%%%%%%%%%%%%%%%%%%%%%%%%%%%%%%%%%%%%%%%%%%%%%%%%%%%%%%%%%%%%%%%%%%%%%%%%%%%%%%%%%%%%%%%%%%%%%%%%%%%%%%%%%%%%%%%%%%%%%%%%%%%%%%%%%%%%%%%%%%%%%%%%%%%%%%%%%%%%%%%%%%%%%%%%%%%%%%%%%%%%

By \cref{lem:negative-cone-fold}, \(m''(0)/2\) is \(\kappa_{\mathrm f}\) from \eqref{eq:corner-fold-coefficient}. The derivative of \(m(t)\) has opposite strict signs on the two sides of zero, so \(t=\pm\kappa_{\mathrm f}^{-1/2}\sqrt{\delta}+O(\delta)\). This gives both expansions. The convergence \eqref{eq:corner-all-convergence} puts all natural solutions at close masses in the fold neighborhood.
%%%%%%%%%%%%%%%%%%%%%%%%%%%%%%%%%%%%%%%%%%%%%%%%%%%%%%%%%%%%%%%%%%%%%%%%%%%%%%%%%%%%%%%%%%%%%%%%%%%%%%%%%%%%%%%%%%%%%%%%%%%%%%%%%%%%%%%%%%%%%%%%%%%%%%%%%%%%%%%%%%%%%%%%%%%%%%%%%%%%
\vskip0.12in
%%%%%%%%%%%%%%%%%%%%%%%%%%%%%%%%%%%%%%%%%%%%%%%%%%%%%%%%%%%%%%%%%%%%%%%%%%%%%%%%%%%%%%%%%%%%%%%%%%%%%%%%%%%%%%%%%%%%%%%%%%%%%%%%%%%%%%%%%%%%%%%%%%%%%%%%%%%%%%%%%%%%%%%%%%%%%%%%%%%%

At a fixed interior point the minus solution is below \(U_*\) by \eqref{eq:corner-fold-expansion} and \(\phi>0\). Thus \(U_m\ge U_*+\delta K\) is the plus solution. Its potential is pointwise strictly below \(V_*\), so its nonnegative form is strictly positive on nonzero tests and coercive by compactness of the energy potential operator. The minus potential converges to \(V_*\) in form norm by the expansion and \(U_*\ge cK\). The endpoint gap on the energy orthogonal complement of \(\phi\) bounds the Morse index of the minus solution by one. It is unstable by the stable solution characterization in \cite[Theorem~1.4]{HuLongTaoWuZou}, so its index is one.
\end{proof}

%%%%%%%%%%%%%%%%%%%%%%%%%%%%%%%%%%%%%%%%%%%%%%%%%%%%%%%%%%%%%%%%%%%%%%%%%%%%%%%%%%%%%%%%%%%%%%%%%%%%%%%%%%%%%%%%%%%%%%%%%%%%%%%%%%%%%%%%%%%%%%%%%%%%%%%%%%%%%%%%%%%%%%%%%%%%%%%%%%%%
\vskip0.12in
%%%%%%%%%%%%%%%%%%%%%%%%%%%%%%%%%%%%%%%%%%%%%%%%%%%%%%%%%%%%%%%%%%%%%%%%%%%%%%%%%%%%%%%%%%%%%%%%%%%%%%%%%%%%%%%%%%%%%%%%%%%%%%%%%%%%%%%%%%%%%%%%%%%%%%%%%%%%%%%%%%%%%%%%%%%%%%%%%%%%

\subsection{Solution counts without a positive slope hypothesis}

Let \(u_j,\widetilde u_j\) have masses tending to \(m_*\), and write \(V_j=V_{u_j,\widetilde u_j}\). The following estimate measures the part of their forms that may concentrate at the boundary. By \cref{lem:corner-compactness}, \(V_j\to V_*\) locally uniformly. At the pole the same lemma gives a uniform lower bound \(u_j,\widetilde u_j\ge cK\), hence \(V_j\lesssim d^{-1-p}\) there. On the other boundary portions the solutions are uniformly bounded locally, by their Green representations. The uniform estimate of \cref{thm:sharp-boundary-transition} gives \(d^2V_j\le h_p+o(1)\) along these portions, using a finite chart cover away from the pole. Together these estimates give \(V_j\le C d^{-2}\) on \(\Omega\).

If \(z_j\rightharpoonup0\) in \(E\), then
\begin{equation}\label{eq:corner-weak-null}
 \limsup_j\int_\Omega V_j z_j^2\le
           \beta_p\limsup_j\|z_j\|_E^2 .
\end{equation}
For the proof, fix a small neighborhood of \(O\) where the \(d^{-1-p}\) bound holds. Its contribution in a collar \(d<\delta\) is at most \(C\delta^{1-p}\|z_j\|_E^2\) by the global Dirichlet Hardy inequality on the Lipschitz domain. On the rest of the boundary use a fixed smooth cutoff \(\chi\), equal to one in the smaller collar away from this neighborhood and supported away from \(O\). Its support in the normal direction can be taken in an arbitrarily thin tubular neighborhood of the ordinary arcs. One-dimensional Hardy along these normals bounds \(\int d^{-2}(\chi z_j)^2\) by \(4(1+\varepsilon)\|\chi z_j\|_E^2\), with \(\varepsilon\) arbitrarily small by the Jacobian bound. Since \(z_j\to0\) in \(L^2\), \(\|\chi z_j\|_E^2\le\|z_j\|_E^2+o(1)\) for \(0\le\chi\le1\). For \(\delta\) small, the secant estimate gives \(d^2V_j\le h_p+\varepsilon\) in this collar. The remaining fixed interior has a uniformly bounded potential and gives a vanishing integral. First let \(j\to\infty\) with these cutoffs fixed, then take the collar and Jacobian tolerances to zero to obtain \eqref{eq:corner-weak-null}.
%%%%%%%%%%%%%%%%%%%%%%%%%%%%%%%%%%%%%%%%%%%%%%%%%%%%%%%%%%%%%%%%%%%%%%%%%%%%%%%%%%%%%%%%%%%%%%%%%%%%%%%%%%%%%%%%%%%%%%%%%%%%%%%%%%%%%%%%%%%%%%%%%%%%%%%%%%%%%%%%%%%%%%%%%%%%%%%%%%%%
\vskip0.12in
%%%%%%%%%%%%%%%%%%%%%%%%%%%%%%%%%%%%%%%%%%%%%%%%%%%%%%%%%%%%%%%%%%%%%%%%%%%%%%%%%%%%%%%%%%%%%%%%%%%%%%%%%%%%%%%%%%%%%%%%%%%%%%%%%%%%%%%%%%%%%%%%%%%%%%%%%%%%%%%%%%%%%%%%%%%%%%%%%%%%

For \(f_j\rightharpoonup f\) in \(E\) with \(\|f_j\|_E=1\), Hardy domination and local convergence give \(\int V_j f^2\to\int V_* f^2\). The cross term \(\int V_j f(f_j-f)\) tends to zero by interior Rellich convergence and Cauchy--Schwarz on a collar, since the fixed function \(f/d\) is square integrable. Applying \eqref{eq:corner-weak-null} to \(f_j-f\) gives
\begin{equation}\label{eq:corner-form-compactness}
 \limsup_j\int_\Omega V_j f_j^2\le
    \int_\Omega V_* f^2+\beta_p(1-\|f\|_E^2).
\end{equation}

For \(Z_*=\{P\in\partial\Omega\setminus\{O\}:\partial_\nu U_*(P)=0\}\), we obtain
\begin{equation}\label{eq:corner-ess}
 r_{\mathrm{ess},E}(B_*)=
   \begin{cases}
      0,& Z_*=\varnothing,\\
      \beta_p,& Z_*\ne\varnothing .
   \end{cases}
\end{equation}
By \cref{thm:negative-two-space-radius,cor:negative-essential-local-forms}, the essential radius is the supremum of the pointwise local form limits. Indeed, for any sequence of shrinking patches \(B(Q_j,r_j)\) with \(Q_j\to Q\in\partial\Omega\), the patch \(B(Q,r_j+|Q_j-Q|)\) contains it. Compactness of the boundary then gives the upper bound by the supremum of the pointwise limits, and fixed-center tests give the lower bound. The local limit is zero at the pole by the bound \(V_*\lesssim d^{-1-p}\) and Hardy's inequality. At every ordinary point, \cref{thm:sharp-boundary-transition} gives zero at positive slope and \(\beta_p\) at zero slope. This proves \eqref{eq:corner-ess}.

%%%%%%%%%%%%%%%%%%%%%%%%%%%%%%%%%%%%%%%%%%%%%%%%%%%%%%%%%%%%%%%%%%%%%%%%%%%%%%%%%%%%%%%%%%%%%%%%%%%%%%%%%%%%%%%%%%%%%%%%%%%%%%%%%%%%%%%%%%%%%%%%%%%%%%%%%%%%%%%%%%%%%%%%%%%%%%%%%%%%
\vskip0.12in
%%%%%%%%%%%%%%%%%%%%%%%%%%%%%%%%%%%%%%%%%%%%%%%%%%%%%%%%%%%%%%%%%%%%%%%%%%%%%%%%%%%%%%%%%%%%%%%%%%%%%%%%%%%%%%%%%%%%%%%%%%%%%%%%%%%%%%%%%%%%%%%%%%%%%%%%%%%%%%%%%%%%%%%%%%%%%%%%%%%%

\begin{theorem}\label{thm:corner-count}
For \(p\ne1/3\), every mass sufficiently close to \(m_*\) from above has at most two natural solutions. If \(\ker_E(I-B_*)=\{0\}\), there is exactly one solution \(U_m\) at each such mass, and \(m\mapsto U_m-mK\) is right differentiable in \(E\) at \(m_*\). If \(\ker_E(I-B_*)\ne\{0\}\), then \(1\) is a simple isolated eigenvalue of \(B_*\). Whenever a second solution exists at a sufficiently close mass, it lies strictly below \(U_m\) in the interior and has Morse index one.
\end{theorem}

%%%%%%%%%%%%%%%%%%%%%%%%%%%%%%%%%%%%%%%%%%%%%%%%%%%%%%%%%%%%%%%%%%%%%%%%%%%%%%%%%%%%%%%%%%%%%%%%%%%%%%%%%%%%%%%%%%%%%%%%%%%%%%%%%%%%%%%%%%%%%%%%%%%%%%%%%%%%%%%%%%%%%%%%%%%%%%%%%%%%
\vskip0.12in
%%%%%%%%%%%%%%%%%%%%%%%%%%%%%%%%%%%%%%%%%%%%%%%%%%%%%%%%%%%%%%%%%%%%%%%%%%%%%%%%%%%%%%%%%%%%%%%%%%%%%%%%%%%%%%%%%%%%%%%%%%%%%%%%%%%%%%%%%%%%%%%%%%%%%%%%%%%%%%%%%%%%%%%%%%%%%%%%%%%%

\begin{proof}
Set \(M=\ker_E(I-B_*)\). The nonnegative-form argument of \cref{lem:negative-cone-fold} gives \(\dim M\le1\). Since \(p\ne1/3\), \(\beta_p<1\) and \eqref{eq:corner-ess} implies
\begin{equation}\label{eq:corner-core-gap}
 \int_\Omega V_* f^2\le q\|f\|_E^2\quad(f\in M^\perp)
 \qquad\text{with some }q<1 .
\end{equation}
For sequences of unit-energy \(f_j\in M^\perp\), \eqref{eq:corner-core-gap} bounds the right side of \eqref{eq:corner-form-compactness} by \(\max(q,\beta_p)<1\). This proves by contradiction a uniform coercivity bound on \(M^\perp\) for all secants and diagonals at sufficiently close masses.

%%%%%%%%%%%%%%%%%%%%%%%%%%%%%%%%%%%%%%%%%%%%%%%%%%%%%%%%%%%%%%%%%%%%%%%%%%%%%%%%%%%%%%%%%%%%%%%%%%%%%%%%%%%%%%%%%%%%%%%%%%%%%%%%%%%%%%%%%%%%%%%%%%%%%%%%%%%%%%%%%%%%%%%%%%%%%%%%%%%%
\vskip0.12in
%%%%%%%%%%%%%%%%%%%%%%%%%%%%%%%%%%%%%%%%%%%%%%%%%%%%%%%%%%%%%%%%%%%%%%%%%%%%%%%%%%%%%%%%%%%%%%%%%%%%%%%%%%%%%%%%%%%%%%%%%%%%%%%%%%%%%%%%%%%%%%%%%%%%%%%%%%%%%%%%%%%%%%%%%%%%%%%%%%%%

Now take two distinct solutions \(u,\widetilde u\) of the same mass and put \(w=u-\widetilde u\in E\). The harmonic terms cancel. Its equation is \(-\Delta w=V_{u,\widetilde u}w\), with the continuous form just estimated. If \(M=\{0\}\), testing with \(w\) contradicts uniform coercivity. Existence is given by the maximal branch. If \(\dim M=1\), any sign-changing \(w\) would supply a two-dimensional space on which the secant form is zero: testing with its positive and negative parts gives zero on each part, and the mixed pairing vanishes by disjoint supports. This space meets \(M^\perp\) nontrivially, contradicting the gap. Thus the solutions are ordered, strictly in the interior by the minimum principle for the nonnegative difference.

%%%%%%%%%%%%%%%%%%%%%%%%%%%%%%%%%%%%%%%%%%%%%%%%%%%%%%%%%%%%%%%%%%%%%%%%%%%%%%%%%%%%%%%%%%%%%%%%%%%%%%%%%%%%%%%%%%%%%%%%%%%%%%%%%%%%%%%%%%%%%%%%%%%%%%%%%%%%%%%%%%%%%%%%%%%%%%%%%%%%
\vskip0.12in
%%%%%%%%%%%%%%%%%%%%%%%%%%%%%%%%%%%%%%%%%%%%%%%%%%%%%%%%%%%%%%%%%%%%%%%%%%%%%%%%%%%%%%%%%%%%%%%%%%%%%%%%%%%%%%%%%%%%%%%%%%%%%%%%%%%%%%%%%%%%%%%%%%%%%%%%%%%%%%%%%%%%%%%%%%%%%%%%%%%%

If three solutions of the same mass existed, write them as \(u_1<u_2<u_3\) and set \(w_1=u_2-u_1\), \(w_2=u_3-u_2\). Testing each difference equation against the other gives
\begin{equation*}
 \int_\Omega (V_{u_1,u_2}-V_{u_2,u_3})w_1w_2=0 .
\end{equation*}
The integrand is strictly positive in the interior by the strict decrease of \(s^{-1-p}\). This contradiction gives the solution count. The diagonal bound gives index at most one for every nearby solution. Nonmaximal solutions are unstable by \cite[Theorem~1.4]{HuLongTaoWuZou}, so their index is one.

%%%%%%%%%%%%%%%%%%%%%%%%%%%%%%%%%%%%%%%%%%%%%%%%%%%%%%%%%%%%%%%%%%%%%%%%%%%%%%%%%%%%%%%%%%%%%%%%%%%%%%%%%%%%%%%%%%%%%%%%%%%%%%%%%%%%%%%%%%%%%%%%%%%%%%%%%%%%%%%%%%%%%%%%%%%%%%%%%%%%
\vskip0.12in
%%%%%%%%%%%%%%%%%%%%%%%%%%%%%%%%%%%%%%%%%%%%%%%%%%%%%%%%%%%%%%%%%%%%%%%%%%%%%%%%%%%%%%%%%%%%%%%%%%%%%%%%%%%%%%%%%%%%%%%%%%%%%%%%%%%%%%%%%%%%%%%%%%%%%%%%%%%%%%%%%%%%%%%%%%%%%%%%%%%%

When \(M=\{0\}\), \eqref{eq:corner-core-gap} gives \(\|B_*\|<1\). Thus the spectral integral in \cref{thm:negative-energy-endpoint} is finite, proving right differentiability of \(U_m-mK\) in \(E\).
\end{proof}

When two solutions occur along \(m_j\downarrow m_*\) with \(p\ne1/3\), their energy-normalized positive differences converge strongly in \(E\) to the positive unit-energy eigenfunction at one. In fact, the secant equation and \eqref{eq:corner-form-compactness} give \(1\le\|f\|_E^2+\beta_p(1-\|f\|_E^2)\) at a weak limit \(f\). Hence \(\|f\|_E=1\) and the nonnegative endpoint form vanishes at \(f\). Positivity and simplicity identify the limit for the whole sequence. When \(Z_*\ne\varnothing\), the existence of a second solution at a unit-kernel endpoint is open.

%%%%%%%%%%%%%%%%%%%%%%%%%%%%%%%%%%%%%%%%%%%%%%%%%%%%%%%%%%%%%%%%%%%%%%%%%%%%%%%%%%%%%%%%%%%%%%%%%%%%%%%%%%%%%%%%%%%%%%%%%%%%%%%%%%%%%%%%%%%%%%%%%%%%%%%%%%%%%%%%%%%%%%%%%%%%%%%%%%%%
\vskip0.12in
%%%%%%%%%%%%%%%%%%%%%%%%%%%%%%%%%%%%%%%%%%%%%%%%%%%%%%%%%%%%%%%%%%%%%%%%%%%%%%%%%%%%%%%%%%%%%%%%%%%%%%%%%%%%%%%%%%%%%%%%%%%%%%%%%%%%%%%%%%%%%%%%%%%%%%%%%%%%%%%%%%%%%%%%%%%%%%%%%%%%

\begin{remark}
If \(Z_*\ne\varnothing\) and \(p=1/3\), \eqref{eq:corner-ess} places one in the essential spectrum, so \(I-B_*\) has no coercive restriction to a finite-codimensional subspace. The solution count at such endpoints is open.

%%%%%%%%%%%%%%%%%%%%%%%%%%%%%%%%%%%%%%%%%%%%%%%%%%%%%%%%%%%%%%%%%%%%%%%%%%%%%%%%%%%%%%%%%%%%%%%%%%%%%%%%%%%%%%%%%%%%%%%%%%%%%%%%%%%%%%%%%%%%%%%%%%%%%%%%%%%%%%%%%%%%%%%%%%%%%%%%%%%%
\vskip0.12in
%%%%%%%%%%%%%%%%%%%%%%%%%%%%%%%%%%%%%%%%%%%%%%%%%%%%%%%%%%%%%%%%%%%%%%%%%%%%%%%%%%%%%%%%%%%%%%%%%%%%%%%%%%%%%%%%%%%%%%%%%%%%%%%%%%%%%%%%%%%%%%%%%%%%%%%%%%%%%%%%%%%%%%%%%%%%%%%%%%%%

At an ordinary zero slope, the divergence in \eqref{eq:negative-smooth-unbounded-row} prevents using the bounded Martin-space derivative from the fold proof.
\end{remark}
\subsection{Strict radial inclusion and a degenerate example}

\begin{proposition}\label{prop:corner-radial}
Suppose \(O=0\) and
\begin{equation}\label{eq:corner-radial-inclusion}
 t(\overline\Omega\setminus\{0\})\subset\Omega
 \qquad(0<t<1).
\end{equation}
There is a unique natural solution at every admissible mass. The solution map \(m\mapsto U_m\) has a right derivative at \(m_*\) in \(C^2_{\mathrm{loc}}(\Omega)\). For each nonempty compact \(D\Subset\Omega\) there are \(c_D,C_D>0\) with \(c_D(m-m_*)\le\|U_m-U_*\|_{C^0(D)}\le C_D(m-m_*)\) near \(m_*\). The endpoint has a zero ordinary slope and \(\ker_E(I-B_*)=\{0\}\). Consequently for \(p\ne1/3\) it has a strict energy gap. For \(p=1/3\), one belongs to its essential energy spectrum and is absent from its point spectrum.
\end{proposition}

\begin{proof}
We compare solutions by dilation about the pole. Since \(t\Omega\subset\Omega\), two-dimensional scaling and domain monotonicity give
\begin{equation*}
 G_\Omega(tx,ty)\ge G_{t\Omega}(tx,ty)=G_\Omega(x,y).
\end{equation*}
Along a fixed interior corner ray as \(y\to0\), the regular harmonic first mode of \(G_\Omega(x_*,\cdot)\) gives \(G_\Omega(x_*,ty)/G_\Omega(x_*,y)\to t^{\beta_O}\). This uses the corner expansion and the normalized chart derivative at \(O\) in \cite[Lemma~2.2]{HuLongTaoWuZou}. Divide the Green inequality by \(G_\Omega(x_*,y)\). The Martin limit of its left side is \(t^{\beta_O} K(tx)\), after inserting the intermediate denominator \(G_\Omega(x_*,ty)\). Hence \(K(tx)\ge t^{-\beta_O}K(x)\).
For any solution \(u\) of mass \(m_1\), its dilate \(u_t(x)=t^{-\gamma}u(tx)\) solves the equation on \(\Omega\). It dominates a solution \(v\) of mass \(m_2\) for small \(t\), by the uniform pole trace and \(K(tx)\ge t^{-\beta_O}K(x)\). Slide \(t\) as long as \(m_1t^{-\gamma-\beta_O}>m_2\). At a first stopping parameter the difference is nonnegative superharmonic and strictly positive in the interior. Strict radial inclusion gives positivity along the ordinary boundary. The strict mass inequality maintains comparison uniformly near the pole for nearby \(t\), and compactness maintains it on the remaining closed set. This excludes a stopping point. The comparison holds at equality by continuity. Taking \(t\uparrow1\) and interchanging solutions of the same mass proves uniqueness.

%%%%%%%%%%%%%%%%%%%%%%%%%%%%%%%%%%%%%%%%%%%%%%%%%%%%%%%%%%%%%%%%%%%%%%%%%%%%%%%%%%%%%%%%%%%%%%%%%%%%%%%%%%%%%%%%%%%%%%%%%%%%%%%%%%%%%%%%%%%%%%%%%%%%%%%%%%%%%%%%%%%%%%%%%%%%%%%%%%%%
\vskip0.12in
%%%%%%%%%%%%%%%%%%%%%%%%%%%%%%%%%%%%%%%%%%%%%%%%%%%%%%%%%%%%%%%%%%%%%%%%%%%%%%%%%%%%%%%%%%%%%%%%%%%%%%%%%%%%%%%%%%%%%%%%%%%%%%%%%%%%%%%%%%%%%%%%%%%%%%%%%%%%%%%%%%%%%%%%%%%%%%%%%%%%

With \(t_\delta=(m_*/(m_*+\delta))^{1/(\gamma+\beta_O)}\), compare \(U_{m_*+\delta}\) with the dilated endpoint. The branch order gives the lower bound \(\delta K\) on their increment from \(U_*\); Taylor expansion gives an upper bound \(C\delta\) on each interior compact set. The difference quotients are pointwise monotone by concavity and precompact in \(C^2\) on compact subsets by their linear secant equations. Their limit \(W\) therefore satisfies \(-\Delta W=V_*W\) and \(W\ge K\). Since \(K\gtrsim r^{-\beta_O}\) on an interior corner subsector, \(K/d\notin L^2\) and the Hardy inequality gives \(W\notin E\). The ground-state identity \cite[Lemma~2.4]{PinchoverTintarev} implies
\[
 \int_\Omega (|\nabla\zeta|^2-V_*\zeta^2)
   =\int_\Omega W^2|\nabla(\zeta/W)|^2
 \qquad(\zeta\in C_c^\infty(\Omega)).
\]
Approximation of an energy unit eigenmode by compact tests and lower semicontinuity on interior compacts would make that mode a multiple of \(W\). This is impossible in \(E\) for a nonzero mode. If every ordinary slope were positive, \cref{thm:corner-fold} would give such a mode. Thus \(Z_*\ne\varnothing\), and \eqref{eq:corner-ess} yields the remaining spectral assertions.

\end{proof}

%%%%%%%%%%%%%%%%%%%%%%%%%%%%%%%%%%%%%%%%%%%%%%%%%%%%%%%%%%%%%%%%%%%%%%%%%%%%%%%%%%%%%%%%%%%%%%%%%%%%%%%%%%%%%%%%%%%%%%%%%%%%%%%%%%%%%%%%%%%%%%%%%%%%%%%%%%%%%%%%%%%%%%%%%%%%%%%%%%%%
\vskip0.12in
%%%%%%%%%%%%%%%%%%%%%%%%%%%%%%%%%%%%%%%%%%%%%%%%%%%%%%%%%%%%%%%%%%%%%%%%%%%%%%%%%%%%%%%%%%%%%%%%%%%%%%%%%%%%%%%%%%%%%%%%%%%%%%%%%%%%%%%%%%%%%%%%%%%%%%%%%%%%%%%%%%%%%%%%%%%%%%%%%%%%

For example, for \(a_0>0\) the domain
\begin{equation*}
 \Omega_{a_0}=\{(x,y):0<x<1,\ |y|<a_0x\sqrt{1-x}\}
\end{equation*}
has its only corner at the origin with angle \(2\arctan a_0\). At \((1,0)\) the two graphs meet on the smooth level set \(y^2-a_0^2x^2(1-x)=0\), whose \(x\)-derivative is nonzero there. For \(0<t<1\) and \(0<x\le1\), the inequality \(t a_0x\sqrt{1-x}<a_0tx\sqrt{1-tx}\) proves \eqref{eq:corner-radial-inclusion}. It therefore realizes the degenerate, same-mass unique case of \cref{prop:corner-radial} on an ordinary one-corner domain.

%%%%%%%%%%%%%%%%%%%%%%%%%%%%%%%%%%%%%%%%%%%%%%%%%%%%%%%%%%%%%%%%%%%%%%%%%%%%%%%%%%%%%%%%%%%%%%%%%%%%%%%%%%%%%%%%%%%%%%%%%%%%%%%%%%%%%%%%%%%%%%%%%%%%%%%%%%%%%%%%%%%%%%%%%%%%%%%%%%%%
\vskip0.32in
%%%%%%%%%%%%%%%%%%%%%%%%%%%%%%%%%%%%%%%%%%%%%%%%%%%%%%%%%%%%%%%%%%%%%%%%%%%%%%%%%%%%%%%%%%%%%%%%%%%%%%%%%%%%%%%%%%%%%%%%%%%%%%%%%%%%%%%%%%%%%%%%%%%%%%%%%%%%%%%%%%%%%%%%%%%%%%%%%%%%

\section{\bf Flat-tangent domains and normalized Green bounds}\label{sec:flat-geometry}

Fix a smooth bounded Jordan domain \(D\subset\mathbb H\) whose boundary includes \([-1,1]\) and which contains \(B^+_{3/4}(0)\). Fix \(z_*\in D\). Set \(a=1+p\), \(\beta=1-p\), and \(\xi=\beta/2\). Choose positive numbers \(x_{n+1}\le x_n/2\), with \(x_1\) small, and a sufficiently small fixed \(\epsilon>0\). Set
\[
 \rho_n=\epsilon x_n^2,\qquad a_n=\rho_n^\beta,\qquad
 b_n=\rho_n\exp(-\mu_n/a_n),
\]
where the depths \(\mu_n\) are in a fixed compact subinterval of \((0,\infty)\). Let \(E_n=B_{\rho_n}(x_n)\cap\mathbb H\). We use the analytic product \(F'\) below and its primitive with \(F(0)=0\). The truncated maps \(F_N\) use the first \(N\) factors with the same normalization. The domains are \(\Omega=F(D)\) and \(\Omega_N=F_N(D)\). In the preimage the Martin kernel is the same \(k\) at \(-1/2\) for all these domains, normalized by \(k(z_*)=1\). We work with \(v=u\circ F\), \(J=|F'|^2\); for a prefix replace \(F,J\) by \(F_N,J_N=|F_N'|^2\). Conformal invariance gives \(v+\lambda G_D[Jv^{-p}]=mk\).

\begin{theorem}\label{thm:negative-critical-angular-counterexample}
The domain \(\Omega\) is Lipschitz and smooth except at \(Q=F(0)\), where its boundary has a tangent line. For the pole \(O=F(-1/2)\),
\[
 \Lambda_p(K)<\infty,\qquad
 \liminf_{\delta\downarrow0}\big\|G_\Omega[\boldsymbol1_{\{d<\delta\}}K^{-p}]\big\|_K>0,\qquad
 \sup_\Omega P_\sigma/K=\infty .
\]
The localized potentials with \(u^{-p}\) in place of \(K^{-p}\) also fail to vanish in the Martin norm for every positive solution. The bound on \(\Lambda_p\) is uniform for the prefixes.
\end{theorem}

\begin{proof}

Define a conformal chart by prescribing its derivative
\[
 F'(z)=\prod_{n\ge1}
 \left(\frac{z-x_n+i\rho_n}{z-x_n+ib_n}\right)^\xi .
\]
Use the logarithm of each ratio with argument in
\((-\pi/2,\pi/2)\). In fact the ratio has real part greater than
one in the upper half-plane and modulus at least one. At every \(z\)
at most one disk \(B_{x_n/8}(x_n)\) contains \(z\). The logarithms of all
other factors have absolute sum bounded by
\(C\sum_n\rho_n/x_n\), by expanding
\(\log(1+i(\rho_n-b_n)/(z-x_n+ib_n))\) there. This proves local
holomorphic convergence and gives \(|\arg F'|\le\xi(\pi/2+C\epsilon)
<\pi/2\). Also \(|F'|\ge1\) and \(\operatorname{Re}F'\ge c>0\).

%%%%%%%%%%%%%%%%%%%%%%%%%%%%%%%%%%%%%%%%%%%%%%%%%%%%%%%%%%%%%%%%%%%%%%%%%%%%%%%%%%%%%%%%%%%%%%%%%%%%%%%%%%%%%%%%%%%%%%%%%%%%%%%%%%%%%%%%%%%%%%%%%%%%%%%%%%%%%%%%%%%%%%%%%%%%%%%%%%%%
\vskip0.12in
%%%%%%%%%%%%%%%%%%%%%%%%%%%%%%%%%%%%%%%%%%%%%%%%%%%%%%%%%%%%%%%%%%%%%%%%%%%%%%%%%%%%%%%%%%%%%%%%%%%%%%%%%%%%%%%%%%%%%%%%%%%%%%%%%%%%%%%%%%%%%%%%%%%%%%%%%%%%%%%%%%%%%%%%%%%%%%%%%%%%

Let \(F\) be a primitive. Integrating on line segments gives
\(|F(z)-F(w)|\ge c|z-w|\), and thus injectivity. The factor bounds
give \(|F'(z)|\le C(1+(\operatorname{Im}z)^{-\xi})\). Integration
on vertical segments and on a horizontal segment at height
\(\max\{\operatorname{Im}z,\operatorname{Im}w,|z-w|\}\)
gives a continuous extension to the closed half-plane, and the
lower distance bound persists there. Normalize \(F(0)=0\).
The argument bound shows that \(F(\mathbb R)\) is a Lipschitz graph
over the entire real axis. The lower distance bound makes the
map proper on the closed half-plane, so the open image of the
upper half-plane is a full component of the graph complement.
At a real point away from zero the derivative extends
holomorphically and has positive real part. Thus the image lies
above the graph, by the orientation of the boundary image.

%%%%%%%%%%%%%%%%%%%%%%%%%%%%%%%%%%%%%%%%%%%%%%%%%%%%%%%%%%%%%%%%%%%%%%%%%%%%%%%%%%%%%%%%%%%%%%%%%%%%%%%%%%%%%%%%%%%%%%%%%%%%%%%%%%%%%%%%%%%%%%%%%%%%%%%%%%%%%%%%%%%%%%%%%%%%%%%%%%%%
\vskip0.12in
%%%%%%%%%%%%%%%%%%%%%%%%%%%%%%%%%%%%%%%%%%%%%%%%%%%%%%%%%%%%%%%%%%%%%%%%%%%%%%%%%%%%%%%%%%%%%%%%%%%%%%%%%%%%%%%%%%%%%%%%%%%%%%%%%%%%%%%%%%%%%%%%%%%%%%%%%%%%%%%%%%%%%%%%%%%%%%%%%%%%

Restrict \(F\) to the fixed cap \(D\) and identify the harmonic normalization by \(x_*=F(z_*)\).
Away from zero, local normal convergence of the logarithm series
across the real line follows by omitting a finite number of factors
and using \(\sum_n\rho_n<\infty\) on compact sets separated from zero.
The finitely many factors have no singularities on the real line.
Hence the image boundary is smooth off \(Q\), and it is Lipschitz
there by the graph construction.

%%%%%%%%%%%%%%%%%%%%%%%%%%%%%%%%%%%%%%%%%%%%%%%%%%%%%%%%%%%%%%%%%%%%%%%%%%%%%%%%%%%%%%%%%%%%%%%%%%%%%%%%%%%%%%%%%%%%%%%%%%%%%%%%%%%%%%%%%%%%%%%%%%%%%%%%%%%%%%%%%%%%%%%%%%%%%%%%%%%%
\vskip0.12in
%%%%%%%%%%%%%%%%%%%%%%%%%%%%%%%%%%%%%%%%%%%%%%%%%%%%%%%%%%%%%%%%%%%%%%%%%%%%%%%%%%%%%%%%%%%%%%%%%%%%%%%%%%%%%%%%%%%%%%%%%%%%%%%%%%%%%%%%%%%%%%%%%%%%%%%%%%%%%%%%%%%%%%%%%%%%%%%%%%%%

The exceptional point has a tangent line. Indeed the nonzero
product
\[
 C_0=\prod_n\left(\frac{-x_n+i\rho_n}{-x_n+ib_n}\right)^\xi
\]
converges. On the real axis, exclude the intervals
\(|x-x_n|<\sqrt{\rho_n x_n}\). Then \(F'(x)\to C_0\) as
\(x\to0\) outside the intervals: the fixed factors converge to
their values at zero, the nonexceptional tail logarithms are
uniformly summable, and any remaining logarithm is
\(O(\sqrt{\rho_n/x_n})\). On each excluded interval the integral
of \(|F'(x)|\) is \(O(\sqrt{\rho_n x_n})\), since
\(\int_0^{\rho_n}(\rho_n/(s+b_n))^\xi\dd s=O(\rho_n)\).
The sum of these bounds for intervals meeting \([-r,r]\) is
\(O(r^{3/2})=o(r)\). The same integrable power bounds dominate
\(F'(x+iy)\) on horizontal lines down to the boundary, so
integration on the real axis gives \(F(x)=C_0 x+o(|x|)\).
Thus, after rotation, the tangent wedge at \(Q\) is a half-plane.

%%%%%%%%%%%%%%%%%%%%%%%%%%%%%%%%%%%%%%%%%%%%%%%%%%%%%%%%%%%%%%%%%%%%%%%%%%%%%%%%%%%%%%%%%%%%%%%%%%%%%%%%%%%%%%%%%%%%%%%%%%%%%%%%%%%%%%%%%%%%%%%%%%%%%%%%%%%%%%%%%%%%%%%%%%%%%%%%%%%%
\vskip0.12in
%%%%%%%%%%%%%%%%%%%%%%%%%%%%%%%%%%%%%%%%%%%%%%%%%%%%%%%%%%%%%%%%%%%%%%%%%%%%%%%%%%%%%%%%%%%%%%%%%%%%%%%%%%%%%%%%%%%%%%%%%%%%%%%%%%%%%%%%%%%%%%%%%%%%%%%%%%%%%%%%%%%%%%%%%%%%%%%%%%%%

Write \(k=K\circ F\), \(z=x+iy\), \(w=s+iv\), and
\(E_n=B_{\rho_n}(x_n)\cap\mathbb H\).
The function \(k\) is a positive multiple of the Martin kernel of \(D\)
at \(-1/2\). This follows either from the Poisson representation
in \(D\) or from conformal invariance of the positive harmonic
pole class. By reflection and the Hopf lemma,
\(k(w)\asymp v\) on all \(E_n\), uniformly in \(n\).
Also \(k(z)\ge cy\) on \(B_{2/3}(0)\cap\mathbb H\).
For the latter estimate one can take the harmonic comparison with
\(y\) on a larger upper half-disc, using the boundary ratio on
its cap and the nonnegative pole.

%%%%%%%%%%%%%%%%%%%%%%%%%%%%%%%%%%%%%%%%%%%%%%%%%%%%%%%%%%%%%%%%%%%%%%%%%%%%%%%%%%%%%%%%%%%%%%%%%%%%%%%%%%%%%%%%%%%%%%%%%%%%%%%%%%%%%%%%%%%%%%%%%%%%%%%%%%%%%%%%%%%%%%%%%%%%%%%%%%%%
\vskip0.12in
%%%%%%%%%%%%%%%%%%%%%%%%%%%%%%%%%%%%%%%%%%%%%%%%%%%%%%%%%%%%%%%%%%%%%%%%%%%%%%%%%%%%%%%%%%%%%%%%%%%%%%%%%%%%%%%%%%%%%%%%%%%%%%%%%%%%%%%%%%%%%%%%%%%%%%%%%%%%%%%%%%%%%%%%%%%%%%%%%%%%

The modulus of the \(n\)-th ratio is comparable to
\(\rho_n/(|w-x_n|+b_n)\) on \(E_n\). Off \(\bigcup E_n\) the
derivative is bounded, since the possible remaining logarithm
in \(B_{x_n/8}(x_n)\setminus E_n\) is bounded. Hence the
transformed forcing \(|F'|^2 k^{-p}\) is bounded above by
\[
 Ck^{-p}+C\sum_n f_n,\qquad
 f_n(w)=\boldsymbol1_{E_n}(w)\,
 \rho_n^\beta(|w-x_n|+b_n)^{-\beta}v^{-p}.
\]
It is bounded below by \(c f_n\) on each \(E_n\).
The normalized potential of the first term is bounded on \(D\):
indeed \(k\ge c\,d_D\), and the smooth-boundary distance-potential
bound follows from Proposition~\ref{prop:negative-distance-potential-chart}
on a finite cover of \(D\).

%%%%%%%%%%%%%%%%%%%%%%%%%%%%%%%%%%%%%%%%%%%%%%%%%%%%%%%%%%%%%%%%%%%%%%%%%%%%%%%%%%%%%%%%%%%%%%%%%%%%%%%%%%%%%%%%%%%%%%%%%%%%%%%%%%%%%%%%%%%%%%%%%%%%%%%%%%%%%%%%%%%%%%%%%%%%%%%%%%%%
\vskip0.12in
%%%%%%%%%%%%%%%%%%%%%%%%%%%%%%%%%%%%%%%%%%%%%%%%%%%%%%%%%%%%%%%%%%%%%%%%%%%%%%%%%%%%%%%%%%%%%%%%%%%%%%%%%%%%%%%%%%%%%%%%%%%%%%%%%%%%%%%%%%%%%%%%%%%%%%%%%%%%%%%%%%%%%%%%%%%%%%%%%%%%

We estimate the other terms on \(B_{2/3}(0)\cap\mathbb H\)
using \(G_D\le G_{\mathbb H}\). The half-plane Green kernel obeys
\[
 \int_{\mathbb R} G_{\mathbb H}(x+iy,s+iv)\dd s=\min(y,v).
\]
This follows by integrating the reflected logarithmic kernel in \(s\).
For \(g\ge0\) supported in a half-disc of radius \(r\) centered at \(x_n\), the horizontal integral and the reflected formula give
\begin{equation}\label{eq:negative-linear-shell-bounds}
 \begin{aligned}
 \frac{G_{\mathbb H}[g](z)}{y}
 &\le C B &&\text{if }g\le B r^{-\beta}v^{-p},\\
 \frac{G_{\mathbb H}[g](z)}{y}
 &\le \frac{C}{|z-x_n|^2}\int v g(w)\dd A(w)
 &&\text{if }|z-x_n|\ge2r .
 \end{aligned}
\end{equation}
The first bound uses \(r^{-\beta}\int_0^r v^{-p}\dd v=1/\beta\). Split \(E_n\) into dyadic shells down to \(b_n\) and the innermost half-disc. On each piece \(f_n\) satisfies the first bound with \(B=C\rho_n^\beta\). Summing gives \(C\rho_n^\beta(1+\log(\rho_n/b_n))\le C(1+\sup_n\mu_n)\). As \(vf_n\le\rho_n^\beta\), the second bound on the whole \(E_n\) gives \(C\rho_n^{2+\beta}/|z-x_n|^2\) for \(|z-x_n|\ge2\rho_n\).
The disks of radius \(2\rho_n\) are disjoint, and outside
each such disk its contribution is at most \(C\rho_n^\beta\).
Summation gives a uniform bound in the half-disc. On the
rest of \(D\) compare with \(k\) by harmonicity, using the
bound on the cap. More explicitly, this comparison applies to
each finite partial sum of the nonnegative potentials with the
same constant, and monotone convergence gives the full bound.
Conformal invariance of the Green kernel now yields
\(\Lambda_p(K)<\infty\).

%%%%%%%%%%%%%%%%%%%%%%%%%%%%%%%%%%%%%%%%%%%%%%%%%%%%%%%%%%%%%%%%%%%%%%%%%%%%%%%%%%%%%%%%%%%%%%%%%%%%%%%%%%%%%%%%%%%%%%%%%%%%%%%%%%%%%%%%%%%%%%%%%%%%%%%%%%%%%%%%%%%%%%%%%%%%%%%%%%%%
\vskip0.12in
%%%%%%%%%%%%%%%%%%%%%%%%%%%%%%%%%%%%%%%%%%%%%%%%%%%%%%%%%%%%%%%%%%%%%%%%%%%%%%%%%%%%%%%%%%%%%%%%%%%%%%%%%%%%%%%%%%%%%%%%%%%%%%%%%%%%%%%%%%%%%%%%%%%%%%%%%%%%%%%%%%%%%%%%%%%%%%%%%%%%

Choose \(R_0\) common to the centers \(x_n\) in the flat chart and take \(\epsilon\) small enough that \(\rho_n/2<cR_0\). Use \eqref{eq:negative-half-disc-green-lower} at \(z_n=x_n+ib_n\) for \(2b_n<r<\rho_n/2\) and \(\arg(w-x_n)\in[\pi/4,3\pi/4]\).
Since \(k(z_n)\asymp b_n\), integration on the displayed
sector gives a lower bound for the normalized contribution of \(E_n\):
\[
 \frac{G_D[\boldsymbol1_{E_n}|F'|^2 k^{-p}](z_n)}{k(z_n)}
 \ge c\rho_n^\beta \int_{2b_n}^{\rho_n/2}\frac{\dd r}{r}
 \ge c'>0.
\]
The images \(F(E_n)\) lie uniformly close to \(Q\) as \(n\to\infty\).
Every sufficiently small physical collar therefore contains a
whole \(F(E_n)\), proving the lower bound for the localized
potential. Every positive solution satisfies \(u\le mK\)
by Theorem~\ref{thm:negative-full-representation}, so its
corresponding localized potential is bounded below by \(m^{-p}\)
times that of \(K^{-p}\). Existence follows from
Theorem~\ref{thm:negative-branch}.

%%%%%%%%%%%%%%%%%%%%%%%%%%%%%%%%%%%%%%%%%%%%%%%%%%%%%%%%%%%%%%%%%%%%%%%%%%%%%%%%%%%%%%%%%%%%%%%%%%%%%%%%%%%%%%%%%%%%%%%%%%%%%%%%%%%%%%%%%%%%%%%%%%%%%%%%%%%%%%%%%%%%%%%%%%%%%%%%%%%%
\vskip0.12in
%%%%%%%%%%%%%%%%%%%%%%%%%%%%%%%%%%%%%%%%%%%%%%%%%%%%%%%%%%%%%%%%%%%%%%%%%%%%%%%%%%%%%%%%%%%%%%%%%%%%%%%%%%%%%%%%%%%%%%%%%%%%%%%%%%%%%%%%%%%%%%%%%%%%%%%%%%%%%%%%%%%%%%%%%%%%%%%%%%%%

Finally \(d(F(w))\asymp v|F'(w)|\) on \(E_n\) by the Koebe
distance theorem on \(D\). On the same sector the pullback
density of the distance potential is thus bounded below by
\(c\rho_n^{\xi\gamma}r^{-\xi\gamma}v^{-\sigma}\).
The relation \(\xi\gamma+\sigma=1\) gives
\[
 \frac{P_\sigma(F(z_n))}{K(F(z_n))}
 \ge c\rho_n^{\xi\gamma}\log(\rho_n/(4b_n))
 \longrightarrow\infty,
\]
because \(\xi\gamma=\beta/(1+p)<\beta\).
\end{proof}
\begin{proposition}\label{prop:negative-micro-essential-zero}
On \(\Omega\), if \(0\le VK\le C K^{-p}\), then \(A=G_\Omega[V\,\cdot]\) is bounded on \(X\) with \(r_{\mathrm{ess},X}(A)=0\). If also \(VK\ge cK^{-p}\) on the cells' images, \(A\) is noncompact on \(X\).
\end{proposition}

\begin{proof}
With \(y=\operatorname{Im}w\) the pullback \(W=|F'|^2V\circ F\) satisfies \(d_D^2W\lesssim a_n(y/(|w-x_n|+b_n))^\beta\le a_n\) on \(E_n\), and \(d_D^2W\lesssim d_D^\beta\) off the cells. It is bounded and tends uniformly to zero at the boundary of \(D\), since a fixed cell has bounded derivative. The Hardy inequality and conformal invariance imply \(b(V)<\infty\) and \(b_\partial(V)=0\); the boundary homeomorphism takes shrinking collars to shrinking collars. The row bound follows from \(\Lambda_p<\infty\), so \cref{thm:negative-two-space-radius} gives \(r_{\mathrm{ess},X}(A)=0\).

For the lower assertion choose continuous \(0\le g_n\le1\) supported in \(E_n\), equal to one on \(4b_n\le |w-x_n|\le\rho_n/4\) and \(\pi/3\le\arg(w-x_n)\le2\pi/3\). The normalized image \(k^{-1}G_D[|F'|^2(VK)\circ F\,g_n]\) has a fixed positive lower bound at \(x_n+ib_n\) by \eqref{eq:negative-half-disc-green-lower} and the positive lower depth bound. It tends to zero at every fixed interior point by dominated convergence. No subsequence can converge in the normalized norm.

The shell estimates used below also follow from \eqref{eq:negative-linear-shell-bounds}. The forcing \(f=|F'|^2(VK)\circ F\) on \(E_n\) is bounded by \(Ca_n(|w-x_n|+b_n)^{-\beta}y^{-p}\). Each dyadic shell down to \(b_n\), and the innermost half-disc, contributes at most \(Ca_n\) to the normalized row. There are \(O(1+\log(\rho_n/b_n))\) such pieces. For a target on a shell of scale \(r\), the same and next shells each cost \(O(a_n)\). The remaining inner shells lie in a half-disc of radius at most half the target distance from \(x_n\). Since \(yf\le Ca_n\), their combined first moment is \(O(a_n r^2)\); the separated kernel bound gives \(O(a_n)\) for their combined row. This also covers targets in the innermost half-disc by using its single-shell bound.
\end{proof}

%%%%%%%%%%%%%%%%%%%%%%%%%%%%%%%%%%%%%%%%%%%%%%%%%%%%%%%%%%%%%%%%%%%%%%%%%%%%%%%%%%%%%%%%%%%%%%%%%%%%%%%%%%%%%%%%%%%%%%%%%%%%%%%%%%%%%%%%%%%%%%%%%%%%%%%%%%%%%%%%%%%%%%%%%%%%%%%%%%%%
\vskip0.32in
%%%%%%%%%%%%%%%%%%%%%%%%%%%%%%%%%%%%%%%%%%%%%%%%%%%%%%%%%%%%%%%%%%%%%%%%%%%%%%%%%%%%%%%%%%%%%%%%%%%%%%%%%%%%%%%%%%%%%%%%%%%%%%%%%%%%%%%%%%%%%%%%%%%%%%%%%%%%%%%%%%%%%%%%%%%%%%%%%%%%

\subsection{Actual inverse-power densities and large-mass uniqueness}

\begin{theorem}\label{thm:negative-micro-nonlinear}
On the full product and on all its finite prefixes, every sufficiently large mass has exactly one positive solution. The threshold is uniform for bounded depths and the fixed cap.
\end{theorem}

\begin{proof}
Write \(f_u=\lambda|F'|^2(u\circ F)^{-p}\) and note that \(\gamma(1-\xi)=1\). All images of the fixed cap, including the prefixes, have a common bound on their Lipschitz character and diameter. Near zero the argument and real-part estimates for \(F'\) give a uniform graph, and the lower Lipschitz bound separates this patch from the curved part of the cap. Away from a fixed neighborhood of zero, the locally convergent products extend across the real axis with uniformly bounded derivatives and \(|F'|\) bounded below. They give uniform charts there. Thus the physical Morrey bound \eqref{eq:negative-global-morrey} has uniform constants for the entire family. We use it to estimate the normalized Green potentials in \(D\).

Suppose a source set \(E\) is in \(B_R(x)\cap\mathbb H\),
where \(x\) is real, the rectangle
\((x-4R,x+4R)\times(0,4R)\) lies in \(D\), and the source
is in the upper half-disc of radius \(2/3\) used in the construction. If \(|F'|\le L\) on \(E\)
and \(\operatorname{diam} F(E)\le C L R\), then
\begin{equation}\label{eq:negative-micro-morrey-row}
 \sup_D \frac{G_D[\boldsymbol 1_E f_u]}{k}
 \lesssim L^\gamma R^{\gamma-1}.
\end{equation}
Indeed, the total source mass is at most \(C L^\gamma R^\gamma\)
by physical Morrey. For targets within the upper half-disc of
radius \(2/3\) and at distance at least \(2R\) from \(x\),
the reflected half-plane kernel gives the stronger bound
\(C L^\gamma R^{\gamma+1}|z-x|^{-2}\).
In the local rectangle \((x-4R,x+4R)\times(0,4R)\),
interior balls \(B_r\) with \(B_{4r}\) still interior that
meet \(E\) have \(\sup_{B_r}|F'|\lesssim L\) by distortion.
Thus their transformed source mass is at most \(C L^\gamma r^\gamma\).
Lemma~\ref{lem:morrey-green-domination}, applied on the scaled
rectangle, bounds its Green potential by \(C L^\gamma\) times
the potential of the distance density with exponent
\(\sigma=2-\gamma\). On \(B_{2R}(x)\cap\mathbb H\) this latter
potential is at most \(C(\operatorname{Im}z)R^{1-\sigma}\).
For the distance to the bottom side, this follows by
horizontal integration of the half-plane kernel. The other
sides give the same bound by splitting their weights into a
bounded part \(O(R^{-\sigma})\) and an integrable part
separated from the target by a distance comparable to \(R\).
The half-plane potential differs from the rectangle potential
by a harmonic lift. The far-target estimate on the outer
rectangle sides bounds this lift by
\(C(\operatorname{Im}z)L^\gamma R^{\gamma-1}\).
All the potentials have zero trace on the bottom side:
the density on \(E\) is in \(L^s\), \(s>1\), by the
change of variables, the derivative bound, and
Theorem~\ref{thm:negative-full-representation}.
Since \(k\gtrsim\operatorname{Im}z\) in the comparison
half-disc, these estimates prove the bound there.
Harmonic comparison with \(k\) extends it to \(D\).
This proves \eqref{eq:negative-micro-morrey-row}.

On \(E_n\), the derivative is bounded by
\(C a_n^{1/2}(|w-x_n|+b_n)^{-\xi}\).
Integration on radial segments from \(x_n\) gives
\[
 \operatorname{diam} F(B_r(x_n)\cap\mathbb H)
 \lesssim a_n^{1/2}r^{1-\xi}\qquad (b_n\le r\le\rho_n).
\]
A radial shell of radius \(r\), including the innermost disk
at scale \(b_n\), therefore contributes at most
\(C a_n^{\gamma/2}\) by \eqref{eq:negative-micro-morrey-row}.
For targets \(z\) on the \(i\)-th shell in \(E_n\), the
\emph{total} contribution from the source shells with index
at least \(i\) is also bounded by \(C a_n^{\gamma/2}\).
For the same and next shell use the estimate just proved.
All deeper shells are within half the target radius
\(r=|z-x_n|\). Their physical image has diameter
\(O(a_n^{1/2}r^{1-\xi})\), so their total source mass is
\(O(a_n^{\gamma/2}r)\), using \(\gamma(1-\xi)=1\).
The reflected kernel bound \(G_D(z,w)/k(z)\lesssim
(\operatorname{Im}w)/r^2\) then proves the assertion.
The same argument for \(z\notin E_n\), using the
outermost shells when \(|z-x_n|\) is comparable to \(\rho_n\)
and the far-target estimate otherwise, gives
\begin{equation}\label{eq:negative-micro-exterior-row}
 \sup_{D\setminus E_n}\frac{G_D[\boldsymbol 1_{E_n} f_u]}{k}
 \lesssim a_n^{\gamma/2}.
\end{equation}
The bounds outside the comparison half-disc follow by harmonic comparison. By lacunarity, \(\sum_n a_n^{1/a}<\infty\).

%%%%%%%%%%%%%%%%%%%%%%%%%%%%%%%%%%%%%%%%%%%%%%%%%%%%%%%%%%%%%%%%%%%%%%%%%%%%%%%%%%%%%%%%%%%%%%%%%%%%%%%%%%%%%%%%%%%%%%%%%%%%%%%%%%%%%%%%%%%%%%%%%%%%%%%%%%%%%%%%%%%%%%%%%%%%%%%%%%%%
\vskip0.12in
%%%%%%%%%%%%%%%%%%%%%%%%%%%%%%%%%%%%%%%%%%%%%%%%%%%%%%%%%%%%%%%%%%%%%%%%%%%%%%%%%%%%%%%%%%%%%%%%%%%%%%%%%%%%%%%%%%%%%%%%%%%%%%%%%%%%%%%%%%%%%%%%%%%%%%%%%%%%%%%%%%%%%%%%%%%%%%%%%%%%

For each fixed \(N\), the sources off \(\bigcup_{n\ge N}E_n\)
have uniformly bounded potentials in the normalized norm, and
those potentials converge in that norm under local convergence
of the solutions to a positive solution. To verify this, first
restrict to sets where \(|F'|\) is bounded. The same interior
ball transfer of Morrey and
Lemma~\ref{lem:morrey-green-domination} give domination by
a multiple of the distance potential in the smooth domain \(D\).
Use the localized version of that lemma in the \(D\) collar.
Choosing a larger distance exponent strictly below \(1\)
gives a uniform vanishing collar bound in the \(k\) norm.
Local convergence and interior truncation now give norm
convergence. This applies to the baseline and every finite set of cells.

%%%%%%%%%%%%%%%%%%%%%%%%%%%%%%%%%%%%%%%%%%%%%%%%%%%%%%%%%%%%%%%%%%%%%%%%%%%%%%%%%%%%%%%%%%%%%%%%%%%%%%%%%%%%%%%%%%%%%%%%%%%%%%%%%%%%%%%%%%%%%%%%%%%%%%%%%%%%%%%%%%%%%%%%%%%%%%%%%%%%
\vskip0.12in
%%%%%%%%%%%%%%%%%%%%%%%%%%%%%%%%%%%%%%%%%%%%%%%%%%%%%%%%%%%%%%%%%%%%%%%%%%%%%%%%%%%%%%%%%%%%%%%%%%%%%%%%%%%%%%%%%%%%%%%%%%%%%%%%%%%%%%%%%%%%%%%%%%%%%%%%%%%%%%%%%%%%%%%%%%%%%%%%%%%%

For large mass \(m\), the normalized deficit
\(m-(u\circ F)/k\) is the potential of \(f_u\) divided by \(k\).
Outside the microscopic disks it is bounded by a constant
independent of \(m,u\), by the baseline estimate and
\eqref{eq:negative-micro-exterior-row}.
For targets in \(E_n\), the same constant bounds the contributions
from the baseline, other disks, and all own shells with index
at least the target shell. If \(u\circ F\ge (m/2)k\) on
the preceding shells in \(E_n\), the remaining density is bounded
by \(C m^{-p}|F'|^2k^{-p}\) there. The linear shell estimate
in Proposition~\ref{prop:negative-micro-essential-zero} bounds
its normalized contribution by \(C m^{-p}a_n(1+\log(\rho_n/b_n))\le C' m^{-p}\).
Taking \(m\) large enough, induction from the outermost shell
gives \(u\ge (m/2)K\) throughout every disk as well.
On this interval the normalized Lipschitz constant of the Green iteration is at most \(\lambda p(m/2)^{-p-1}\Lambda_p\), which is strictly smaller than one for large \(m\). The maximal-branch theorem gives existence, and the contraction proves uniqueness among all positive solutions at those masses.

\end{proof}
\section{\bf Microscopic loss and endpoint classification}\label{sec:negative-microscopic}

\subsection{Notation for the microscopic cells}

Let \(\kappa=\partial_y k(0)>0\) and put
\[
 C_0=\prod_j\left(\frac{-x_j+i\rho_j}{-x_j+ib_j}\right)^\xi,\qquad
 c_p=\frac2\pi\int_0^\pi\sin^{1-p}\theta\,\dd\theta,\qquad
 D_p=\frac a2\,\lambda |C_0|^2c_p .
\]
For a prefix, \(C_0\) and \(D_p\) use its retained factors.

\subsection{One quantitative cell law and its boundary consequences}

Put \(a=1+p\) and \(\beta=1-p\). In a retained microscopic
cell use \(r=\rho_n e^{-t}\), \(z=x_n+re^{i\theta}\), and
\[
 W(t,\theta)=v(z)/r,\qquad
 Y(t)=\frac2\pi\int_0^\pi W(t,\theta)\sin\theta\,\dd\theta,
 \qquad q=a_n,\qquad L=\mu_n/q.
\]
Fix a mass cap and a compact range of positive depths \(\mu_n\). All constants below are uniform over the full products and their prefixes with the same cap \(D\) and kernel normalization. Estimates for sufficiently late cells mean \(x_n\le x_0\) for a fixed \(x_0>0\), which can be chosen after the comparison widths are fixed.

\begin{theorem}\label{thm:negative-cell-master}
For every positive solution on the full product or a prefix,
\(Y\) decreases to a limit \(s_n\ge0\). For all sufficiently
late retained cells,
\begin{equation}\label{eq:negative-cell-partial-loss}
 \sup_{t\ge0}\left|Y(t)^a-Y(0)^a+D_p a_n\min(t,L)\right|
 \le C\mathcal E_n,\qquad
 \mathcal E_n=a_n\log(e/a_n)+x_n\log(e/x_n).
\end{equation}
In particular,
\begin{equation}\label{eq:negative-quantitative-cell-loss}
 \left|Y(0)^a-s_n^a-D_p\mu_n\right|\le C\mathcal E_n.
\end{equation}
Here \(D_p\) is evaluated for the full product or the prefix in use. In particular
\begin{equation}\label{eq:negative-flat-retained-inlet}
 Y(0)^a\ge D_p\mu_n
       -C\left\{a_n\log(e/a_n)+x_n\log(e/x_n)\right\}.
\end{equation}
At each fixed smooth center, \(s_n\) is the inward derivative
of \(v\), by boundary \(C^1\) regularity. The estimate is
uniform over the stated bounded-mass class.
\end{theorem}

\begin{proof}
Put \(H=\lambda r^{1-p}|F'(z)|^2\). The equation of \(v=rW\) in
the polar strip gives \((\partial_t^2-2\partial_t+\partial_\theta^2+1)W=HW^{-p}\).
Write
\(f=(2/\pi)\int_0^\pi HW^{-p}\sin\theta\,\dd\theta\).
The mass bound gives \(0<W\le C\sin\theta\). The geometric
shell envelope and first-mode equation are
\begin{equation}\label{eq:negative-micro-plateau-envelope}
 H\asymp h(t):=q e^{-\beta(t-L)_+},\qquad
 Y''-2Y'=f,\qquad
 -Y'(t)=\int_0^\infty e^{-2s}f(t+s)\,\dd s.
\end{equation}
Boundedness excludes the increasing homogeneous mode.
Thus \(Y\) decreases. For \(0\le t\le L\), Jensen's
inequality and the first unit of the convolution give
\(-(Y^a)'(t)\gtrsim q\), since \(h(t+s)\gtrsim q\)
for \(0\le s\le1\). Integrating gives
\begin{equation}\label{eq:negative-micro-Y-lower}
 Y(t)^a\ge cq(L-t)_+.
\end{equation}
At any fixed smooth center, the local Morrey estimate
transfers to the preimage and to the odd reflection of
the forcing. The Newtonian estimate of
Lemma~\ref{lem:negative-interior-compact-index} and
harmonic reflection give \(C^1\) regularity up to the
diameter there, identifying \(s_n\) with the inward
derivative.

%%%%%%%%%%%%%%%%%%%%%%%%%%%%%%%%%%%%%%%%%%%%%%%%%%%%%%%%%%%%%%%%%%%%%%%%%%%%%%%%%%%%%%%%%%%%%%%%%%%%%%%%%%%%%%%%%%%%%%%%%%%%%%%%%%%%%%%%%%%%%%%%%%%%%%%%%%%%%%%%%%%%%%%%%%%%%%%%%%%%
\vskip0.12in
%%%%%%%%%%%%%%%%%%%%%%%%%%%%%%%%%%%%%%%%%%%%%%%%%%%%%%%%%%%%%%%%%%%%%%%%%%%%%%%%%%%%%%%%%%%%%%%%%%%%%%%%%%%%%%%%%%%%%%%%%%%%%%%%%%%%%%%%%%%%%%%%%%%%%%%%%%%%%%%%%%%%%%%%%%%%%%%%%%%%

Use a logarithmic rectangle of fixed half-width \(M\ge1\) about \(t\). The solution \(v\) has nonnegative cap values and zero data on the angular sides. Its source in the rectangle is \(\lambda r^2|F'|^2 v^{-p}\). The derivative scale \(r|F'|\) is comparable throughout the rectangle to its value at the center by the product envelope. The image of a ball of radius \(\ell\) therefore has diameter at most \(C_M r|F'|\ell\); the integral of the source equals the corresponding physical source integral. The physical Morrey estimate transfers with exponent \(\gamma=2/a\) and constant \(C_M(r|F'|)^\gamma=C_M\lambda^{-1/a}rH^{1/a}\), with \(H\) evaluated at the center.

%%%%%%%%%%%%%%%%%%%%%%%%%%%%%%%%%%%%%%%%%%%%%%%%%%%%%%%%%%%%%%%%%%%%%%%%%%%%%%%%%%%%%%%%%%%%%%%%%%%%%%%%%%%%%%%%%%%%%%%%%%%%%%%%%%%%%%%%%%%%%%%%%%%%%%%%%%%%%%%%%%%%%%%%%%%%%%%%%%%%
\vskip0.12in
%%%%%%%%%%%%%%%%%%%%%%%%%%%%%%%%%%%%%%%%%%%%%%%%%%%%%%%%%%%%%%%%%%%%%%%%%%%%%%%%%%%%%%%%%%%%%%%%%%%%%%%%%%%%%%%%%%%%%%%%%%%%%%%%%%%%%%%%%%%%%%%%%%%%%%%%%%%%%%%%%%%%%%%%%%%%%%%%%%%%

Apply \cref{lem:morrey-green-domination} on the rectangle. On its central slice the Green potential of the distance density of exponent \(-\sigma\) is bounded by \(C_M\sin\theta\). For the weight from the angular sides, integrate the infinite-strip Green kernel in the time variable to get the one-dimensional angular Green operator. The weights from the two caps satisfy the same bound by separation from the central slice near the caps and boundedness away from the caps. This gives a Green correction bounded by \(C_Mr h(t)^{1/a}\sin\theta\).

%%%%%%%%%%%%%%%%%%%%%%%%%%%%%%%%%%%%%%%%%%%%%%%%%%%%%%%%%%%%%%%%%%%%%%%%%%%%%%%%%%%%%%%%%%%%%%%%%%%%%%%%%%%%%%%%%%%%%%%%%%%%%%%%%%%%%%%%%%%%%%%%%%%%%%%%%%%%%%%%%%%%%%%%%%%%%%%%%%%%
\vskip0.12in
%%%%%%%%%%%%%%%%%%%%%%%%%%%%%%%%%%%%%%%%%%%%%%%%%%%%%%%%%%%%%%%%%%%%%%%%%%%%%%%%%%%%%%%%%%%%%%%%%%%%%%%%%%%%%%%%%%%%%%%%%%%%%%%%%%%%%%%%%%%%%%%%%%%%%%%%%%%%%%%%%%%%%%%%%%%%%%%%%%%%

In the harmonic lift of the cap values, the \(j\)-th sine coefficient at the central slice is the sum of its cap coefficients divided by \(2\cosh(jM)\). Nonnegativity of the cap data and \(|\sin(j\theta)|\le j\sin\theta\) bound the higher modes relative to the first by \(C\sum_{j\ge2}j^2\cosh(M)/\cosh(jM)\le Ce^{-M}\). Divide by the central radius and use the first coefficient of \(W\) to compare the lift to \(Y\). For the block estimate apply the source mass bound on a unit strip and divide by the radius, which is comparable to its central value there. We obtain

\begin{align}
 \int_t^{t+1}f(s)\,\dd s&\le Ch(t)^{1/a},
 \label{eq:negative-micro-plateau-block}\\
 \left|\frac{W(t,\theta)}{Y(t)\sin\theta}-1\right|
 &\le Ce^{-M}+C_M\frac{h(t)^{1/a}}{Y(t)}.
 \label{eq:negative-micro-plateau-profile}
\end{align}
For fixed \(M\), taking \(x_n\) small makes \(e^M\rho_n\) smaller than \(x_n/8\). The comparison rectangles then avoid the other cells and extend across \(t=0\). The estimates apply to every retained cell once \(x_1\) is this small.
The block bound and the convolution imply
\(-Y'\le C_1q^{1/a}\). Fix a large rectangle width to make its
higher-mode error small and choose \(D\ge2C_1\) sufficiently large.
Then
\begin{equation}\label{eq:negative-cell-threshold-profile}
 Y(t)\ge Dq^{1/a}\quad\Longrightarrow\quad
                       W(t,\theta)\ge cY(t)\sin\theta.
\end{equation}

The following sharper estimates are the only additional
strip argument:
\begin{equation}\label{eq:negative-cell-sharp-profile}
 \begin{split}
 -Y'(t)&\le CqY(t)^{-p},\\
 |W(t,\theta)-Y(t)\sin\theta|
 &\le C\bigl(qY(t)^{-p}+e^{-t}\bigr)\sin\theta.
 \end{split}
\end{equation}
They hold when \(Y(t)\ge2Dq^{1/a}\), with \(t\) above a
fixed constant for the second estimate. Here are the kernel
details. Subtract the first target sine mode from the Green
representation of \(e^{-t}W\) on the half-strip. The left
harmonic lift contributes \(O(e^{-t}\sin\theta)\) after
conjugation, and the right lift vanishes at infinity. The
higher-mode kernel before conjugation is bounded, for
\(|s-t|\ge1\), by
\(Ce^{-2|s-t|}\sin\theta\sin\varphi\). For
\(|s-t|<1\), its absolute value is bounded by the full
strip Green kernel plus \(C\sin\theta\sin\varphi\);
its integral against \(\sin^{-p}\varphi\) over this
unit band is \(O(\sin\theta)\).
The speed bound makes the first future halving of \(Y(t)\)
at least \(Y(t)/(2C_1 q^{1/a})\ge2\) strip units away, so
the separated kernel estimate applies beyond it.
Up to that time, and at all earlier times beyond the fixed
left strip, \eqref{eq:negative-cell-threshold-profile}
bounds the forcing by
\(CqY(t)^{-p}\sin^{-p}\varphi\).
Separated tails need only the angular first-moment block
bound in \eqref{eq:negative-micro-plateau-block}; their
contribution is at most
\(Cq^{1/a}e^{-cY(t)/q^{1/a}}\sin\theta\), bounded by
\(CqY(t)^{-p}\sin\theta\). The fixed left strip costs
\(Ce^{-t}\sin\theta\). This proves the second estimate.
The same split in the first-mode convolution proves the
first estimate, without a higher-mode kernel.

%%%%%%%%%%%%%%%%%%%%%%%%%%%%%%%%%%%%%%%%%%%%%%%%%%%%%%%%%%%%%%%%%%%%%%%%%%%%%%%%%%%%%%%%%%%%%%%%%%%%%%%%%%%%%%%%%%%%%%%%%%%%%%%%%%%%%%%%%%%%%%%%%%%%%%%%%%%%%%%%%%%%%%%%%%%%%%%%%%%%
\vskip0.12in
%%%%%%%%%%%%%%%%%%%%%%%%%%%%%%%%%%%%%%%%%%%%%%%%%%%%%%%%%%%%%%%%%%%%%%%%%%%%%%%%%%%%%%%%%%%%%%%%%%%%%%%%%%%%%%%%%%%%%%%%%%%%%%%%%%%%%%%%%%%%%%%%%%%%%%%%%%%%%%%%%%%%%%%%%%%%%%%%%%%%

Let \(\mathcal P_{\ne n}\) be the product with the \(n\)-th factor omitted and set
\begin{equation}\label{eq:negative-tuned-coefficient}
 C_n=\lambda|\mathcal P_{\ne n}(x_n)|^2,\qquad D_n=ac_p C_n/2.
\end{equation}
Factoring the central ratio at \(z=x_n+\rho_ne^{-t+i\theta}\) gives
\begin{equation}\label{eq:negative-local-cell-coefficient}
 H(t,\theta)=C_nq(1+e(t,\theta)),\qquad
 |e(t,\theta)|\le C(e^{-t}+e^{t-L}),\quad 0\le t\le L .
\end{equation}
The earlier factors have logarithmic derivatives \(O(\rho_j/x_j^2)\) each, and the later derivatives sum to \(O(\sum_{j>n}\rho_j/x_n^2)=O(1)\). Their variation from \(x_n\) is \(O(\rho_n\log(e/x_n)e^{-t})\). Finally the earlier factors vary from zero to \(x_n\) by \(O(nx_n)\) in the logarithm, and the \(n\)-th and later values at zero contribute \(O(x_n)\). This holds also for omitted tail factors. Lacunarity gives the comparison
\begin{equation}\label{eq:negative-micro-plateau}
 H=\lambda|C_0|^2q\left(1+O(e^{-t}+e^{t-L}+\delta_n)\right),
 \qquad \delta_n=x_n\log(e/x_n),\qquad 0\le t\le L .
\end{equation}

%%%%%%%%%%%%%%%%%%%%%%%%%%%%%%%%%%%%%%%%%%%%%%%%%%%%%%%%%%%%%%%%%%%%%%%%%%%%%%%%%%%%%%%%%%%%%%%%%%%%%%%%%%%%%%%%%%%%%%%%%%%%%%%%%%%%%%%%%%%%%%%%%%%%%%%%%%%%%%%%%%%%%%%%%%%%%%%%%%%%
\vskip0.12in
%%%%%%%%%%%%%%%%%%%%%%%%%%%%%%%%%%%%%%%%%%%%%%%%%%%%%%%%%%%%%%%%%%%%%%%%%%%%%%%%%%%%%%%%%%%%%%%%%%%%%%%%%%%%%%%%%%%%%%%%%%%%%%%%%%%%%%%%%%%%%%%%%%%%%%%%%%%%%%%%%%%%%%%%%%%%%%%%%%%%

Choose fixed large \(T,M_1\), and put \(U=L-M_1\).
On \([T,U]\), the moment lower bound makes
\eqref{eq:negative-cell-sharp-profile} applicable and
its relative error smaller than \(1/2\). Hence
\[
 Y^p f=\lambda|C_0|^2c_pq(1+O(E(t))),\qquad
 E(t)=qY^{-a}+e^{-t}/Y+e^{-t}+e^{t-L}+\delta_n,
\]
and
\[
 \int_T^U E(t)\,\dd t\le C(1+\log L+L\delta_n).
\]
For the first term use \(Y^a\gtrsim q(L-t)\); for the
second split at \(L/2\), using \(Y\gtrsim1\) first and
\(Y\gtrsim q^{1/a}\) afterwards.
Now set \(X=Y^a\) and use the exact identity
\begin{equation}\label{eq:negative-cell-power-identity}
 X''-2X'=aY^pf+apY^{p-1}(Y')^2.
\end{equation}
The speed estimate gives \(|X'|\le Cq\) throughout \([T,U]\), and
\[
 \int_T^U Y^{p-1}(Y')^2\,\dd t
 \le Cq\int_T^U\frac{-Y'}Y\,\dd t
 \le Cq\log(e/q).
\]
Integrating \eqref{eq:negative-cell-power-identity}
from \(T\) to any \(t\in[T,U]\) gives
\(X(T)-X(t)=D_pq(t-T)+O(q\log(e/q)+\delta_n)\),
uniformly in \(t\).

%%%%%%%%%%%%%%%%%%%%%%%%%%%%%%%%%%%%%%%%%%%%%%%%%%%%%%%%%%%%%%%%%%%%%%%%%%%%%%%%%%%%%%%%%%%%%%%%%%%%%%%%%%%%%%%%%%%%%%%%%%%%%%%%%%%%%%%%%%%%%%%%%%%%%%%%%%%%%%%%%%%%%%%%%%%%%%%%%%%%
\vskip0.12in
%%%%%%%%%%%%%%%%%%%%%%%%%%%%%%%%%%%%%%%%%%%%%%%%%%%%%%%%%%%%%%%%%%%%%%%%%%%%%%%%%%%%%%%%%%%%%%%%%%%%%%%%%%%%%%%%%%%%%%%%%%%%%%%%%%%%%%%%%%%%%%%%%%%%%%%%%%%%%%%%%%%%%%%%%%%%%%%%%%%%

Both end pieces cost only \(O(q)\) in the power variable.
The initial piece has fixed length, \(Y\gtrsim1\), and
the first estimate in \eqref{eq:negative-cell-sharp-profile}
applies using the negative strip extension. At the terminal
end the block bound gives \(Y(U)-s_n\le Cq^{1/a}\).
If \(s_n\le Kq^{1/a}\), then \(Y(U)^a-s_n^a\le Cq\).
If \(s_n>Kq^{1/a}\), choose \(K\) so large that the fixed
profile gives \(W\ge cY\sin\theta\) for all \(t\ge U\).
Then \(f\le Ch(t)s_n^{-p}\), and \(\int_U^\infty h\le Cq\)
gives \(Y(U)-s_n\le Cq s_n^{-p}\) and again \(Y(U)^a-s_n^a\le Cq\).
Monotonicity bounds every partial terminal loss by the
same constant. Combining the three intervals proves
\eqref{eq:negative-cell-partial-loss}.
\end{proof}

%%%%%%%%%%%%%%%%%%%%%%%%%%%%%%%%%%%%%%%%%%%%%%%%%%%%%%%%%%%%%%%%%%%%%%%%%%%%%%%%%%%%%%%%%%%%%%%%%%%%%%%%%%%%%%%%%%%%%%%%%%%%%%%%%%%%%%%%%%%%%%%%%%%%%%%%%%%%%%%%%%%%%%%%%%%%%%%%%%%%
\vskip0.32in
%%%%%%%%%%%%%%%%%%%%%%%%%%%%%%%%%%%%%%%%%%%%%%%%%%%%%%%%%%%%%%%%%%%%%%%%%%%%%%%%%%%%%%%%%%%%%%%%%%%%%%%%%%%%%%%%%%%%%%%%%%%%%%%%%%%%%%%%%%%%%%%%%%%%%%%%%%%%%%%%%%%%%%%%%%%%%%%%%%%%

\subsection{A matched logarithmic microscopic loss law}

For the logarithmic correction use \(q=a_n\), \(L=\mu_n/q\), and \eqref{eq:negative-local-cell-coefficient} with \(C_n\) and \(D_n\) as in the cell proof.

Put \(\phi=\sin\theta\),
\(\langle f,g\rangle=(2/\pi)\int_0^\pi fg\).
 \(\Pi_\perp\) denotes projection off \(\phi\), given by
\(\Pi_\perp g=g-\langle g,\phi\rangle\phi\) for \(g\in L^1(0,\pi)\).
On this complement set
\[
 \mathscr L_\perp=\partial_\theta^2+1,\qquad
 g=\Pi_\perp\phi^{-p},\qquad \chi=\mathscr L_\perp^{-1}g,\qquad
 J_p=-\langle g,\chi\rangle
       =\langle\chi',\chi'\rangle-\langle\chi,\chi\rangle>0.
\]
The weak Dirichlet inverse is well defined; the interval Green
estimate gives \(\chi\in H^1_0\cap W^{2,1}\),
\(\chi\perp\phi\), and \(|\chi|\lesssim\phi\).
Positivity of \(J_p\) follows from the spectral gap on the
complement and \(g\ne0\). Define
\[
 \Theta_p=\frac{2p}{a}\left(\frac14+\frac{J_p}{c_p^2}\right),
 \qquad \kappa_n=D_n\Theta_p
       =pC_n\left(\frac{c_p}{4}+\frac{J_p}{c_p}\right)>0.
\]

\begin{theorem}\label{thm:negative-cell-signed-log}
For every bounded-mass positive solution, uniformly over late
retained cells in full products and finite prefixes, with \(\mu_n\) in a fixed positive compact interval,
\begin{equation}\label{eq:negative-cell-matched-log}
 Y_n(0)^a-s_n^a-D_n\mu_n
 =\kappa_nq\log\left(1+\frac1{s_n^a+q}\right)+O(q).
\end{equation}
The error is uniform over all terminal slopes \(s_n\ge0\).
In particular, for each \(0<\vartheta<1\),
\begin{equation}\label{eq:negative-cell-signed-log}
 Y_n(0)^a-s_n^a-D_n\mu_n
 \ge\vartheta\frac{pC_nc_p}{4}\,q
       \log\left(1+\frac1{s_n^a+q}\right)-C_\vartheta q,
\end{equation}
and, for some uniform \(c>0\),
\begin{equation}\label{eq:negative-cell-critical-inlet-log}
                 Y_n(0)^a\ge D_n\mu_n+cq\log(1/q)
\end{equation}
for all sufficiently late cells. Stability is not required.
\end{theorem}

\begin{proof}
Omit the cell index, and write
\(P=-Y'\), \(X=Y^a\), \(Z=W-Y\phi\), and
\(\mathscr D=\partial_t^2-2\partial_t+\partial_\theta^2+1\).
Choose fixed large \(T,K\), and put \(I=[T,U]\), \(U=L-K\),
with unit-width buffers. The master estimates give, on these intervals,
\[
 X\ge cq(L-t),\quad P\le CqY^{-p},\quad
 W\asymp Y\phi,\quad |Y''|\le CqY^{-p},
\]
and
\begin{equation}\label{eq:negative-cell-leading-forcing}
 Y^pf=C_nc_pq(1+O(E)),\qquad
 E=q/X+e^{-t}/Y+e^{-t}+e^{t-L}.
\end{equation}
Increasing \(T,K\) and then taking the cell late makes the
relative angular errors as small as needed. The moment lower
bound implies
\[
 q^3\int_I X^{-2}\,\dd t+q^4\int_I X^{-3}\,\dd t=O(q).
\]
All terms containing the two exponential errors have the same
bound when integrated with the weights below: split at \(L/2\),
using \(Y\gtrsim1\) first and \(Y\gtrsim q^{1/a}\) afterwards.

Separate the higher-mode homogeneous lift from the inlet:
\(Z=Z_h+u\), with \(\mathscr D Z_h=0\), \(Z_h(0)=Z(0)\),
and increasing modes excluded. Its sine series gives
\(|Z_h|\le Ce^{-t}\phi\). The forced component satisfies
\[
 \mathscr D u=\Pi_\perp(HW^{-p}),\qquad \Pi_\perp u=u,\qquad
                         |u|\le CqY^{-p}\phi.
\]
For the last estimate use the higher-mode kernel from the master
proof, after removing its homogeneous lift. Up to the first future
halving of \(Y(t)\), the source is bounded by
\(CqY(t)^{-p}\phi^{-p}\). The separated later blocks cost at most
\(Cq^{1/a}e^{-cY(t)/q^{1/a}}\phi\), which obeys the same bound.
The fixed initial strip has source \(O(q)\phi^{-p}\).
Thus no independent \(e^{-t}\) error remains in \(u\).

%%%%%%%%%%%%%%%%%%%%%%%%%%%%%%%%%%%%%%%%%%%%%%%%%%%%%%%%%%%%%%%%%%%%%%%%%%%%%%%%%%%%%%%%%%%%%%%%%%%%%%%%%%%%%%%%%%%%%%%%%%%%%%%%%%%%%%%%%%%%%%%%%%%%%%%%%%%%%%%%%%%%%%%%%%%%%%%%%%%%
\vskip0.12in
%%%%%%%%%%%%%%%%%%%%%%%%%%%%%%%%%%%%%%%%%%%%%%%%%%%%%%%%%%%%%%%%%%%%%%%%%%%%%%%%%%%%%%%%%%%%%%%%%%%%%%%%%%%%%%%%%%%%%%%%%%%%%%%%%%%%%%%%%%%%%%%%%%%%%%%%%%%%%%%%%%%%%%%%%%%%%%%%%%%%

Test the angular equation against \(\chi\), with time weight
\(\omega=q/Y\) times the square of a smooth cutoff equal to one
on \(I\) and supported in its buffers. Since
\(\mathscr L_\perp u=\Pi_\perp(HW^{-p})-u_{tt}+2u_t\), two time
integrations by parts transfer the last terms to
\(-[(\zeta\omega)''+2(\zeta\omega)']\langle\chi,u\rangle\).
On the bulk their absolute value is at most
\[
 C\left(q^3/X^2+q^4/X^3\right).
\]
Indeed \(|\omega'|+|\omega''|\lesssim
q^2Y^{-p-2}+q^3Y^{-2p-3}\). Buffer terms are
\(O(q^2/X)=O(q)\). The source in this pairing may be replaced
by \(C_nqY^{-p}g\); its integrated error is bounded by
\(Cq^2\int_I E/X=O(q)\). Removing the cutoff gives
\begin{equation}\label{eq:negative-cell-angular-energy}
 \int_I\frac qY\langle g,u\rangle\,\dd t
       =-C_nJ_pq^2\int_I X^{-1}\,\dd t+O(q).
\end{equation}
The same weighted pairing with \(Z_h\) is \(O(q)\).
All integrations are legitimate by weak testing: the local reflected
Morrey estimate gives \(C^1\) regularity up to the angular boundary,
and \(u,\chi=O(\phi)\) make the angular endpoint terms vanish.

%%%%%%%%%%%%%%%%%%%%%%%%%%%%%%%%%%%%%%%%%%%%%%%%%%%%%%%%%%%%%%%%%%%%%%%%%%%%%%%%%%%%%%%%%%%%%%%%%%%%%%%%%%%%%%%%%%%%%%%%%%%%%%%%%%%%%%%%%%%%%%%%%%%%%%%%%%%%%%%%%%%%%%%%%%%%%%%%%%%%
\vskip0.12in
%%%%%%%%%%%%%%%%%%%%%%%%%%%%%%%%%%%%%%%%%%%%%%%%%%%%%%%%%%%%%%%%%%%%%%%%%%%%%%%%%%%%%%%%%%%%%%%%%%%%%%%%%%%%%%%%%%%%%%%%%%%%%%%%%%%%%%%%%%%%%%%%%%%%%%%%%%%%%%%%%%%%%%%%%%%%%%%%%%%%

Taylor expansion of \(W^{-p}\) around \(Y\phi\) now yields
\begin{equation}\label{eq:negative-cell-forcing-signed}
 \int_IY^pf\,\dd t
 =C_nc_pq|I|+pC_n^2J_pq^2\int_I X^{-1}\,\dd t+O(q).
\end{equation}
Its quadratic error is bounded by
\(Cq\int_I(q/X+e^{-t}/Y)^2=O(q)\), and the coefficient
rounding error integrates to \(O(q)\). The positive angular
correction is therefore retained with its exact coefficient.

%%%%%%%%%%%%%%%%%%%%%%%%%%%%%%%%%%%%%%%%%%%%%%%%%%%%%%%%%%%%%%%%%%%%%%%%%%%%%%%%%%%%%%%%%%%%%%%%%%%%%%%%%%%%%%%%%%%%%%%%%%%%%%%%%%%%%%%%%%%%%%%%%%%%%%%%%%%%%%%%%%%%%%%%%%%%%%%%%%%%
\vskip0.12in
%%%%%%%%%%%%%%%%%%%%%%%%%%%%%%%%%%%%%%%%%%%%%%%%%%%%%%%%%%%%%%%%%%%%%%%%%%%%%%%%%%%%%%%%%%%%%%%%%%%%%%%%%%%%%%%%%%%%%%%%%%%%%%%%%%%%%%%%%%%%%%%%%%%%%%%%%%%%%%%%%%%%%%%%%%%%%%%%%%%%

The radial contribution needs no further profile expansion.
The first-mode equation says \(2P=f+P'\). Hence, with
\(R_I=\int_IY^{p-1}P^2\), integration by parts gives
\[
 R_I=\frac12\int_IY^{p-1}Pf+O(q)
     =\frac{C_nc_p}{2a}\,q\log\frac{X(T)}{X(U)}+O(q).
\]
The first error has endpoint terms \(O(q^2/X)=O(q)\)
and interior terms \(O(q^3\int_I X^{-2})=O(q)\).
The error from \eqref{eq:negative-cell-leading-forcing}
is \(O(q\int_I(P/Y)E)=O(q)\).
Similarly,
\[
 \log\frac{X(T)}{X(U)}
 =\frac a2\int_I\frac fY+\frac a2\int_I\frac{P'}Y
 =D_nq\int_I X^{-1}+O(1).
\]
Here \(\int_I P'/Y=[P/Y]_T^U-\int_I P^2/Y^2=O(1)\),
and \(q\int_I E/X=O(1)\). Consequently
\begin{equation}\label{eq:negative-cell-log-time}
 q^2\int_I X^{-1}\,\dd t
         =\frac q{D_n}\log\frac{X(T)}{X(U)}+O(q).
\end{equation}

Integrate the exact power identity to obtain
\[
 X(T)-X(U)=\frac a2\int_IY^pf
                    +\frac{ap}{2}R_I-\frac{X'(U)-X'(T)}2.
\]
The endpoint derivatives and the fixed entrance and terminal
power losses are \(O(q)\), by the master proof. Combining
\eqref{eq:negative-cell-forcing-signed} and
\eqref{eq:negative-cell-log-time} therefore gives the logarithmic
coefficient
\[
 \frac{apC_n^2J_p}{2D_n}+\frac{pC_nc_p}{4}
                  =D_n\Theta_p.
\]
Finally \(X(T)\) is bounded above and below uniformly, while
\(X(U)\asymp s_n^a+q\): use monotonicity, the moment lower
bound, and the terminal \(O(q)\) power loss.
Thus \(\log(X(T)/X(U))\) differs from
\(\log(1+1/(s_n^a+q))\) by \(O(1)\).
This proves \eqref{eq:negative-cell-matched-log} uniformly
over all terminal slopes, and implies
\eqref{eq:negative-cell-signed-log}.

If \(s_n^a\ge\sqrt q\), that term dominates \(q\log(1/q)\).
Otherwise the logarithm in the matched formula is at least
\(\tfrac12\log(1/q)-C\). The uniformly positive coefficient
then proves \eqref{eq:negative-cell-critical-inlet-log}.
\end{proof}

%%%%%%%%%%%%%%%%%%%%%%%%%%%%%%%%%%%%%%%%%%%%%%%%%%%%%%%%%%%%%%%%%%%%%%%%%%%%%%%%%%%%%%%%%%%%%%%%%%%%%%%%%%%%%%%%%%%%%%%%%%%%%%%%%%%%%%%%%%%%%%%%%%%%%%%%%%%%%%%%%%%%%%%%%%%%%%%%%%%%
\vskip0.32in
%%%%%%%%%%%%%%%%%%%%%%%%%%%%%%%%%%%%%%%%%%%%%%%%%%%%%%%%%%%%%%%%%%%%%%%%%%%%%%%%%%%%%%%%%%%%%%%%%%%%%%%%%%%%%%%%%%%%%%%%%%%%%%%%%%%%%%%%%%%%%%%%%%%%%%%%%%%%%%%%%%%%%%%%%%%%%%%%%%%%

\subsection{The exterior trace and exact limiting loss}

\begin{theorem}\label{thm:negative-micro-exact-drop}
On every flat product with \(\mu_n\to\mu>0\),
every positive solution \(u\) has a finite exterior trace
\[
 A(u)=\lim_{\substack{z\to0\\z\notin\bigcup_n E_n}}
                      \frac{u(F(z))}{\operatorname{Im}z}.
\]
If \(s_n(u)=\partial_y(u\circ F)(x_n)\) is the slope at the smooth
center of the \(n\)-th microscopic cap, then
\begin{equation}\label{eq:negative-micro-exact-drop}
 \lim_{n\to\infty}s_n(u)^a=A(u)^a-D_p\mu .
\end{equation}
In particular \(A(u)\ge(D_p\mu)^{1/a}\).
A positive lower bound for \(u/K\) in a neighborhood of \(Q\)
holds exactly when \(A(u)>(D_p\mu)^{1/a}\).
The global lower bound \(u\ge cK\) holds exactly when this strict
inequality and positive inward slopes at every smooth off-pole
boundary point both hold.
\end{theorem}

\begin{proof}
Use the source decomposition in
Theorem~\ref{thm:negative-micro-nonlinear}.
For fixed \(N\), the sources outside \(\bigcup_{n\ge N}E_n\) have normalized potentials uniformly approximable by compact interior sources, using collar truncation for the bounded derivative part. The Morrey bound controls their integrals uniformly. For fixed interior support, reflection and the harmonic derivative estimate give a common modulus at zero for the quotient by \(k\). At targets outside the cells the remaining normalized potential is bounded by \(C\sum_{n\ge N}a_n^{1/a}\) by \eqref{eq:negative-micro-exterior-row}. The Green representation, with \(N\to\infty\), gives the exterior limit of \(u\circ F/k\), uniformly over bounded masses on the fixed product. Multiplication by \(k/y\to\kappa\) gives \(A(u)\). On the inlet semicircle \(W_n(0,\theta)=(u\circ F/y)\sin\theta\), so \(Y_n(0)\to A(u)\) uniformly. The cell law of \cref{thm:negative-cell-master} and \(\mu_n\to\mu\) give \eqref{eq:negative-micro-exact-drop}.

%%%%%%%%%%%%%%%%%%%%%%%%%%%%%%%%%%%%%%%%%%%%%%%%%%%%%%%%%%%%%%%%%%%%%%%%%%%%%%%%%%%%%%%%%%%%%%%%%%%%%%%%%%%%%%%%%%%%%%%%%%%%%%%%%%%%%%%%%%%%%%%%%%%%%%%%%%%%%%%%%%%%%%%%%%%%%%%%%%%%
\vskip0.12in
%%%%%%%%%%%%%%%%%%%%%%%%%%%%%%%%%%%%%%%%%%%%%%%%%%%%%%%%%%%%%%%%%%%%%%%%%%%%%%%%%%%%%%%%%%%%%%%%%%%%%%%%%%%%%%%%%%%%%%%%%%%%%%%%%%%%%%%%%%%%%%%%%%%%%%%%%%%%%%%%%%%%%%%%%%%%%%%%%%%%

If \(A(u)>(D_p\mu)^{1/a}\), the slopes \(s_n(u)\) have a
positive lower bound for all sufficiently large \(n\).
Monotonicity gives the same lower bound for \(Y_n\) on each
whole microscopic strip. Use
\eqref{eq:negative-micro-plateau-profile} with fixed large \(M\)
and \(h_n\le a_n\to0\) to get \(W_n\gtrsim\sin\theta\).
Since \(k\asymp r\sin\theta\) there, this gives the desired
uniform lower bound within the late \(E_n\). Outside the
\(E_n\) the exterior trace is positive. Conversely, equality
\(A(u)=(D_p\mu)^{1/a}\) forces \(s_n(u)\to0\), and the smooth
boundary limits of \(u\circ F/k\) at \(x_n\) tend to zero
because \(\partial_y k(x_n)\to\kappa>0\).
This rules out a positive lower bound in every neighborhood of
\(Q\). Finally the pole has a positive relative trace by its
smoothness, by Proposition~\ref{prop:negative-distance-potential-chart}
and \eqref{eq:negative-large-potential}.
On compact smooth off-pole arcs the ratio extends continuously
by the inward derivatives. This and interior positivity give
the global equivalence.

Write \(v=u\circ F\), \(J=|F'|^2\), \(k=K\circ F\),
and \(P_Q(w)=\partial_{y_z}G_D(0,w)\).
The exterior slope is the Green flux
\begin{equation}\label{eq:negative-kernel-exterior-identity}
 A(u)=\kappa m-\lambda\int_D P_Q Jv^{-p}.
\end{equation}
Indeed, for a compact interior source the boundary derivative
of its potential is the integral against \(P_Q\).
The source approximation used for the exterior trace in
Theorem~\ref{thm:negative-micro-exact-drop} converges uniformly
in the quotient by \(y\) along the imaginary axis.
In that approximation the normalized exterior contribution
of the microscopic tail is bounded by
\(C\sum_{n\ge N}a_n^{1/(1+p)}\). Monotone exhaustion proves
\eqref{eq:negative-kernel-exterior-identity}.
For a finite prefix the same trace and flux formula hold with \(J_N\), by the same compact-source approximation.
\end{proof}

%%%%%%%%%%%%%%%%%%%%%%%%%%%%%%%%%%%%%%%%%%%%%%%%%%%%%%%%%%%%%%%%%%%%%%%%%%%%%%%%%%%%%%%%%%%%%%%%%%%%%%%%%%%%%%%%%%%%%%%%%%%%%%%%%%%%%%%%%%%%%%%%%%%%%%%%%%%%%%%%%%%%%%%%%%%%%%%%%%%%
\vskip0.32in
%%%%%%%%%%%%%%%%%%%%%%%%%%%%%%%%%%%%%%%%%%%%%%%%%%%%%%%%%%%%%%%%%%%%%%%%%%%%%%%%%%%%%%%%%%%%%%%%%%%%%%%%%%%%%%%%%%%%%%%%%%%%%%%%%%%%%%%%%%%%%%%%%%%%%%%%%%%%%%%%%%%%%%%%%%%%%%%%%%%%

\subsection{Cancellation and the flat inlet error}

Put \(\beta=1-p\), and use the flat product with lacunary centers
\(x_{j+1}\le x_j/2\), \(\rho_j=\epsilon x_j^2\), and depths
\(\mu_j=a_j\log(\rho_j/b_j)\) in a fixed compact
subinterval of \((0,\infty)\). Here \(a_j=\rho_j^\beta\).
The estimates concern all positive solutions of bounded mass.

\begin{lemma}\label{lem:negative-flat-inlet-refined}
Let \(v=u\circ F\) have exterior slope \(A>0\), and let
\(Y_n(0)\) be its inlet moment. Then
\begin{equation}\label{eq:negative-flat-sharp-inlet}
             |Y_n(0)-A|\le C(x_n+a_n).
\end{equation}
The constant is uniform for bounded masses and exterior
slopes bounded below by a positive constant. In particular,
for \(p\ge1/2\),
\[
       Y_n(0)^a-A^a=O(a_n)
                      =o\bigl(a_n\log(1/a_n)\bigr).
\]
\end{lemma}

\begin{proof}
Work in a fixed pole-free half-disc \(D_R\) centered at zero within the flat chart. We choose \(R\) below using the uniform exterior modulus. Let \(G\) be its Dirichlet kernel, \(P_0(w)=\partial_{y_z}G(0,w)\), \(y_z=\operatorname{Im}z\), and \(y_w=\operatorname{Im}w\). Sources in the Green integrals are restricted to \(D_R\). Write
\[
 K_z(w)=G(z,w)/y_z-P_0(w),
 \qquad r=|z|.
\]
For \(|w|\ge4r\), harmonic reflection and the separated
half-plane kernel bound give
\begin{equation}\label{eq:negative-flat-kernel-increment}
             |K_z(w)|\le Cr\,y_w/|w|^3.
\end{equation}
For sources with \(|w|\le4r\), horizontal integration
of the normalized half-plane Green kernel is at most
one. Together with \(P_0(w)\lesssim y_w/|w|^2\), this
proves, for every \(\alpha\ge0\),
\begin{equation}\label{eq:negative-flat-weighted-increment}
 \int |K_z(w)|\,|w|^\alpha y_w^{-p}\,\dd w
 \le C_\alpha
 \begin{cases}
 r^{\alpha+\beta},&\alpha<p,\\
 r\log(e/r),&\alpha=p,\\
 r,&\alpha>p.
 \end{cases}
\end{equation}
Indeed the near contribution is \(O(r^{\alpha+\beta})\),
whereas the far contribution is bounded by
\(Cr\int_r^R s^{\alpha-p-1}\,\dd s\).
The same argument with the additional weight
\(|w|\log(e/|w|)\) gives an \(O(r)\) bound.

We first improve the exterior estimate for an actual
microscopic source. Put \(f=\lambda Jv^{-p}\).
For each sufficiently late cell,
\begin{equation}\label{eq:negative-flat-cell-force-moment}
 \int_{E_j}y_w f\le Ca_j\rho_j^2,
 \qquad
 \sup_{z\notin E_j}\frac{G[\mathbf1_{E_j}f](z)}{
                                    y_z}
       \le Ca_j .
\end{equation}
To prove the first estimate, use the strip variables
and let \(f_j(t)\) be the first sine moment of
\(H_jW_j^{-p}\). The moment lower bound gives
\(Y_j\ge c>0\) on \(0\le t\le L_j/2\).
The fixed-width relative profile estimate therefore
gives \(W_j\ge c\sin\theta\) there, and hence
\(f_j(t)\le Ca_j\). On the remaining strip the
unit-block source estimate gives
\(\int_t^{t+1}f_j\le Ca_j^{1/a}\), with the terminal
exponential cutoff. Consequently
\[
 \int_{E_j}y_w f
 =\frac\pi2\rho_j^2\int_0^\infty e^{-2t}f_j(t)\,\dd t
 \le C\rho_j^2(a_j+a_j^{1/a}e^{-L_j})
 \le Ca_j\rho_j^2.
\]
For the exterior row, the inner half of the disc is
separated from the target and costs at most its dipole
moment divided by \(\rho_j^2\). On the outer annulus,
\(f\le Ca_j\rho_j^{-\beta}y_w^{-p}\); horizontal
integration gives the same \(O(a_j)\) bound.
These arguments use the all-positive strip estimates
in Theorem~\ref{thm:negative-cell-master}.

%%%%%%%%%%%%%%%%%%%%%%%%%%%%%%%%%%%%%%%%%%%%%%%%%%%%%%%%%%%%%%%%%%%%%%%%%%%%%%%%%%%%%%%%%%%%%%%%%%%%%%%%%%%%%%%%%%%%%%%%%%%%%%%%%%%%%%%%%%%%%%%%%%%%%%%%%%%%%%%%%%%%%%%%%%%%%%%%%%%%
\vskip0.12in
%%%%%%%%%%%%%%%%%%%%%%%%%%%%%%%%%%%%%%%%%%%%%%%%%%%%%%%%%%%%%%%%%%%%%%%%%%%%%%%%%%%%%%%%%%%%%%%%%%%%%%%%%%%%%%%%%%%%%%%%%%%%%%%%%%%%%%%%%%%%%%%%%%%%%%%%%%%%%%%%%%%%%%%%%%%%%%%%%%%%

We record the analogous bounds for the geometric error.
Set \(J_0=|C_0|^2\). Outside the true microdiscs,
\begin{equation}\label{eq:negative-flat-coefficient-defect}
 |J(w)-J_0|
 \le C|w|\log(e/|w|)
   +C\sum_{x_j/2<|w|<2x_j}
       \frac{\rho_j y_w+\rho_j^2}{|w-x_j|^2}.
\end{equation}
For completeness, apply the logarithm to each product
factor. If \(x_j\ge2|w|\), its change from its value
at zero is at most \(C|w|\rho_j/x_j^2\); summing gives
\(O(|w|\log(e/|w|))\). If \(x_j\le |w|/2\), the
sum at \(w\) is \(O(y_w+|w|^2)\), and the sum at zero
is \(O(|w|^2)\). The remaining bounded number of
factors have logarithms bounded by the displayed
fractions. The product is uniformly bounded outside
the microdiscs, so these logarithmic estimates imply
\eqref{eq:negative-flat-coefficient-defect}.

%%%%%%%%%%%%%%%%%%%%%%%%%%%%%%%%%%%%%%%%%%%%%%%%%%%%%%%%%%%%%%%%%%%%%%%%%%%%%%%%%%%%%%%%%%%%%%%%%%%%%%%%%%%%%%%%%%%%%%%%%%%%%%%%%%%%%%%%%%%%%%%%%%%%%%%%%%%%%%%%%%%%%%%%%%%%%%%%%%%%
\vskip0.12in
%%%%%%%%%%%%%%%%%%%%%%%%%%%%%%%%%%%%%%%%%%%%%%%%%%%%%%%%%%%%%%%%%%%%%%%%%%%%%%%%%%%%%%%%%%%%%%%%%%%%%%%%%%%%%%%%%%%%%%%%%%%%%%%%%%%%%%%%%%%%%%%%%%%%%%%%%%%%%%%%%%%%%%%%%%%%%%%%%%%%

Let \(h_j\) be the \(j\)-th displayed fraction,
multiplied by \(y_w^{-p}\), restricted to
\(x_j/2<|w|<2x_j\) and \(w\notin E_j\). Then
\begin{equation}\label{eq:negative-flat-geometric-load}
 \int y_w h_j\le C\rho_j x_j^{2-p},\qquad
 \sup_z G[h_j](z)/y_z\le C\rho_j^\beta.
\end{equation}
To check these estimates, decompose into half-annuli
about \(x_j\) of radii \(s\ge\rho_j\). There
\(h_j\lesssim(\rho_j/s)y_w^{-p}\). Its dipole moment
costs at most \(C\rho_j s^{2-p}\); its normalized
Green row costs at most \(C\rho_j s^{-p}\).
The first sum is controlled by the largest radius
\(O(x_j)\), and the second by the smallest radius
\(\rho_j\), since \(p>0\).

%%%%%%%%%%%%%%%%%%%%%%%%%%%%%%%%%%%%%%%%%%%%%%%%%%%%%%%%%%%%%%%%%%%%%%%%%%%%%%%%%%%%%%%%%%%%%%%%%%%%%%%%%%%%%%%%%%%%%%%%%%%%%%%%%%%%%%%%%%%%%%%%%%%%%%%%%%%%%%%%%%%%%%%%%%%%%%%%%%%%
\vskip0.12in
%%%%%%%%%%%%%%%%%%%%%%%%%%%%%%%%%%%%%%%%%%%%%%%%%%%%%%%%%%%%%%%%%%%%%%%%%%%%%%%%%%%%%%%%%%%%%%%%%%%%%%%%%%%%%%%%%%%%%%%%%%%%%%%%%%%%%%%%%%%%%%%%%%%%%%%%%%%%%%%%%%%%%%%%%%%%%%%%%%%%

Equations \eqref{eq:negative-flat-cell-force-moment}
and \eqref{eq:negative-flat-geometric-load} imply
the following increment bound for the cell sources,
for the geometric loads, and for their sum:
\begin{equation}\label{eq:negative-flat-local-load-increment}
        \sum_j\int |K_z(w)|\,h_j(w)\,\dd w
                       \le C(r+r^{2\beta}),
\end{equation}
with the exterior-target restriction for a true cell
source. Here the notation \(h_j\) may denote either
type of load. To verify the assertion, split the centers
into \(x_j\ge8r\), \(x_j\le r/8\), and the bounded
number of comparable centers. In the first group use
\eqref{eq:negative-flat-kernel-increment}; the sum of
the dipole moments divided by \(x_j^3\) is finite.
For the second group, both the separated target kernel
and \(P_0\) give a bound \(O(r^{2-p})=O(r^{1+\beta})\).
In the comparable group the row bounds are
\(O(a_j)=O(r^{2\beta})\). The same estimates apply to
\(y_w^{-p}\mathbf1_{E_j}\), whose dipole moment is
\(O(a_j\rho_j^2)\) and whose row is \(O(a_j)\).

Now let \(U_A\) be the one-dimensional profile determined by
\[
 U_A''=\lambda J_0U_A^{-p},\qquad U_A(0)=0,
                  \qquad U_A'(0)=A.
\]
It is equivalently defined by
\[
 y=\int_0^{U_A(y)}
       \left(A^2+\frac{2\lambda J_0}{1-p}t^{1-p}\right)^{-1/2}
       \,\dd t.
\]
In particular \(U_A(y)\asymp y\) and
\(U_A(y)/y=A+O(y^\beta)\), uniformly for \(A\) in
a compact subinterval of \((0,\infty)\).
Put \(q(z)=(v(z)-U_A(y_z))/y_z\) outside the microdiscs.
Both quotients have exterior limit \(A\), so \(q\to0\).
The exterior trace has a common modulus on the stated families. For \(r=|z|\), the baseline potential has \(C^{1,\gamma-1}\) control by the Morrey estimate and odd reflection. Cells with \(x_j\ge\sqrt r\) give harmonic potentials near zero with normalized variation \(O(\sqrt r\,a_j^{1/a})\) by reflection and \eqref{eq:negative-micro-exterior-row}. The remaining exterior contributions are bounded by \(C\sum_{x_j<\sqrt r}a_j^{1/a}\). These estimates also hold for prefixes. Choose a common \(R\) small enough that \(v\asymp y_z\) outside the cells for \(|z|<R\), using the uniform lower bound on \(A\). The nonnegative load bounds above remain valid with every source intersected with \(D_R\).

%%%%%%%%%%%%%%%%%%%%%%%%%%%%%%%%%%%%%%%%%%%%%%%%%%%%%%%%%%%%%%%%%%%%%%%%%%%%%%%%%%%%%%%%%%%%%%%%%%%%%%%%%%%%%%%%%%%%%%%%%%%%%%%%%%%%%%%%%%%%%%%%%%%%%%%%%%%%%%%%%%%%%%%%%%%%%%%%%%%%
\vskip0.12in
%%%%%%%%%%%%%%%%%%%%%%%%%%%%%%%%%%%%%%%%%%%%%%%%%%%%%%%%%%%%%%%%%%%%%%%%%%%%%%%%%%%%%%%%%%%%%%%%%%%%%%%%%%%%%%%%%%%%%%%%%%%%%%%%%%%%%%%%%%%%%%%%%%%%%%%%%%%%%%%%%%%%%%%%%%%%%%%%%%%%

Let \(f_{\mathrm{dif}}=\lambda Jv^{-p}-\lambda J_0U_A^{-p}\) on \(D_R\). The Green representation for \(v-U_A\) is the harmonic lift of its cap values minus \(G[f_{\mathrm{dif}}]\). The first density is integrable by change of variables. Its potential tends to zero at the origin by comparison with the full physical potential. At every other boundary point its local part has \(L^q\) density for some \(q>1\), since \(|F'|\) is bounded on closed patches separated from the origin. Its separated part tends to zero by dominated convergence in the Green integral. The second density is bounded by \(C y_w^{-p}\), so its potential also has zero trace. The normalized harmonic lift is a constant plus \(O(r)\) near zero by reflection across the diameter.

%%%%%%%%%%%%%%%%%%%%%%%%%%%%%%%%%%%%%%%%%%%%%%%%%%%%%%%%%%%%%%%%%%%%%%%%%%%%%%%%%%%%%%%%%%%%%%%%%%%%%%%%%%%%%%%%%%%%%%%%%%%%%%%%%%%%%%%%%%%%%%%%%%%%%%%%%%%%%%%%%%%%%%%%%%%%%%%%%%%%
\vskip0.12in
%%%%%%%%%%%%%%%%%%%%%%%%%%%%%%%%%%%%%%%%%%%%%%%%%%%%%%%%%%%%%%%%%%%%%%%%%%%%%%%%%%%%%%%%%%%%%%%%%%%%%%%%%%%%%%%%%%%%%%%%%%%%%%%%%%%%%%%%%%%%%%%%%%%%%%%%%%%%%%%%%%%%%%%%%%%%%%%%%%%%

Off the cells the source difference consists of \(\lambda(J-J_0)U_A^{-p}\) and a term bounded by \(C y_w^{-p}|q(w)|\). On each cell its absolute value is bounded by \(f+C y_w^{-p}\). These bounds make \(\int P_0|f_{\mathrm{dif}}|\) finite. The kernel estimates with \(\alpha=0\) and the load bounds give \(\int |K_z||f_{\mathrm{dif}}|\lesssim r^\beta+r+r^{2\beta}\) for exterior targets, since \(q\) is bounded. As \(q\to0\), the constant in the normalized harmonic lift is therefore \(\int P_0 f_{\mathrm{dif}}\). Subtracting this flux and using \eqref{eq:negative-flat-coefficient-defect}--\eqref{eq:negative-flat-local-load-increment} gives
\begin{equation}\label{eq:negative-flat-inlet-bootstrap}
 |q(z)|\le C(r+r^{2\beta})
   +C\int_{D_R\setminus\cup_j E_j}
                     |K_z(w)|y_w^{-p}|q(w)|\,\dd w .
\end{equation}
Initially \(q\) is bounded. Applying
\eqref{eq:negative-flat-weighted-increment} with
\(\alpha=0\) gives \(q=O(r^\beta)\).
Apply it again with \(\alpha=\beta\). If \(p>1/2\),
this gives \(q=O(r^{2\beta})\); if \(p<1/2\), it
gives \(q=O(r)\). At \(p=1/2\) it first gives
\(q=O(r\log(e/r))\), hence \(q=O(r^{3/4})\);
a further application gives \(q=O(r)\).
Consequently in every case
\[
       |v(z)/y_z-U_A(y_z)/y_z|\le C(r+r^{2\beta})
                    \quad(z\notin\cup_j E_j).
\]
For inlets with \(x_n<R/4\), the inlet semicircle has \(r\asymp x_n\) and \(y_z\le\rho_n\). Thus \(|v(z)/y_z-A|\le C(x_n+a_n)\), uniformly in angle. When \(x_n\ge R/4\), bounded mass bounds the inlet moment and the exterior slope, and \(x_n+a_n\ge R/4\).
Average with the probability measure
\((2/\pi)\sin^2\theta\,\dd\theta\) to obtain
\eqref{eq:negative-flat-sharp-inlet}.
\end{proof}
For the matched loss law let \(D_n\) denote the local product coefficient in \eqref{eq:negative-tuned-coefficient}. Uniformly for full or truncated products containing \(E_n\),
\begin{equation}\label{eq:flat-coefficient-comparison}
 |D_n-D_p|\lesssim x_n .
\end{equation}
In the logarithm of the squared modulus of the product, the derivative of a factor with \(j<n\) on \([0,x_n]\) is bounded by \(C\rho_j^2/x_j^3\), a summable sequence. The \(n\)-th factor at zero differs from one by \(O(x_n^2)\); the factors with \(j>n\) have total logarithmic modulus \(O(\sum_{j>n}x_j^2)\) both at zero and at \(x_n\), using \(|x_n-x_j|\gtrsim x_n\). These bounds imply \eqref{eq:flat-coefficient-comparison}.

\subsection{Local energy and the logarithmic depth threshold}

The local form \(b_Q(V)\) of \cref{thm:sharp-boundary-transition} can be computed on shrinking
patches in the preimage, with the pulled-back potential
\(V_D=|F'|^2(V\circ F)\). This follows from conformal
invariance of the Dirichlet integral and from the two
nested systems of neighborhoods defined by the boundary
homeomorphism. The test spaces can be treated by taking
suprema over smooth compactly supported tests and using
density for the bounded form.

\begin{proposition}\label{prop:negative-local-energy-slopes}
On a flat product of \cref{thm:negative-micro-exact-drop}, let \(u\) be a stable solution. At \(Q\),
\begin{equation}\label{eq:negative-local-energy-slopes}
 b_Q(\lambda p u^{-a})
 \asymp
 \limsup_{n\to\infty}\frac{a_n}{a_n+s_n(u)^a}.
\end{equation}
In particular, a strict exterior slope margin implies
\(b_Q(\lambda p u^{-a})=0\). At a saturated solution,
the local energy concentration vanishes exactly when
\(s_n(u)/a_n^{1/a}\to\infty\).
\end{proposition}

\begin{proof}
For \(|t-L_n|\le1\), the convolution formula and
\eqref{eq:negative-micro-plateau-block} give
\[
 0\le Y_n(t)-s_n(u)
 \le \tfrac12\int_t^\infty f_n(s)\,\dd s
 \lesssim a_n^{1/a}.
\]
The upper estimate in
\eqref{eq:negative-micro-plateau-profile}, with fixed
half-width, gives \(W_n(t,\theta)\lesssim
(s_n(u)+a_n^{1/a})\sin\theta\).
On a fixed interior ball in the rescaled variable
\((z-x_n)/b_n\) with radius and height comparable to one
and \(|t-L_n|\le1\), the two-sided envelope then gives
\[
 V_D(z)=\frac{pH_n(t,\theta)}{r^2 W_n(t,\theta)^a}
 \gtrsim b_n^{-2}\frac{a_n}{a_n+s_n(u)^a}.
\]
A fixed nonzero bump scaled to that ball gives the lower
bound in \eqref{eq:negative-local-energy-slopes}.

Write \(q_0\) for the limsup on the right.
If \(q_0\) is below a sufficiently small fixed threshold,
choose the half-width in
\eqref{eq:negative-micro-plateau-profile} large and use
\(Y_n\ge s_n(u)\). On all sufficiently late discs,
\(W_n\ge c s_n(u)\sin\theta\), since
\(a_n^{1/a}/s_n(u)\) is small. Thus
\[
 d_D^2 V_D\lesssim \frac{a_n}{s_n(u)^a}\sin^{2-a}\theta
 \quad\text{on }E_n.
\]
Outside the discs, the exterior trace is positive, so
\(u\circ F\gtrsim y\) near zero. In the flat-tangent
case \(d_D^2V_D\lesssim y^{1-p}\) there, which tends to zero. Applying the Hardy inequality
in the smooth preimage gives the upper bound by \(Cq_0\).
If \(q_0\) is above the fixed threshold, stability gives
\(b_Q(V)\le1\lesssim q_0\). The last assertions follow
from the exact loss law.
\end{proof}

\begin{corollary}\label{cor:flat-depth-threshold}
Suppose \(p\ge 1/2\), \(\mu_n\to\mu>0\), and \(A^a=D_p\mu\). Set
\[
 S_n=s_n^a/a_n,\qquad R_n=(\mu-\mu_n)/a_n-\Theta_p\log(1/a_n).
\]
Then
\begin{equation}\label{eq:flat-depth-threshold}
 S_n-D_n\Theta_p\log(1+S_n)=D_n R_n+O(1).
\end{equation}
In particular \(R_n\) is bounded below. On subsequences with \(R_n\to\infty\), \(S_n\sim D_pR_n\); on subsequences with \(R_n\) bounded above, \(S_n\) is bounded. For a stable solution, \(b_Q=0\) if and only if \(R_n\to\infty\).
\end{corollary}

\begin{proof}
The refined inlet estimate and \eqref{eq:flat-coefficient-comparison} give \(Y_n(0)^a=A^a+O(a_n)\) and \(D_n-D_p=O(a_n)\) since \(x_n\lesssim a_n\). In the signed loss identity, use
\[
 \log(1+(s_n^a+a_n)^{-1})
 =\log(1/a_n)-\log(1+S_n)+\log(1+s_n^a+a_n).
\]
The last logarithm is bounded on bounded masses. This gives \eqref{eq:flat-depth-threshold}. The function \(S-D_n\Theta_p\log(1+S)\) is uniformly bounded below for \(S\ge0\) and is asymptotic to \(S\) at infinity. The last assertion follows from \cref{prop:negative-local-energy-slopes}.
\end{proof}

\subsection{A continuation criterion at the minimum}

Suppose \(\mu_n\to\mu>0\). If the endpoint \(U_*\) has positive slopes at every smooth off-pole point and \(\|B_*\|_E<1\), then \(A(U_*)^a=D_p\mu\). A strict exterior margin would give \(U_*\ge cK\) by \cref{thm:negative-micro-exact-drop}, and hence a finite linearized Martin row. The spectral transfer of \cref{thm:negative-two-space-radius} implies \(r_X(A_*)\le\|B_*\|_E<1\). The Green map is analytic in a relative ball by domination by \(K^{-p}\), as in \cref{lem:negative-cone-fold}. The implicit function theorem would continue positive solutions below \(m_*\).

%%%%%%%%%%%%%%%%%%%%%%%%%%%%%%%%%%%%%%%%%%%%%%%%%%%%%%%%%%%%%%%%%%%%%%%%%%%%%%%%%%%%%%%%%%%%%%%%%%%%%%%%%%%%%%%%%%%%%%%%%%%%%%%%%%%%%%%%%%%%%%%%%%%%%%%%%%%%%%%%%%%%%%%%%%%%%%%%%%%%
\vskip0.32in
%%%%%%%%%%%%%%%%%%%%%%%%%%%%%%%%%%%%%%%%%%%%%%%%%%%%%%%%%%%%%%%%%%%%%%%%%%%%%%%%%%%%%%%%%%%%%%%%%%%%%%%%%%%%%%%%%%%%%%%%%%%%%%%%%%%%%%%%%%%%%%%%%%%%%%%%%%%%%%%%%%%%%%%%%%%%%%%%%%%%

\section{\bf Saturated endpoints and finite-product approximation}\label{sec:negative-realizations}

We give two choices of the depths \(\mu_n\). Both yield a strict energy gap at a saturated endpoint. The first gives a sharp Martin-norm response, and the second has positive energy concentration and a quantitative truncation law.

\begin{corollary}\label{cor:negative-micro-realized-marginal}
For every \(0<p<1\) and \(\lambda>0\) there is a modified
flat-tangent domain of Theorem~\ref{thm:negative-micro-exact-drop}
with unique positive solution at each admissible mass, including
the endpoint, and with
\[
 A(U_{m_*})=(D_p\mu)^{1/(1+p)}.
\]
The minimum-mass solution loses its positive Martin lower bound
only at the exceptional point \(Q\). It has a positive relative
boundary trace at every smooth boundary point. Its energy linearization is compact and has norm strictly less than one. The Green correction has a right derivative in energy, and the Martin norm response is bounded below by \(c(m-m_*)^{1/a}\).
\end{corollary}

\begin{proof}
Keep the cap \(D\) and kernel \(k\) fixed and vary the centers and depths in \cref{sec:flat-geometry}. The uniform geometry in the proof of \cref{thm:negative-micro-nonlinear} makes the estimates independent of the chosen centers.

%%%%%%%%%%%%%%%%%%%%%%%%%%%%%%%%%%%%%%%%%%%%%%%%%%%%%%%%%%%%%%%%%%%%%%%%%%%%%%%%%%%%%%%%%%%%%%%%%%%%%%%%%%%%%%%%%%%%%%%%%%%%%%%%%%%%%%%%%%%%%%%%%%%%%%%%%%%%%%%%%%%%%%%%%%%%%%%%%%%%
\vskip0.12in
%%%%%%%%%%%%%%%%%%%%%%%%%%%%%%%%%%%%%%%%%%%%%%%%%%%%%%%%%%%%%%%%%%%%%%%%%%%%%%%%%%%%%%%%%%%%%%%%%%%%%%%%%%%%%%%%%%%%%%%%%%%%%%%%%%%%%%%%%%%%%%%%%%%%%%%%%%%%%%%%%%%%%%%%%%%%%%%%%%%%

The actual-force bound
\eqref{eq:negative-micro-exterior-row} and the bounded
derivative baseline in the proof of
Theorem~\ref{thm:negative-micro-nonlinear} give
\[
 u\circ F\ge(m-C)k\quad\hbox{outside }\bigcup_n E_n
\]
with \(C\) independent of \(\mu_n\). The linear inverse
potential has norm at most \(C(1+\sup_n\mu_n)\) by the
shell estimate in the concentrating example. For \(\mu_n\in[\mu/2,\mu]\),
the threshold in Theorem~\ref{thm:negative-branch} and the
large-mass induction of
Theorem~\ref{thm:negative-micro-nonlinear} therefore give
a number \(M=M(\mu)\) independent of the centers. The
endpoint is below \(M\), and all masses above \(M\) have
a unique positive solution. By the exact drop, every
admissible mass is at least
\((D_p\mu)^{1/a}/\kappa\), since \(A(u)\le m\kappa\).
The factor \(|C_0|\) is bounded below by one uniformly.
Thus one can choose \(\mu\) so large that \(m-C\) is
positive and, outside \(\bigcup_n E_n\), the pulled-back
linearized potential
\[
 V_D=\lambda p|F'|^2(u\circ F)^{-1-p}
\]
satisfies \(d_D^2 V_D\le\eta_0\), where \(\eta_0>0\) is
any prescribed small number. This uses \(|F'|\lesssim1\)
off the disks, \(k\gtrsim d_D\), and \(1+p<2\).

%%%%%%%%%%%%%%%%%%%%%%%%%%%%%%%%%%%%%%%%%%%%%%%%%%%%%%%%%%%%%%%%%%%%%%%%%%%%%%%%%%%%%%%%%%%%%%%%%%%%%%%%%%%%%%%%%%%%%%%%%%%%%%%%%%%%%%%%%%%%%%%%%%%%%%%%%%%%%%%%%%%%%%%%%%%%%%%%%%%%
\vskip0.12in
%%%%%%%%%%%%%%%%%%%%%%%%%%%%%%%%%%%%%%%%%%%%%%%%%%%%%%%%%%%%%%%%%%%%%%%%%%%%%%%%%%%%%%%%%%%%%%%%%%%%%%%%%%%%%%%%%%%%%%%%%%%%%%%%%%%%%%%%%%%%%%%%%%%%%%%%%%%%%%%%%%%%%%%%%%%%%%%%%%%%

We need a common inlet error as \(x_n\to0\) for bounded masses and depths in \([\mu/2,\mu]\). The baseline potential divided by \(y\) has a common modulus \(O(x_n^{\gamma-1})\) from zero to the inlet. Indeed \(|F'|\) is bounded on the baseline, so the physical Morrey estimate transfers to preimage balls, including boundary balls by Whitney summation. Odd reflection and the Newtonian kernel bound give \(C^{1,\gamma-1}\) control up to the flat arc, with the harmonic remainder uniformly controlled by the baseline potential estimate. For \(x_j\ge\sqrt{x_n}\), the potential of \(E_j\) is harmonic in \(B_{x_j/3}^+\), with normalized norm \(O(a_j^{1/a})\) by \eqref{eq:negative-micro-exterior-row}. Reflection bounds its normalized variation by \(C\sqrt{x_n}a_j^{1/a}\). For \(x_j<\sqrt{x_n}\) the same row bound and lacunarity give total contribution \(O(x_n^{(1-p)/a})\) on the inlet and in the exterior limit. Adding the harmonic mass term gives \(Y_n(0)-A(u)=O(x_n^{\gamma-1}+\sqrt{x_n}+x_n^{(1-p)/a})\). These estimates hold also for retained prefixes, since the baseline bounds are uniform and the sums are restricted to the retained cells. Applying \eqref{eq:negative-quantitative-cell-loss} proves
\begin{equation}\label{eq:negative-micro-uniform-error}
 \left|s_n(u)^a-A(u)^a+D_p\mu_n\right|\le\omega(x_n),\qquad
 \omega(r)\longrightarrow0,
\end{equation}
with a common modulus. For a prefix the coefficient and the exterior trace are those of its own retained product.

Now take, for example, \(\mu_n=\mu(1-1/(n+1))\).
Choose the lacunary centers sufficiently small so that
\(\omega(x_n)\le \frac12(\inf D_p)(\mu-\mu_n)\)
and \(a_n/(\mu-\mu_n)\) is smaller than a prescribed constant
for every \(n\). The infimum here is positive over the family
because \(|C_0|\ge1\).
The exact drop and \eqref{eq:negative-micro-uniform-error}
give \(s_n(u)^a\ge \frac12 D_p(\mu-\mu_n)\)
for every solution of mass \(m\le M\).
With the fixed comparison half-width \(M_0\) in
\eqref{eq:negative-micro-plateau-profile} chosen large,
the same smallness of \(a_n/(\mu-\mu_n)\) gives
\(W_n\gtrsim s_n(u)\sin\theta\) on each whole microdisc.
The centers can also be taken so small that the fixed
outer extension for this estimate is available for every \(n\).
As \(d_D(z)=r\sin\theta\) on \(E_n\), we conclude there that
\[
 d_D(z)^2 V_D(z)
 \lesssim \frac{a_n}{\mu-\mu_n}\sin^{\,1-p}\theta
 \le\eta_0
\]
after choosing that smallness sufficiently strong.

%%%%%%%%%%%%%%%%%%%%%%%%%%%%%%%%%%%%%%%%%%%%%%%%%%%%%%%%%%%%%%%%%%%%%%%%%%%%%%%%%%%%%%%%%%%%%%%%%%%%%%%%%%%%%%%%%%%%%%%%%%%%%%%%%%%%%%%%%%%%%%%%%%%%%%%%%%%%%%%%%%%%%%%%%%%%%%%%%%%%
\vskip0.12in
%%%%%%%%%%%%%%%%%%%%%%%%%%%%%%%%%%%%%%%%%%%%%%%%%%%%%%%%%%%%%%%%%%%%%%%%%%%%%%%%%%%%%%%%%%%%%%%%%%%%%%%%%%%%%%%%%%%%%%%%%%%%%%%%%%%%%%%%%%%%%%%%%%%%%%%%%%%%%%%%%%%%%%%%%%%%%%%%%%%%

The Hardy inequality on the fixed \(D\), and conformal
invariance of Dirichlet energy, now give a uniform strict
stability inequality for all solutions of mass at most \(M\).
Theorem~\ref{thm:negative-branch} gives uniqueness among stable
solutions, and the large-mass argument covers \(m>M\).
The exterior and individual cell bounds give positive slopes at every smooth off-pole point. The stability gap and the continuation criterion above imply saturation. Since \(a_n/(\mu-\mu_n)\to0\), the local form criterion gives \(b_Q=0\). All other local forms vanish by the positive slopes and the estimate \(V_*^D\lesssim k^{-a}\) near the smooth pole. Thus the energy operator is compact by \cref{thm:negative-two-space-radius,cor:negative-essential-local-forms}. The strict gap gives the energy derivative by \cref{thm:negative-energy-endpoint}. The lower Martin bound follows by taking the difference of the limiting center slopes in \eqref{eq:negative-micro-exact-drop}, using \(A(U_{m_*+\delta})\ge A(U_{m_*})+\kappa\delta\) and \(\partial_y k(x_n)\to\kappa\).
\end{proof}

\begin{lemma}\label{lem:flat-collar-response}
Fix \(\ell\in(\max(p,1/2),1)\). On the full product or a prefix, let \(v_1,v_2\) be maximal solutions at \(m_1,m_2=m_1+\delta\), \(\delta>0\), and set \(V_1^D=\lambda pJv_1^{-a}\). Suppose \(d_D^2 V_1^D\le\eta_0\) with \(\eta_0\) sufficiently small and \(V_1^D\lesssim k^{-a}\) off the retained cells. Then
\begin{equation}\label{eq:flat-collar-response}
 0\le v_2-v_1-\delta k\le C\delta d_D^\ell,\qquad
 \kappa\delta\le A(v_2)-A(v_1)\le C\delta .
\end{equation}
The constants are uniform if the off-cell bound is uniform. For a prefix we take \(J=J_N\). Here \(A(v)\) denotes the exterior trace of the solution in the preimage.
\end{lemma}

\begin{proof}
Let \(C_{12}\) be the preimage secant potential. Branch order and the mean value formula give \(0\le C_{12}\le V_1^D\) and \(Z=v_2-v_1-\delta k\ge0\). The physical Green corrections have \(L^q\) sources for some \(q>1\) by \cref{thm:negative-full-representation}, hence finite Dirichlet energy by the planar Sobolev embedding. Conformal invariance gives \(Z\in H^1_0(D)\). It solves \((-\Delta-C_{12})Z=\delta C_{12}k\). On the cells \(k\lesssim d_D\), and away from the cells \(V_1^D k^2\lesssim k^{1-p}\). Thus \(\int_D V_1^D k^2\) is uniformly bounded, using the smooth-pole estimate and \(1-p<1\). The Hardy inequality gives a coercive form for \(\eta_0\) small. Testing with \(Z\) and using Cauchy--Schwarz in the potential form gives \(\|Z\|_{H^1_0(D)}\lesssim\delta\). Interior estimates give \(Z\lesssim\delta\) on the inner crosscut of a fixed smooth collar, since \(C_{12}\le\eta_0 d_D^{-2}\) and \(k\) is bounded there.

In the collar, \((-\Delta-C_{12})d_D^\ell\gtrsim d_D^{\ell-2}\), by \(\eta_0<\ell(1-\ell)\) and boundedness of \(\Delta d_D\). The right side dominates \(C_{12}k\) up to a constant, since \(C_{12}k\lesssim d_D^{-1}\) on the cells and \(C_{12}k\lesssim k^{-p}\lesssim d_D^{-p}\) elsewhere. As \(\ell>1/2\), \(d_D^\ell\) has finite collar energy; the positive-part test of \(Z-C\delta d_D^\ell\) has zero trace on both collar boundaries for \(C\) large. Coercivity gives the pointwise bound, extended to the interior by the interior estimates.

By the flux formula \eqref{eq:negative-kernel-exterior-identity}, the slope increment is
\[
 \kappa\delta+\int_D P_Q C_{12}(\delta k+Z).
\]
Near zero off the cells, \(P_Q\lesssim y/|w|^2\), \(k\asymp y\) and \(C_{12}\lesssim y^{-a}\); the integral there is bounded by a constant times \(\delta\int_0^{r_0}(r^{-p}+r^{\ell-p-1})\dd r\). On \(E_j\) the bound \(C_{12}\lesssim y^{-2}\) gives at most \(C\delta x_j^{-2}\int_{E_j}y^{\ell-1}\lesssim\delta x_j^{2\ell}\). Away from zero, \(P_Q\lesssim d_D\) and \(k\gtrsim d_D\); the remaining integral is \(O(\delta)\). This proves \eqref{eq:flat-collar-response}.
\end{proof}

\begin{corollary}\label{cor:negative-micro-matching-response}
The centers in Corollary~\ref{cor:negative-micro-realized-marginal}
can be chosen so that the response has matching upper and
lower bounds:
\[
 c\delta^{1/(1+p)}
 \le \|U_{m_*+\delta}-U_{m_*}\|_K
 \le C\delta^{1/(1+p)}
\]
for all sufficiently small \(\delta>0\).
\end{corollary}

\begin{proof}
Retain the notation of the tuned construction and put
\(g_n=\mu-\mu_n\), \(u_*=U_{m_*}\circ F\), and
\(w=(U_{m_*+\delta}-U_{m_*})\circ F\).
All estimates concern \(m_*+\delta\le M\). Apply \cref{lem:flat-collar-response} at \(v_1=u_*\) and write \(b=\ell\) and \(V_\delta=C_{12}\). We have \(Z=w-\delta k\le C\delta d_D^b\) and \(0\le V_\delta\le V_*^D=\lambda p J u_*^{-a}\). By \cref{thm:negative-energy-endpoint} and conformal invariance, \(Z/\delta\) converges in energy to the solution \(Z_*\) of \((-\Delta-V_*^D)Z_*=V_*^Dk\) in \(D\). Dominated convergence in the flux formula, with the bounds in the collar proof and local elliptic convergence, gives
\begin{equation}\label{eq:flat-exterior-derivative}
 A'_*(+)=\kappa+\int_D P_Q V_*^D(k+Z_*)>0 .
\end{equation}
The exact loss law thus also gives
\[
 \lim_{\delta\downarrow0}\delta^{-1/a}\lim_{n\to\infty}
       s_n(U_{m_*+\delta})
   = \big(a A(U_*)^{a-1} A'_*(+)\big)^{1/a}.
\]

We estimate the normalized potential of
\(-\Delta Z=V_\delta w\).
Off the disks, the preceding claim and \(u_*\gtrsim k\)
give the bound \(C\delta(d_D^{-p}+d_D^{b-a})\) for
this density. Both exponents are integrable in the smooth
distance potential criterion, since \(b>p\).
On \(E_n\) one has \(w\lesssim\delta v^b\), where
\(v=\operatorname{Im}z\), and \(V_*^D\lesssim (a_n/g_n)v^{-2}\).
For targets outside \(E_n\) the contribution of its inner
half-disc of radius \(\rho_n/2\) is bounded by
\(C\delta(a_n/g_n)\rho_n^{b-1}\), using the separated
half-plane kernel and
\(\int_{E_n}v^{b-1}\dd x\dd v\lesssim\rho_n^{b+1}\).
On the outer annulus the sharper pointwise bound
\(V_*^Dw\lesssim\delta(a_n/g_n)\rho_n^{-(1-p)}v^{b-a}\)
follows from \(|F'|^2\lesssim a_n\rho_n^{-(1-p)}\) and
\(u_*\gtrsim g_n^{1/a}v\).
Horizontal integration of the half-plane kernel gives the same
normalized bound on this annulus, since
\(\int_0^{\rho_n}v^{b-a}\dd v\lesssim\rho_n^{b-p}\).
Outside the comparison half-disc the estimates extend
by harmonic comparison with \(k\), as in
Theorem~\ref{thm:negative-micro-nonlinear}.
Take the centers sufficiently small to make
\(\sum_n\rho_n^{b-p}/g_n<\infty\).
Then the contribution of all other disks and the baseline
is bounded by \(C\delta\) at targets in any disk, and the
whole potential has this bound at exterior targets.
In particular,
\begin{equation}\label{eq:negative-micro-exterior-lipschitz}
 0\le A(U_{m_*+\delta})-A(U_{m_*})\le C\delta .
\end{equation}

Inside \(E_n\), denote by \(M_i\) the supremum of \(w/k\)
on shell \(i\) and by \(M^{(n)}\) the supremum on the
whole disk. The lower bound for the endpoint first moment is
\[
 Y_{*,n}(t)^a\gtrsim g_n+a_n(L_n-t)_+ .
\]
Indeed, its slope power is bounded below by \(c g_n\),
and integration of the first-mode Jensen lower bound in the proof of \cref{thm:negative-cell-master} gives the second term. On each shell the profile gives
\(u_*\gtrsim kY_{*,n}\).
If \(t_h\) is the largest time on an outer shell \(h<i\),
the normalized linear shell estimate therefore gives
\[
 M_i\le C\delta+C(a_n/g_n)M^{(n)}
       + C\sum_{h<i}
          \frac{a_n M_h}{g_n+a_n(L_n-t_h)_+}.
\]
Here the same and deeper shells are bounded together by
\(C(a_n/g_n)M^{(n)}\), using the linear inner-shell
row estimate in Proposition~\ref{prop:negative-micro-essential-zero}.
The earlier coefficients have sum at most
\(C[1+\log(1+\mu/g_n)]\). This follows by comparison with
the integral of \(a_n/(g_n+a_n(L_n-t))\), using the fixed
dyadic shell width and \(a_n/g_n\lesssim1\).
Discrete Gronwall therefore gives, with a fixed \(B>0\),
\[
 M^{(n)}
 \le C g_n^{-B}\bigl(\delta+(a_n/g_n)M^{(n)}\bigr).
\]
The bound \(0\le w\le(m_*+\delta)k\) ensures \(M^{(n)}<\infty\).
Decrease the centers so that the coefficient of \(M^{(n)}\)
on the right is at most \(1/2\). This proves
\(M^{(n)}\le C\delta g_n^{-B}\).

For the late disks, the inlet modulus just used in \eqref{eq:negative-micro-uniform-error} and the partial loss \eqref{eq:negative-cell-partial-loss} give a common error modulus \(\widetilde\omega(x_n)\to0\) for
\[
 \sup_{t\ge0}\big|Y_n(t)^a-A(u)^a+D_p a_n\min(t,L_n)\big|.
\]
The same modulus can dominate \(\sup_{t,\theta}|W_n(t,\theta)/\sin\theta-Y_n(t)|\). Indeed, multiply \eqref{eq:negative-micro-plateau-profile} by \(Y_n\le CM\), then let \(x_n\to0\) at fixed strip width before sending the width to infinity.

Choose a fixed \(\zeta>0\) such that
\(1-\zeta B\ge1/a\), and take the centers so small that
\(\widetilde\omega(x_n)\lesssim g_n^N\) for some
\(N>1/\zeta\), simultaneously with the preceding size
conditions. These choices can be made in a single recursion.
First fix \(b\) and \(\eta_0\), then \(\mu\) and the uniform mass
bound \(M\). The error moduli and the shell constants are
uniform over the lacunary configurations satisfying the
earlier bounds, so \(B,\zeta,N\) can be fixed before the final
choice of the centers. With \(g_n=\mu/(n+1)\) fixed, each
condition is attained by making \(x_n\) sufficiently small
relative to \(x_{n-1}\).
On \(g_n\ge\delta^\zeta\) the shell Gronwall
estimate gives \(M^{(n)}\lesssim\delta^{1/a}\).
On \(g_n<\delta^\zeta\), subtract the two approximate
\(a\)-th power laws for \(Y_n\); the common deterministic
term cancels. By \eqref{eq:negative-micro-exterior-lipschitz},
the difference of the initial powers is \(O(\delta)\).
The error is \(O(g_n^N)=O(\delta)\).
The \(1/a\)-H\"older continuity of the nonnegative \(a\)-th
root gives \(|Y_{\delta,n}-Y_{*,n}|\lesssim\delta^{1/a}\).
The additive angular approximation and
\(k\asymp r\sin\theta\) give the same bound on \(w/k\).
The exterior estimate covers the remaining targets.
This proves the upper bound. The lower bound was proved in \cref{cor:negative-micro-realized-marginal}.
\end{proof}

\subsection{Limits of finite products}

Let \(m_N\) be the minimum mass in \(F_N(D)\), and write \(v_{N,m}=U_{N,m}\circ F_N\) and \(v_m=U_m\circ F\). Let \(A_N(m)\) denote the exterior trace of \(v_{N,m}\) at zero and \(A(m)\) that of \(v_m\).

\begin{proposition}\label{prop:negative-product-truncation}
For a flat product with \(\mu_n\to\mu>0\), \(m_N\uparrow m_*\). Both \(v_{N,m_N}\) and \(v_{N,m_*}\) converge to \(v_{m_*}\) in \(C^2_{\mathrm{loc}}(D)\), with convergence of their exterior slopes. For either of these sequences \(v_N\),
\[
 \liminf_{N\to\infty}\|v_N-v_{m_*}\|_k
 \ge \kappa^{-1}\big[A(m_*)-(A(m_*)^a-D_p\mu)^{1/a}\big]>0 .
\]
\end{proposition}

\begin{proof}
Since \(J_N\uparrow J\), lower-barrier iteration yields \(m_N\uparrow\bar m\le m_*\), and \(v_{N,m_*}\ge v_{m_*}\) decreases with \(N\). For the two sequences, stability and \cref{lem:negative-stable-distance-seed}, pulled back using Koebe's distance theorem, give
\[
 v_{N,m}\gtrsim (d_D |F_N'|)^\gamma,\qquad
 f_{N,m}:=\lambda J_N v_{N,m}^{-p}\lesssim |F'|^\gamma d_D^{-\sigma}
 \quad(m=m_N\text{ or }m_*).
\]
Local elliptic compactness follows from the uniform positive lower bound on compacts and \(v_{N,m}\le m_* k\). The displayed majorant is integrable against every fixed interior Green row and against \(P_Q\). Off the cells this follows from \(d_D^{-\sigma}\) and \(P_Q\lesssim y/|w|^2\) near zero. On \(E_j\) its first moment in \(y\) is \(O(a_j^{1/a}\rho_j^2)\), by polar integration using \(\xi\gamma+\sigma=1\). Since \(P_Q\lesssim y/x_j^2\) there, the resulting sum is finite. A fixed interior Green row is \(O(y)\) on the late cells, and the remaining finitely many cells have bounded derivatives. In particular, for the limiting stable solution \(f_*=\lambda Jv_{m_*}^{-p}\),
\begin{equation}\label{eq:flat-prefix-near-flux}
 \int_{B_R(0)\cap D}P_Q f_* \lesssim R^{\beta/a}
 \qquad (R\downarrow0).
\end{equation}
Indeed the baseline integral is \(O(R^{1-\sigma})\); the cells that meet \(B_R(0)\) have \(x_j\lesssim R\) and their first-moment sum weighted by \(x_j^{-2}\) is \(O(R^{2+2\beta/a})\).

%%%%%%%%%%%%%%%%%%%%%%%%%%%%%%%%%%%%%%%%%%%%%%%%%%%%%%%%%%%%%%%%%%%%%%%%%%%%%%%%%%%%%%%%%%%%%%%%%%%%%%%%%%%%%%%%%%%%%%%%%%%%%%%%%%%%%%%%%%%%%%%%%%%%%%%%%%%%%%%%%%%%%%%%%%%%%%%%%%%%
\vskip0.12in
%%%%%%%%%%%%%%%%%%%%%%%%%%%%%%%%%%%%%%%%%%%%%%%%%%%%%%%%%%%%%%%%%%%%%%%%%%%%%%%%%%%%%%%%%%%%%%%%%%%%%%%%%%%%%%%%%%%%%%%%%%%%%%%%%%%%%%%%%%%%%%%%%%%%%%%%%%%%%%%%%%%%%%%%%%%%%%%%%%%%

Dominated convergence passes the representation and exterior flux for each local limit. The limit at mass \(\bar m\) is positive in \(D\); its Green representation and the continuous zero trace of the physical potential give a positive one-pole solution on \(\Omega\). Thus \(\bar m=m_*\). The limits are stable by testing on compacts and are identified by \cref{thm:negative-branch}. This proves local and exterior slope convergence. For fixed \(N\), the boundary near zero contains no retained cells; the quotient \(v_N/k\) has unrestricted limit \(A_N/\kappa\) there, by the ordinary \(C^1\) estimate. At the omitted centers \(x_j\) its traces converge to that value as \(j\to\infty\). The full solution has corresponding traces tending to \((A(m_*)^a-D_p\mu)^{1/a}/\kappa\) by \eqref{eq:negative-micro-exact-drop}. Evaluate the normalized difference at these boundary limits and then send \(N\to\infty\).
\end{proof}

\subsection{Positive concentration and the scale of the mass error}

\begin{theorem}\label{thm:flat-critical-depth}
For every \(1/2\le p<1\) and \(\lambda>0\), one can choose a flat product with
\[
 x_{n+1}\le x_n^5,\qquad
 \mu-\mu_n=\Theta_p a_n\log(1/a_n)+K_0a_n ,
\]
where \(K_0>0\) is a sufficiently large constant, for which the minimum-mass solution is saturated and has positive slopes at all smooth off-pole points. It satisfies
\[
 0<b_Q(V_*)\le b(V_*)<1,\qquad s_n(U_*)^a\asymp a_n,\qquad
 m_*-m_N\asymp a_N .
\]
\end{theorem}

\begin{proof}
Fix the cap \(D\) and let \(v_*=v_{m_*}\). Choose \(\eta_0\) small enough for \cref{lem:flat-collar-response}. Because \(|C_0|\ge1\), \(D_p\) has a fixed positive lower bound \(D_{\min}=a\lambda c_p/2\). Choose \(\mu\) large and set \(\underline m=(D_{\min}\mu)^{1/a}/(2\kappa)\). The exterior estimate \(v\ge(m-C)k\) holds uniformly off the cells, with \(C\) independent of the depths. Thus for \(m\ge\underline m\) one has \(v\gtrsim k\) and \(d_D^2\lambda pJv^{-a}\le\eta_0\) off the cells if \(\mu\) is large, also for prefixes with \(J_N\). The full endpoint has \(m_*\ge 2\underline m\) by \cref{thm:negative-micro-exact-drop}. The branch threshold and the uniform bound for the inverse-power row give an upper mass cap \(M(\mu)\) for all endpoints with depths in \([\mu/2,\mu]\).

Apply the matching and inlet estimates on \([\underline m,M(\mu)]\), where the exterior slopes are uniformly bounded above and bounded away from zero. Fix \(K_0\) large compared with these uniform errors. Choose the centers with \(x_{n+1}\le x_n^5\) and \(x_1\) small enough for every inlet to lie in the fixed half-disc of the inlet estimate, for the strip comparisons to hold on all cells, and for \(\mu_n\in[\mu/2,\mu]\). After the product is fixed, \(m_N\uparrow m_*\) by \cref{prop:negative-product-truncation}, so \(m_N\ge\underline m\) for large \(N\).

%%%%%%%%%%%%%%%%%%%%%%%%%%%%%%%%%%%%%%%%%%%%%%%%%%%%%%%%%%%%%%%%%%%%%%%%%%%%%%%%%%%%%%%%%%%%%%%%%%%%%%%%%%%%%%%%%%%%%%%%%%%%%%%%%%%%%%%%%%%%%%%%%%%%%%%%%%%%%%%%%%%%%%%%%%%%%%%%%%%%
\vskip0.12in
%%%%%%%%%%%%%%%%%%%%%%%%%%%%%%%%%%%%%%%%%%%%%%%%%%%%%%%%%%%%%%%%%%%%%%%%%%%%%%%%%%%%%%%%%%%%%%%%%%%%%%%%%%%%%%%%%%%%%%%%%%%%%%%%%%%%%%%%%%%%%%%%%%%%%%%%%%%%%%%%%%%%%%%%%%%%%%%%%%%%

With \(D_{p,N}\) the prefix coefficient, write \(D_{N,n}\) for its coefficient \(D_n\) from \eqref{eq:negative-tuned-coefficient}. For \(n\le N\), put \(S=s_{N,m,n}^a/a_n\). The signed loss identity, \eqref{eq:flat-coefficient-comparison}, and the refined inlet estimate give
\begin{equation}\label{eq:flat-prefix-matched}
 S-D_{N,n}\Theta_p\log(1+S)
 =\frac{A_N(m)^a-D_{p,N}\mu}{a_n}+D_{N,n}K_0+O(1).
\end{equation}
To see the cancellation, substitute the depths in the loss identity and use the logarithm splitting in \cref{cor:flat-depth-threshold}. The term \(D_{N,n}\Theta_pa_n\log(1/a_n)\) cancels. Replacing \(D_{N,n}\mu\) by \(D_{p,N}\mu\) costs \(O(x_n)=O(a_n)\), and \(Y_n(0)^a=A_N(m)^a+O(a_n)\). The same equation holds on the full product with its coefficients.

%%%%%%%%%%%%%%%%%%%%%%%%%%%%%%%%%%%%%%%%%%%%%%%%%%%%%%%%%%%%%%%%%%%%%%%%%%%%%%%%%%%%%%%%%%%%%%%%%%%%%%%%%%%%%%%%%%%%%%%%%%%%%%%%%%%%%%%%%%%%%%%%%%%%%%%%%%%%%%%%%%%%%%%%%%%%%%%%%%%%
\vskip0.12in
%%%%%%%%%%%%%%%%%%%%%%%%%%%%%%%%%%%%%%%%%%%%%%%%%%%%%%%%%%%%%%%%%%%%%%%%%%%%%%%%%%%%%%%%%%%%%%%%%%%%%%%%%%%%%%%%%%%%%%%%%%%%%%%%%%%%%%%%%%%%%%%%%%%%%%%%%%%%%%%%%%%%%%%%%%%%%%%%%%%%

At the full endpoint \(A(m_*)^a\ge D_p\mu\), hence \(s_n^a/a_n\ge c K_0\). The cell profile bound then gives \(v_*\gtrsim s_n y\) in each cell, so \(d_D^2\lambda pJ v_*^{-a}\lesssim a_n/s_n^a\le\eta_0\) there for \(K_0\) large. Together with the exterior bound this proves positive smooth slopes and a strict energy gap. The continuation criterion gives \(A(m_*)^a=D_p\mu=:T\). The bounded residual in \eqref{eq:flat-prefix-matched} on the full product gives \(s_n^a/a_n\lesssim1\), and \cref{prop:negative-local-energy-slopes} gives \(b_Q(V_*)>0\).

We next estimate the prefix slope at the same mass \(m_*\). Put \(x=x_{N+1}\) and \(e_N=1-J_N/J\), so \(0\le e_N\le1\). For \(|w|\ge2x\), the omitted logarithms are bounded in absolute sum by \(C\sum_{j>N}\rho_j/|w|\le C x\), and \(e_N\le Cx\) there. If \(V_*^D=\lambda pJv_*^{-a}\) and \(\dd\nu=V_*^Dk^2\dd w\), then \(\nu(D)<\infty\) by the same bounds as in \cref{lem:flat-collar-response}: \(V_*^D k^2\) is bounded in the cells and \(O(k^{1-p})\) elsewhere. Also \(\nu(B_{2x})\lesssim x^2\). Thus \(\|e_N\|_{L^2(\nu)}\lesssim x\). The difference \(W=v_{N,m_*}-v_*\ge0\) belongs to \(H^1_0(D)\) and satisfies
\[
 (-\Delta-C_N)W=e_N f_*,\qquad
 0\le C_N\le V_*^D,\qquad f_*=\lambda Jv_*^{-p}\le C V_*^D k .
\]
Testing with \(W\), coercivity and Cauchy--Schwarz give \(\|W\|_{H^1_0(D)}\lesssim x\). The source \(g=-\Delta W=f_*-f_{N,m_*}\) satisfies \(0\le g\le f_*\). Use the cutoff \(\chi_R=1\) on \(B_R(0)\), supported in \(B_{2R}(0)\). By \eqref{eq:flat-prefix-near-flux} and \(\|(1-\chi_R)P_Q\|_{H^1_0(D)}\lesssim R^{-1}\), the flux difference is bounded by
\[
 0\le A_N(m_*)-A(m_*)=\int_D P_Q g
     \lesssim R^{\beta/a}+x/R .
\]
The norm bound for the test follows from the reflected Poisson singularity \(P_Q=O(r^{-1})\), \(\nabla P_Q=O(r^{-2})\) at zero and its smooth zero trace elsewhere. Taking \(R=x^{a/2}\) gives \(O(x^{\beta/2})=o(a_N)\) by \(x\le x_N^5\).

In \eqref{eq:flat-prefix-matched} take \(n=N\), \(m=m_N\). The left side is bounded below, so \(A_N(m_N)^a\ge D_{p,N}\mu-Ca_N\). Moreover \(0\le D_p-D_{p,N}\lesssim\sum_{j>N}\rho_j=o(a_N)\) by the product estimate at zero. The same-mass slope estimate just proved and the branch order give \(m_*-m_N\lesssim a_N\).

%%%%%%%%%%%%%%%%%%%%%%%%%%%%%%%%%%%%%%%%%%%%%%%%%%%%%%%%%%%%%%%%%%%%%%%%%%%%%%%%%%%%%%%%%%%%%%%%%%%%%%%%%%%%%%%%%%%%%%%%%%%%%%%%%%%%%%%%%%%%%%%%%%%%%%%%%%%%%%%%%%%%%%%%%%%%%%%%%%%%
\vskip0.12in
%%%%%%%%%%%%%%%%%%%%%%%%%%%%%%%%%%%%%%%%%%%%%%%%%%%%%%%%%%%%%%%%%%%%%%%%%%%%%%%%%%%%%%%%%%%%%%%%%%%%%%%%%%%%%%%%%%%%%%%%%%%%%%%%%%%%%%%%%%%%%%%%%%%%%%%%%%%%%%%%%%%%%%%%%%%%%%%%%%%%

For the other inequality, consider masses \(m\in[m_N,m_*]\) with \(A_N(m)^a>T-a_N\). Since \(D_p\ge D_{p,N}\) and \(a_n\ge a_N\) on all retained cells, \eqref{eq:flat-prefix-matched} gives \(s_{N,m,n}^a\ge cK_0a_n\) there. The same cell and exterior bounds give \(d_D^2\lambda pJ_N v_{N,m}^{-a}\le\eta_0\) and \(v_{N,m}\gtrsim k\) off the cells. Each fixed prefix in this window has a positive Martin lower bound also on its finitely many retained cells. At such a solution the lower bound by \(k\) and the finite inverse-power row give a finite Martin row for the linearization. The strict form gap bounds its Martin spectral radius below one by \cref{thm:negative-two-space-radius}, so the implicit function theorem permits continuation below its mass. In particular \(m_N\) cannot be in the window, so \(A_N(m_N)^a\le T-a_N\). The branch slope \(A_N(m)\) is continuous on \([m_N,m_*]\). At \(m_N\) this follows by right monotone convergence in the flux with the stable majorant from the prefix proof; at interior masses it follows from smooth branch dependence in the Martin norm and the same flux domination. Since \(A_N(m_*)^a\ge T\), some \(m'\) in the window has \(A_N(m')^a=T-a_N/2\). Apply \cref{lem:flat-collar-response} with \(m_1=m'\), \(m_2=m_*\). Their slope difference is at least \(ca_N\) and at most \(C(m_*-m')\). This proves \(m_*-m_N\gtrsim a_N\).
\end{proof}

\section*{Declaration on the use of artificial intelligence}
Generative artificial intelligence assisted with literature searches and reference checks, the development and review of mathematical arguments, language editing, and \LaTeX\ typesetting. The authors take full responsibility for the manuscript's accuracy and integrity, including its mathematical statements and proofs.

\end{document}